\documentclass[11pt,a4paper]{article}
\usepackage{amsfonts}
\usepackage{amsmath,amssymb}
\usepackage{mathrsfs}
\usepackage{hyperref}
\usepackage{graphicx}
\usepackage{lineno}
\modulolinenumbers[5]

\newtheorem{theorem}{Theorem}[section]

\newtheorem{lemma}[theorem]{Lemma}
\newtheorem{proposition}{Proposition}

\newtheorem{definition}[theorem]{Definition}
\newtheorem{remark}{Remark}

\numberwithin{equation}{section}\allowdisplaybreaks

\begin{document}

\title{\large\bf  Resonant Decompositions and Global Well-posedness for 2D Zakharov-Kuznetsov Equation in Sobolev spaces of Negative Indices}

\author{\normalsize \bf Minjie  Shan$^{\dag}$,  \ Baoxiang Wang$^{\ddag}$   and Liqun Zhang$^{\dag}$ \\
\footnotesize
\it $^\dag$Academy of Mathematics and Systems Science, Chinese Academy of Sciences, Beijing 100190, P.R. China, \\ \footnotesize
\it E-mail: smj@amss.ac.cn, lqzhang@math.ac.cn\\
\footnotesize
\it $^\ddag$LMAM, School of Mathematical Sciences, Peking University, Beijing 100871, China, \\
\footnotesize
\it Emails: wbx@math.pku.edu.cn\\
} \maketitle

\thispagestyle{empty}
\begin{abstract}
The Cauchy problem for Zakharov-Kuznetsov equation on $\mathbb{R}^2$ is shown to be global
well-posed for the initial date in $H^{s}$ provided $s>-\frac{1}{13}$. As conservation
laws are invalid in Sobolev spaces below $L^2$, we construct an almost conserved quantity
using multilinear correction term following the $I$-method  introduced by Colliander,
Keel, Staffilani, Takaoka and Tao. In contrast to KdV equation, the main difficulty is to
handle the resonant interactions which are significant due to the multidimensional and
multilinear setting of the problem. The proof relies upon the bilinear Strichartz estimate
and the nonlinear Loomis-Whitney inequality.\\

 \\

\end{abstract}
\section{Introduction}

We consider the initial value problem for the Zakharov-Kuznetsov equation (ZK)
\begin{equation}
  \left\{
   \begin{aligned}
   &\partial_{t}u +( \partial^3_{x}+\partial^3_{y}) u +(\partial_{x}+\partial_{y}) u^2 = 0,
\ \ \ (x,y)\in\mathbb{R}^2, \ t\geq0     \quad \\
   &u(0,x,y)=u_0(x,y)\in H^{s}(\mathbb{R}^2), \label{ZK} \\
   \end{aligned}
   \right.
  \end{equation}
where  $u=u(t,x,y)$ is a real-valued function.

The Cauchy problem associated with dispersive equations with derivative nonlinearity has been extensively studied for about forty years, for instance see \cite{Kato1983},\cite{Abdelouhabet1989} and \cite{Kenigetal1991},\cite{KPV93},\cite{Ponce1991}.  Kenig, Ponce and Vega first proved new dispersive estimates which enable them to deal with low regularity problem. By using the fixed point argument, they obtained local well-posedness for dispersive equations. Then, Bourgain  \cite{Bourgain1993a} \cite{Bourgain1993b} introduced the now so-called Bourgain spaces to solve a wide class of dispersive equations with very rough initial data. Since the nonlinearity is in general algebraic, the fixed point argument ensures the real analyticity of the solution map. However, Molinet, Saut and Tzvetkov \cite{MST2001} found that a large class of ¡°weakly¡± dispersive equations, including in particular the Benjamin-Ono equation, cannot be solved by a fixed point argument for initial data in any Sobolev spaces $H^s$. This is caused by bad interactions between high frequencies and very low frequencies. In \cite{Mv2015}, Molinet and Vento proposed a  new approach to get local and global well-posedness results for dispersive equations without too-strong resonances. This approach combines classical energy estimates with Bourgain-type estimates on a time interval that does not depend on the space frequency. These approaches have been developed to lower the regularity requirement on initial data.

 The ZK equation was introduced by Zakharov and  Kuznetsov in \cite{ZK74} as a  model to
describe the propagation of ionic-acoustic waves in magnetized plasma. It is a
multidimensional generalization of  the Korteweg-de Vries equation. In  \cite{EWLKHS82},
Laedke and Spatschek derived  the two-dimensional Zakharov-Kuznetsov equation from the basic
hydrodynamic equations.  Lannes, Linares and Saut \cite{LLS13} showed that in two and three
dimensions the ZK equation is a long-wave limit of the Euler-Poisson system.

The ZK equation is not completely integrable, but there are two important invariants,
\begin{align}
M(u)(t)= \int_{\mathbb{R}^2} u^2(t,x,y) dxdy = \int_{\mathbb{R}^2} u^2_0(x,y) dxdy = M(u_0)
\label{M()}
\end{align}
and
\begin{align}
E(u)(t)= \int_{\mathbb{R}^2} \frac{1}{2} | \nabla u|^2 -\frac{1}{2}
\partial_{x}u\partial_{y}u - \frac{1}{3} u^3 dxdy = E(u_0).\label{E()}
\end{align}

Both the local and global-in-time initial value problems for Zakharov-Kuznetsov equations
have attracted a substantial literature  \cite{GH14}, \cite{Kinoshita19}, \cite{LP09},
\cite{MP15},  \cite{Shan1},\cite{Shan2}. In two dimensional case, Faminskii \cite{F95} first
showed the local well-posedness for the two-dimensional ZK equation in energy space
$H^1(\mathbb{R}^2)$ by using the local smoothing effect and  a maximal function estimate for
the linearized equation. This idea is originally from Kenig, Ponce and Vega \cite{KPV93}.
The local solution can be extended to a global one via the $L^2$ and $H^1$ conserved
quantities. In \cite{LP09}, Linares and Pastor proved the local well-posedness for $ s >
\frac{3}{4}$. Applying the Fourier restriction norm method and one kind of sharp Strichartz
estimates, Gr\"unrock and Herr \cite{GH14}, and Molinet and Pilod \cite{MP15} proved
independently that the Cauchy problem for \eqref{ZK} is local well-posed in
$H^s(\mathbb{R}^2)$ for  $ s > \frac{1}{2}$ . In  \cite{Shan1} and \cite{Shan2}, we showed
the local well-posedness in $B^{ \frac{1}{2}}_{2,1}(\mathbb{R}^2)$ and the global
well-posedness in $H^s(\mathbb{R}^2)$ for $ s > \frac{5}{7}$ using atomic spaces  introduced
by Koch and Tataru. Inspired by the ideas of  Bejenaru et al (\cite{BHHT09} \cite{BH11}
\cite{BHT10}), Kinoshita recently obtained the local well-posedness in $H^s(\mathbb{R}^2)$
for $ s > -\frac{1}{4}$ in \cite{Kinoshita19}.
 Meanwhile, he proved that the date-to-solution map fails to be $C^2$ from the unit ball in
$H^s(\mathbb{R}^2)$ to $C([0,T];H^s)$, which suggests that for the Picard iteration approach
the subcritical threshold $ -\frac{1}{4}$ is optimal. The tools utilized to greatly improve
local well-posedness results are the nonlinear version of the classical Loomis-Whitney
inequality and an almost orthogonal decomposition technique of resonant frequencies.
Certainly,  the global well-posedness in $L^2(\mathbb{R}^2)$ is a natural result because of
the conservation of mass.

We now turn our attention to the global-in-time well-posedness problem with  low regularity.
It was Bourgain who first made a breakthrough in improving well-posed results below energy
space in  \cite{Bourgain98} and \cite{Bourgain99} where he obtained global well-posedness
for  the nonlinear Schr\"odinger equation in  $H^s(\mathbb{R}^2)$ when $ s > \frac{3}{5}$
using what is now referred to as the Fourier restriction norm method. Taking use of the
$I$-method developed in  \cite{CKSTT01}, \cite{CKSTT02} and  \cite{CKSTT03}, Colliander et
al established an a-priori bound of the solution to show the global well-posedness for the
cubic nonlinear Schr\"odinger equation in $H^s$  for  $ s > \frac{4}{7}$ . The increment
estimate of the solution in $H^s$ implies the global well-posedness via the local
well-posedness theory and standard limiting arguments. Subsequently, in \cite{CKSTT08} they
refined the global well-posed result from  $ \frac{4}{7}$ to  $ \frac{1}{2}$ by adding a
 "correction term" in order to damp out some oscillations at the expense of causing a
singular symbol which is intractable to estimate. To get around this new difficulty they
employed a resonant decomposition technique appeared previously in \cite{Bourgain04} and
\cite{BDGS07}.

When comes to the global well-posedness for the ZK equation \eqref{ZK} below
$L^2(\mathbb{R}^2)$, we encounter similar difficulty. Due to the multidimensional setting
(in contrast to the KdV equation \cite{CKSTT03}), adding a  "correction term" to the
modified mass functional $E^0(Iu)$ will bring a singular symbol
$|\xi_1\xi_2\xi_3+\eta_1\eta_2\eta_3|$. It is easy to see that even though the frequencies
$|(\xi_j,\eta_j)|\sim N_j$ for $j=1,2,3$ are high frequencies,
$|\xi_1\xi_2\xi_3+\eta_1\eta_2\eta_3|$ can be very small. Not the same as the nonlinear
Schr\"odinger equation case, we directly add an "correction term" which contains
non-resonant interactions. The reason is that for the cubic nonlinear Schr\"odinger equation
one needs to modify the energy functional, but for the ZK equation what we need to modify is
the mass functional. As a consequence,  non-resonant interactions shall make a  cancellation
with parts of oscillations in the functional, then resonant interactions will appear in the
trilinear remainder term.

The purpose of this paper is to obtain the global well-posedness in negative Sobolev spaces.
We now state the main result.
\begin{theorem} \label{gwp}
The initial value problem \eqref{ZK} is globally well-posed in $H^s(\mathbb{R}^2)$ for
$ s >- \frac{1}{13}$.
\end{theorem}
\begin{remark} "Globally well-posed" means that given data  $u_0 \in H^s(\mathbb{R}^2)$ and
any time $T>0$, there exists a unique solution to \eqref{ZK}
\end{remark}
$$ u(t,x,y)\in \mathbb{C}\big([0,T]; H^s(\mathbb{R}^2)\big)\cap  X^{s,\frac{1}{2}+}_T$$
which depends continuously upon  $u_0$.

As in other application of the $I$-method, the matter is reduced to the construction of a
modified mass (or energy) functional $E^0(Iu)$ and making sure that it has almost
conservation properties. The appearance of the singular symbol
$|\xi_1\xi_2\xi_3+\eta_1\eta_2\eta_3|$ impels us to truncate the correction term to
non-resonant interactions. But the key point is how to deal with resonant interactions. The
estimates for the resonant interactions will eventually determine the almost conservation
properties of the modified mass functional. Our strategies of controlling the trilinear
remainder term containing resonant interactions are the bilinear Strichartz estimate and the
nonlinear Loomis-Whitney inequality. More specifically, we employ the bilinear Strichartz
estimate for high modulation and use a convolution estimate on hypersurfaces which was
introduced by Bennet, Carbery and Wright  \cite{BCW05} for low modulation.  This convolution
estimate originates from the classical Loomis-Whitney inequality \cite{LW49}. Bejenare, Herr
and Tataru generalized it in  \cite{BHT10} and used the nonlinear Loomis-Whitney inequality
to deal with the initial value problem for the Zakharov system in the space of optimal
regularity (see  \cite{BHHT09} and \cite{BH11}). It's also worth mentioning that the
Loomis-Whitney inequality is related to the multilinear restriction theorem  \cite{BCW05}.

The upper bound given by the nonlinear Loomis-Whitney inequality relies on the
transversality of three characteristic hypersurfaces. To be precise, the uniform upper bound
is inversely proportional to the square root of the determinant consisting of unit normal
vectors of these three characteristic hypersurfaces. For the Zakaharov system, the
transversality is related to the size of the angle between two frequencies of the wave.
While the transversality for the Zakharov-Kuznetsov is complex. In fact,  the transversality
depends not only on the angle but also the sizes of frequencies. In \cite{Kinoshita19},
Kinoshita used a new almost orthogonal decomposition to estimate the transversality. But
this decomposition technique is very complicated.

Even if the transversality depends on both the angle and the sizes of frequencies, we find
that these two bad factor can not happened at the same time. In other words, one can use the
angular frequency decomposition to estimate the transversality when the frequencies is
controllable, and use the Whitney type decomposition which actually contains certain almost
orthogonality when the angle is controllable. Moreover, due to the restriction on
$|\xi_1\xi_2\xi_3+\eta_1\eta_2\eta_3|$, the remainder terms are easier to estimate in
resonant interactions region.

\textbf{Organization of the paper.} In Section 2, we introduce function spaces and some
estimates which will be used in what follows. In Section 3 we recall the I-method and apply
it to investigate the increment of modified mass. This leads to the main theorem by
contraction mapping principle and iteration argument. Then in Section 4, we show Proposition
\ref{fixedtimep1} that we call the fixed time estimate. Section 5 is devoted to the proof
the quadrilinear estimate. Finally, in Section 6 we prove the most important proposition to
obtain the global well-posedness  in sobolev spaces below $L^2$. This section contains
high-high interactions and high-low interactions which can be divided into non-parallel and
parallel subcases respectively.

We give the notation that will be used throughout this paper. Given $A,B,C \geq 0$, $A\lesssim B$ stands for $A\leq C \cdot B$. $A\sim B$ means that $A\lesssim B$ and $B\lesssim A$, while $A\gg B$ means  $A>C \cdot B$.  We use the notation $c+\equiv
c+\epsilon$ for some  $0<\epsilon\ll 1$ and write  $c++\equiv c+2\epsilon$ and $ c-\equiv
c-\epsilon$. $\chi\in C^\infty_0([-2,2])$ is a  fixed smooth cut-off function satisfying
$\chi$ is even, nonnegative,  and $\chi=1$ on $[-1,1]$. Denote spatial variables by $x,y$
and their dual Fourier variables by $\xi,\eta$. $\tau$ is the dual variable of the time $t$.
$\mathscr F(u)$ or $\hat u$ will denote space-time Fourier transform of $u$, whereas
$\mathscr{F}_{x,y}(u)$ or ${\widehat u}^{xy}$ will denote its Fourier transform in space.
For brevity, we write $\zeta=(\xi,\eta)$ and $\lambda=(\tau,\xi,\eta)$. We also use $N,M$ to
denote dyadic numbers and  write $\sum_{N\geq 1}a_N:=\sum_{n\in \mathbb{N}}a_{2^n}$ ,
$\sum_{N\geq M}a_N:=\sum_{n\in \mathbb{N};2^n\geq M}a_{2^n}$ for dyadic summations.
\section{Function spaces and some estimates}

Denote $\psi(x) :=\chi(x)-\chi(2x)$ and $\psi_N :=\psi(N^{-1}\cdot)$. The Littlewood-Paley
projections for frequency and modulation are defined by
\begin{align}
 \widehat{P_0 u}=\chi(|\zeta|)\widehat{u}  \hspace{5mm} &and \hspace{3mm}
 \ \ \widehat{P_Nu}=\psi_N(|\zeta|)\widehat{u}  for   N\geq 1 \notag \\
 \widehat{Q_0 u}=\chi(|\tau -\xi^3-\eta^3|)\widehat{u}  \hspace{5mm} &and \hspace{3mm}
 \ \ \widehat{Q_Lu}=\psi_N(|\tau -\xi^3-\eta^3|)\widehat{u} for   L\geq 1 ,\nonumber
\end{align}
we often write $u_N=P_Nu$, $u_{N,L}=P_NQ_Lu$ and denote the space-time Fourier support of
$P_NQ_L$ by
\begin{align}
 G_{N,L}=\left\{(\tau ,\xi,\eta)\in\mathbb{R}^3\Big|\psi_N(|(\xi,\eta)|)\psi_N(|\tau
-\xi^3-\eta^3|)\neq 0\right\}\nonumber
\end{align}

Given Bourgain spaces exponents $s,b\in\mathbb{R}$ and a function $u(x,y,t)\in\mathcal
S'(\mathbb{R}^3)$, we use the dyadic frequency localization operators $P_N$ and $Q_L$ to
define
\begin{align}
\|u\|_{X^{s,b}}=\left(\sum_{N,L}N^{2s}L^{2b}\|P_NQ_Lu\|^{2}_{L^{2}_{txy}}\right)^{\frac{1}{2}}.
\label{X Norm}
\end{align}

The truncated versions of  $X^{s,b}$ norm is defined as
\begin{align}
\|u\|_{X^{s,b}_\delta}=\inf \limits_{{\tilde u=u} \ on \ [0,\delta]} \|\tilde u\|_{X^{s,b}}
\label{X1 Norm}.
\end{align}
We also define equidistant partitions of unity on $\mathbb{R}$ and on the unit circle,
\begin{align}1=\sum_{j\in\mathbb{Z}}\beta_j, \hspace{5mm}
\beta_j(t)=\chi(t-j)\left(\sum_{k\in\mathbb{Z}}\chi(t-k)\right)^{-1},\nonumber \end{align}
and for dyadic number $A\in \mathbb{N}$,
\begin{align}1=\sum^{A-1}_{j=0}\beta^A_j, \hspace{5mm}
\beta^A_j(\theta)=\beta_j(\frac{A\theta}{\pi})+\beta_{j-A}(\frac{A\theta}{\pi}).\nonumber
\end{align}

The support sets of $ \beta^A_j$ are
\begin{align}
\Theta^A_j:=\left[\frac{\pi}{A}(j-2),\frac{\pi}{A}(j+2)\right]\cup\left[-\pi+\frac{\pi}
{A}(j-2),-\pi+\frac{\pi}{A}(j+2)\right].\nonumber
\end{align}
Now we use $ \beta^A_j$ to define the angular frequency localization operators $R^A_j$,
\begin{align}
\mathscr F(R^A_ju)(\xi,\eta)=\beta^A_j(\theta)\widehat{u}(\xi,\eta), \nonumber
\end{align}
where $(\xi,\eta)=|(\xi,\eta)|(\cos\theta, \sin\theta) $. The support sets of these
frequency localization functions $R^A_ju $ are $ \tilde {\mathcal D}^A_j=\mathbb{R}\times
{\mathcal D}^A_j$, where
\begin{align}
 {\mathcal D}^A_j=\left\{(\xi,\eta)\in\mathbb{R}^2\Big|(\xi,\eta)=|(\xi,\eta)|(\cos\theta,
\sin\theta), \theta \in\Theta^A_j \right\}.\label{angulardecomp}
\end{align}
Now one can decompose $u$ as
$$u=\sum^{A-1}_{j=0}R^A_ju.$$

Let's recall some well-known estimates about  $X^{s,b}$ spaces used in the well-posedness
theory and the Strichartz estimates associated to the unitary group
$e^{-t( \partial^3_{x}+\partial^3_{y})}$.
\begin{lemma}\label{X Lema1}
Let $s>-\frac{1}{4}$ . Assume that $u_0\in L^2(\mathbb{R}^2)$ and $F\in X^{0,
-\frac{1}{2}++}_\delta$ , then we have
\begin{align}
\|e^{-t(\partial^3_{x}+\partial^3_{y})}u_0\|_{X^{0, \frac{1}{2}+}_\delta} \lesssim
\|u_0\|_{L^2(\mathbb{R}^2)},\label{estimate1}
\end{align}
\begin{align}
\|\int^t_0 e^{-(t-t')(\partial^3_{x}+\partial^3_{y})}F(t')dt'\|_{X^{0, \frac{1}{2}+}_\delta}
\lesssim \|F\|_{X^{0, -\frac{1}{2}+}_\delta},\label{estimate2}
\end{align}
\begin{align}
\|F\|_{X^{0, -\frac{1}{2}+}_\delta} \lesssim \delta^{0+} \|F\|_{X^{0,
-\frac{1}{2}++}_\delta},\label{estimate3}
\end{align}
\begin{align}
\|(\partial_{x}+\partial_{y})(uv)\|_{X^{s, -\frac{1}{2}++}} \lesssim  \|u\|_{X^{s,
\frac{1}{2}+}}\|v\|_{X^{s, \frac{1}{2}+}}.  \label{estimate4}
\end{align}
\end{lemma}
{\bf Proof.}  \eqref{estimate1}, \eqref{estimate2} and \eqref{estimate3} can be found in
\cite{GTV97} and \cite{MP15}. The last bilinear estimate \eqref{estimate4} is gave in
\cite{Kinoshita19} .
\begin{lemma} \label{X Lema2}
Let  $s,b\in\mathbb{R}$ and $u\in X^{s,b}$. Then we have
\begin{align}
\||\nabla_x|^{\frac{1}{8}}|\nabla_y|^{\frac{1}{8}}Q_Lu\|_{L^4(\mathbb{R}^3)} \lesssim
L^{\frac{1}{2}}\|Q_Lu\|_{L^2(\mathbb{R}^3)}, \label{estimate15}
\end{align}
\begin{align}
\|Q_Lu\|_{L^4(\mathbb{R}^3)} \lesssim L^{\frac{5}{12}}\|Q_Lu\|_{L^2(\mathbb{R}^3)}.
\label{estimate16}
\end{align}
\end{lemma}
{\bf Proof.} From \cite{KPV91},
\begin{align}
\||\nabla_x|^{\frac{1}{2p}}|\nabla_y|^{\frac{1}{2p}}e^{-t(\partial^3_{x}+\partial^3_{y})}u_0\|
_{L^p_tL^q_{xy}} \lesssim \|u_0\|_{L^2_{xy}}, \ \ \ if \ \frac{2}{p}+\frac{2}{q}=1, \ p>2,
\nonumber
\end{align}
and
\begin{align}
\ \ \ \  \|e^{-t(\partial^3_{x}+\partial^3_{y})}u_0\|_{L^p_t L^q_{xy}} \lesssim
\|u_0\|_{L^2_{xy}},\ \ \ \ \ \ \ \ if \ \frac{3}{p}+\frac{2}{q}=1, \ p>3.  \nonumber
\end{align}

By the extension principle for $X^{s,b}$ (see Lemma 5.3 in \cite{WaHaHuGu11}), we have
\begin{align}
\||\nabla_x|^{\frac{1}{2p}}|\nabla_y|^{\frac{1}{2p}}u\|_{L^p_tL^q_{xy}} \lesssim
\|u\|_{X^{0, \frac{1}{2}+}}, \ \ \ if \ \frac{2}{p}+\frac{2}{q}=1, \ p>2,  \nonumber
\end{align}
and
\begin{align}
\ \ \ \  \|u\|_{L^p_t L^q_{xy}} \lesssim \|u\|_{X^{0, \frac{1}{2}+}},\ \ \ \ \ \ \ \ if \
\frac{3}{p}+\frac{2}{q}=1, \ p>3.  \nonumber
\end{align}
Hence,
\begin{align}
\||\nabla_x|^{\frac{1}{2p}}|\nabla_y|^{\frac{1}{2p}}Q_Lu\|_{L^p_tL^q_{xy}} \lesssim
L^{\frac{1}{2}}\|Q_Lu\|_{L^2_{txy}}, \ \ \ if \ \frac{2}{p}+\frac{2}{q}=1, \ p>2,
\label{estimate8a1a}
\end{align}
and
\begin{align}
\ \ \ \  \|Q_Lu\|_{L^p_t L^q_{xy}} \lesssim L^{\frac{1}{2}}\|Q_Lu\|_{L^2_{txy}},\ \ \ \ \ \
\ \ if \ \frac{3}{p}+\frac{2}{q}=1, \ p>3. \label{estimate8b1b}
\end{align}

We get \eqref{estimate15} by choosing $p=q=4$ in \eqref{estimate8a1a}. The estimate
\eqref{estimate16} follows from interpolation between \eqref{estimate8b1b} and the trivial
bound $\|Q_Lu\|_{L^2(\mathbb{R}^3)}=\|Q_Lu\|_{L^2(\mathbb{R}^3)}$(see  \cite{Kinoshita19}).

For the Schr\"odinger equation, Bourgain showed the bilinear generalization of the linear
$L^4$ Strichartz estimate ( see \cite{Bourgain98} Lemma 111 ). A similar bilinear Strichartz
estimate holds true for the Zakharov-Kuznetsov interaction (see \cite{Kinoshita19} and
\cite{BHHT09}).
\begin{lemma} \label{X Lema3}
Let $C>0$, $N_1 \sim N_2\gtrsim N_3\geq1$ and $L_1, L_2, L_3\geq1$. Assume that $v_j\in
L^2(\mathbb{R}^3)$ and $supp\  \mathscr F(v_j) \subset G_{N_j,L_j}$ for $j=1,2$ and $3$. If
$\lambda_1+\lambda_2+\lambda_3=0$, $|\xi_1|\sim N_1$, $|\xi_2|\sim N_2$,  and
$\Big|\left\{\eta_j\big|\lambda_j=(\tau_j,\xi_j,\eta_j)\in \text{supp}\  \mathscr F(v_j)
\right\}\Big|\lesssim C$ for $j=1,2$, then
\begin{align}
 \big\|\chi_{G_{N_3,L_3}}(\lambda_3)\int \widehat{v}_1(\lambda_1)
\widehat{v}_2(\lambda_1+\lambda_3)d\lambda_1\big\|_{L^2_{\lambda_3}}\lesssim
C^{\frac{1}{2}}N_1^{-1}(L_1L_2)^{\frac{1}{2}}\|v_1\|_{L^{2}}\|v_2\|_{L^{2}},\label{bles1}
\end{align}
\begin{align}
 \big\|\chi_{G_{N_2,L_2}}(\lambda_2)\int \widehat{v}_1(\lambda_1)
\widehat{v}_3(\lambda_1+\lambda_2)d\lambda_1\big\|_{L^2_{\lambda_2}}\lesssim
C^{\frac{1}{2}}N_1^{-1}(L_1L_3)^{\frac{1}{2}}\|v_1\|_{L^{2}}\|v_3\|_{L^{2}},\label{bles2}
\end{align}
\begin{align}
 \big\|\chi_{G_{N_1,L_1}}(\lambda_1)\int \widehat{v}_2(\lambda_2)
\widehat{v}_3(\lambda_1+\lambda_2)d\lambda_2\big\|_{L^2_{\lambda_1}}\lesssim
C^{\frac{1}{2}}N_1^{-1}(L_2L_3)^{\frac{1}{2}}\|v_2\|_{L^{2}}\|v_3\|_{L^{2}},\label{bles3}
\end{align}
and
\begin{align}
 \|v_1v_3\|_{L^2([0,\delta]\times\mathbb{R}^2)}\lesssim
\frac{N_3^{\frac{1}{2}}}{N_1}\|v_1\|_{X^{0, \frac{1}{2}+}_\delta}\|v_3\|_{X^{0,
\frac{1}{2}+}_\delta},\label{bles4}
\end{align}
where $\chi_{G_{N_j,L_j}}$ are characteristic functions on  $G_{N_j,L_j}$ for $j=1,2,3$.
\end{lemma}
{\bf Proof.} From the Cauchy-Schwartz inequality,
\begin{align}
& \|\chi_{G_{N_3,L_3}}(\lambda_3)\int \widehat{v}_1(\lambda_1)
\widehat{v}_2(\lambda_1+\lambda_3)d\lambda_1\|_{L^2_{\lambda_3}}\notag \\
\lesssim &\sup_{\lambda_3\in
G_{N_3,L_3}}\big|E(\lambda_3)\big|^{\frac{1}{2}}\|v_1\|_{L^{2}}\|v_2\|_{L^{2}},\label{bles1a}
\end{align}
where
\begin{align}
E(\lambda_3)=\left\{(\tau_1,\xi_1,\eta_1)\in \text{supp}\ \widehat{v}_1     \big|
(\tau_1+\tau_3,\xi_1+\xi_3,\eta_1+\eta_3)\in \text{supp}\  \widehat{v}_2 \right\}.\nonumber
\end{align}

For any fixed $\lambda_3=(\tau_3,\xi_3,\eta_3)$, we have
\begin{align}
|\tau_1-\xi^3_1-\eta^3_1|\sim L_1,\nonumber
\end{align}
and
\begin{align}
 |\tau_1+\tau_3-(\xi_1+\xi_3)^3-(\eta_1+\eta_3)^3|\sim L_2 .\nonumber
\end{align}
Hence,
\begin{align}
\Big| \left\{\tau_1\big|  \tau_1\in E(\lambda_3)  \right\}\Big|\lesssim
\min\{L_1,L_2\}.\nonumber
\end{align}

Since
\begin{align}
|\partial_{\xi_1}(\tau_1-\xi^3_1-\eta^3_1)|=3|\xi_1|^2\sim N_1^2,\nonumber
\end{align}
Therefore,
$$\Big| \left\{\xi_1\big|  \xi_1\in E(\lambda_3)  \right\}\Big|\lesssim
N_1^{-2}\max\{L_1,L_2\}, \ \ \Big| \left\{\eta_1\big|  \eta_1\in E(\lambda_3)
\right\}\Big|\lesssim C.$$
This implies
\begin{align}
\sup_{\lambda_3\in G_{N_3,L_3}}\big|E(\lambda_3)\big|\lesssim CN_1^{-2}L_1L_2.\label{bles1b}
\end{align}
Thus, \eqref{bles1} follows from \eqref{bles1a} and \eqref{bles1b}. \eqref{bles2} --
\eqref{bles4} can be dealt with in a similar way.

Next we recall the nonlinear version of the classical Loomis-Whitney inequality which play a
crucial role in our analysis on low modulation.
\begin{lemma}[ see \cite{BHT10}]\label{LW1}
Let $j=1,2,3$. Assume that
\begin{itemize}
 \item[\rm (i)] the oriented surface $S_j^*$ is given as
$$S_j^*=\left\{ \lambda_j \in U_j \ \big|\  \Phi_j(\lambda_j)=0, \ \nabla \Phi_j\neq0,\
\Phi_j\in C^{1,1}(U_j) \right\}$$
for a convex $U_j\in\mathbb{R}^3$ such that $dist(S_j,U^c_j)\geq diam(S_j)$, where $S_j$
is an open and bounded subset of  $S_j^*$;

\item[\rm (ii)] the unit normal vector field $\mathfrak {n}_j$ on $S_j^*$ satisfies the
    H\"older condition
\begin{align}
\sup_{\lambda,\tilde{\lambda}\in S_j^*}\frac{|\mathfrak {n}_j(\lambda)-\mathfrak
{n}_j(\tilde{\lambda})|}{|\lambda-\tilde{\lambda}|}+\frac{|\mathfrak
{n}_j(\lambda)(\lambda-\tilde{\lambda})|}{|\lambda-\tilde{\lambda}|^2}\lesssim 1;\nonumber
\end{align}
\item[\rm (iii)] there exists a positive constant $d>0$ such that the matrix $\mathcal
    {N}(\lambda_1,\lambda_2,\lambda_3)=\big(\mathfrak {n}_1(\lambda_1),\mathfrak
    {n}_2(\lambda_2),\mathfrak {n}_3(\lambda_3)\big)$ satisfies the transversal condition
\begin{align}
d\leq |det\mathcal {N}(\lambda_1,\lambda_2,\lambda_3)|\leq 1 \nonumber
\end{align} for all $(\lambda_1,\lambda_2,\lambda_3)\in S_1^*\times S_2^*\times S_3^*$.
\end{itemize}
We also assume $diam(S_j)\lesssim d$, then for each $f\in L^2(S_1)$ and $g\in L^2(S_2)$ the
restriction of the convolution $f*g$ to $S_3$ is a well-defined $L^2(S_3)$-function which
satisfies
\begin{align}
\|f*g\|_{L^2(S_3)}\lesssim\frac{1}{\sqrt{d}}\|f\|_{L^2(S_1)}\|g\|_{L^2(S_2)}.\label{LW1a}
\end{align}
\end{lemma}                                                                                       \section{Almost conservation law}

Let's recall the $I$-method. It's convenient to introduce some notation for multilinear
operators involving $u$ before our construction of the modified mass functional $E^0$.

Given $M:\mathbb{R}^{2k} \to \mathbb{C}$, we say $M$ is symmetric if
$$M(\zeta_1, \cdots,\zeta_k)=M(\sigma(\zeta_1, \cdots,\zeta_k))$$
for all $\sigma \in S_k$, where $S_k$ is the permutation group for $k$ elements. Define the
symmetrization of $M$ as following
$$[M]_{sym}(\zeta_1, \cdots,\zeta_k)=\frac{1}{k!}\sum_{\sigma \in S_k}M(\sigma(\zeta_1,
\cdots,\zeta_k)).$$

For each $M$, a $k$-linear functional acting on $k$ functions $u_1, \cdots,u_k$ is given by
$$\Lambda_k(M;u_1, \cdots,u_k)=\int_{\zeta_1+\cdots+\zeta_k=0}M(\zeta_1,
\cdots,\zeta_k){\widehat u_1}^{xy}(\zeta_1)\cdots{\widehat u_k}^{xy}(\zeta_k).$$
We abbreviate $\Lambda_k(M):=\Lambda_k(M;u, \cdots,u)$. Note that there is symmetrization
rule for  $\Lambda_k(M)$:
$$\Lambda_k(M)= \Lambda_k([M]_{sym}).$$

Now we investigate the behaviour of these multilinear forms $\Lambda_k(M)$ in time.

From the identity
$$\partial_t{\widehat u}^{xy}(\zeta)=i(\xi^3+\eta^3){\widehat
u}^{xy}(\zeta)-i(\xi+\eta){\widehat {u^2}}^{xy}$$
arising from \eqref{ZK}, together with some Fourier analysis, one can easily get the
following lemma.
\begin{lemma}\label{Lambda_k}
Suppose that  $M_k$ is symmetric and independent of time. If $u$ is a solution to
\eqref{ZK}, then we have the differentiation formula
\begin{align}
\frac{d\Lambda_k(M_k)}{dt}=\Lambda_k(M_kh_k)-ik\Lambda_{k+1}\big(X(M_k)\big),
\label{Lambda_k1}
\end{align}
where
\begin{align}
h_k=i(\xi_1^3+\eta_1^3+\cdots+\xi_k^3+\eta_k^3), \nonumber
\end{align}
and
\begin{align}
X(M_k)=M_k(\zeta_1,
\cdots,\zeta_{k-1},\zeta_k+\zeta_{k+1})\cdot(\xi_k+\xi_{k+1}+\eta_{k}+\eta_{k+1}). \nonumber
\end{align}
\end{lemma}

Given $s<0$ and $N\gg1$, $m^s_N(\zeta)$ is a smooth, radially symmetric, non-increasing
function satisfying
$$ m^s_N(\zeta)=
\left\{
\begin{array}{l} 1 \hspace{9mm} |\zeta|\leq N \\
 (\frac{|\zeta|}{N})^{s} \hspace{3mm} |\zeta|\geq 2N
\end{array}
\right. .
$$
Define the Fourier multiplier operator
$$\widehat{I^s_Nf}^{xy}(\zeta):=m^s_N(\zeta) \widehat f^{xy}(\zeta).$$
Sometimes we drop the $N$ and $s$ by writing $I$ and $m$ for simplicity.

The next lemma is useful in low regularity global well-posedness theory.
\begin{lemma}[Lemma 12.1 in \cite{CKSTT04})]\label{L5}
Let $n$ be positive integer. Suppose that $Z$, $X_1,\cdots, X_n$ are translation invariant
Banach spaces and $T$ is a translation invariant $n$-linear operator such that
$$\|I^s_1T(u_1, \cdots,u_n)\|_Z\lesssim\prod_{j=1}^n\|I^s_1u_j\|_{X_j}$$
for all $u_1, \cdots,u_n$ and all $-\frac{1}{2}\leq s \leq 0$. Then we have
$$\|I^s_NT(u_1, \cdots,u_n)\|_Z\lesssim\prod_{j=1}^n\|I^s_Nu_j\|_{X_j}$$
for  all $u_1, \cdots,u_n$ , all $-\frac{1}{2}\leq s \leq 0$ and $N\geq 1$, with the
implicit constant independent of $N$.
\end{lemma}

Note that $X^{s,b}$ is translation invariant Banach space and $\partial_x+\partial_y$ is a translation
invariant multi-linear operator.

Let's consider the modified energy functional. Using the $I$-operator and the Fourier
inversion formula, we observe that
\begin{align}
E^0(u)=M(Iu)=\|Iu\|^2_2=\Lambda_2\big(m(\zeta_1)m(\zeta_2)\big). \label{E0(u)}
\end{align}
From Lemma \ref{Lambda_k}, one has
\begin{align}
 \frac{dE^0(u)}{dt}&=2i\Lambda_3\big(m^2(\zeta_1)(\xi_1+\eta_1)\big)\notag \\
  &=\frac{2i}{3}\Lambda_3(M_3),  \label{dE0(u)}
\end{align}
where $M_3=\sum^3_{j=1}m^2(\zeta_j)(\xi_j+\eta_j)$.

Modifying the quantity $E^0(u)$ so that the time derivative has less of a
$\Lambda_3$ term, we define
\begin{align}
E^1(u)=E^0(u)+\Lambda_3(\sigma_3)  \nonumber
\end{align}
for some $\sigma_3$ to be chosen shortly. Computing as before we have
\begin{align}
 \frac{dE^1(u)}{dt}=\Lambda_3(\frac{2iM_3}{3}+\sigma_3h_3)-3i\Lambda_4\big(
 X(\sigma_3)\big),  \nonumber
\end{align}
where
$$h_3=i(\xi^3_1+\eta^3_1+\xi^3_2+\eta^3_2+\xi^3_3+\eta^3_3)=3i(\xi_1\xi_2\xi_3+\eta_1\eta_2
\eta_3),$$
and
$$X(\sigma_3)=\sigma_3(\zeta_1,\zeta_2,\zeta_3+\zeta_4)\cdot(\xi_3+\xi_4+\eta_3+\eta_4).$$

An intuitional guess for $\sigma_3$ would thus be
$$\sigma_3=-\frac{2iM_3}{3h_3}.$$
However this choice runs into the problem that $h_3$ can vanish in the resonant interaction
case when $\xi_1\xi_2\xi_3$ equals to $-\eta_1\eta_2\eta_3$. In particular, when all
frequencies are less than $N$, then $M_3=0$ and so the vanishing of the denominator is
cancelled by the numerator. Unfortunately, this cancellation is lost when we have one or
more high frequencies. This is in contrast to the one-dimensional situation in \cite{},
where the resonant interactions is removable.

We shall in fact set
$$\tilde{\sigma}_3=-\frac{2iM_3}{3h_3}1_{\Omega_{nr}},$$
where $1_{\Omega_{nr}}$ is the indicator function to the non-resonant set
\begin{align}\Omega_{nr}:=&\left\{ (\zeta_1,\zeta_2,\zeta_3)\in \Sigma_3 \ \Big|\ \
\max_{1\leq j\leq 3}|\zeta_j|\leq N \right\}\notag \\
&\cup \left\{ (\zeta_1,\zeta_2,\zeta_3)\in \Sigma_3 \ \Big|\ \
|\xi_1\xi_2\xi_3+\eta_1\eta_2\eta_3|\geq \gamma_0N_1N_2N_3 \right\}
.
\end{align}
where
$$\Sigma_3=\left\{ (\zeta_1,\zeta_2,\zeta_3)\in \mathbb{R}^{2\times 3} \ \Big|\ \
\zeta_1+\zeta_2+\zeta_3=0 \right\}$$
and $0<\theta<2^{-30}$ is a parameter to be chosen later.

We now define $\tilde{E}^1(u)$ by
\begin{align}
\tilde{E}^1(u)=E^0(u)+\Lambda_3(\tilde{\sigma}_3),\label{tildeE1(u)}
\end{align}
and consider the growth of $\tilde{E}^1(u)$.

A simple computation implies
\begin{align}
 \frac{d\tilde{E}^1(u)}{dt}=\Lambda_3(\frac{2iM_3}{3}+\tilde{\sigma}_3h_3)-3i\Lambda_4\big(
 X(\tilde{\sigma}_3)\big),  \label{dtildeE^1(u)}
\end{align}
where
$$X(\tilde{\sigma}_3)=\tilde{\sigma}_3(\zeta_1,\zeta_2,\zeta_3+\zeta_4)\cdot(\xi_3+\xi_4+
\eta_3+\eta_4).$$
Integrating in time, one has
\begin{align}
&\tilde{E}^1(u)(\delta)- \tilde{E}^1(u)(0) \notag \\
=&\int^{\delta}_0\Lambda_3(\frac{2iM_3}{3}+\tilde{\sigma}_3h_3)dt-3i\int^{\delta}_0\Lambda_4
\big(
 X(\tilde{\sigma}_3)\big)dt. \label{In1}
\end{align}

To prove Theorem \ref{gwp}, it thus suffices to prove the following three propositions.
\begin{proposition} \label{fixedtimep1}
Let $-\frac{1}{13}<s<0$. We have
\begin{align}
 |\Lambda_3(\tilde{\sigma}_3)|\lesssim
\gamma_0^{-1}N^{-1+}\|Iu\|^3_{L^2(\mathbb{R}^2)}.\label{fixedtimep1a}
\end{align}
\end{proposition}

\begin{proposition}[Trilinear estimate] \label{p3}
Let $-\frac{1}{13}<s<0$. We have
\begin{align}
 \Big|\int^{\delta}_0\Lambda_3(\frac{2iM_3}{3}+\tilde{\sigma}_3h_3)dt\Big|\lesssim
 \big(\gamma_0^{-\frac{1}{2}}N^{-\frac{1}{2}+}+N^{-\frac{1}{4}+}\big)\|Iu\|^3_{X^{0,
\frac{1}{2}+}_{\delta}} .\label{Lambda(3)}
\end{align}
\end{proposition}

\begin{proposition}[Quadrilinear estimate]\label{p4}
Let $-\frac{1}{13}<s<0$. We have
\begin{align}
 \Big|\int^{\delta}_0\Lambda_4\big(X(\tilde{\sigma}_3)\big)dt\Big|\lesssim
\gamma_0^{-1}N^{-2+}\|Iu\|^4_{X^{0,\frac{1}{2}+}_{\delta}} .\label{Lambda(4)}
\end{align}
\end{proposition}

Before extending to a global solution, we need to control the smoothed solution for small
time.
\begin{proposition}[Modified local existence] \label{p2}
Let $-\frac{1}{13}<s<0$. Assume $u_0$ satisfies $\|Iu_0\|_{L^2(\mathbb{R}^2)}\leq 1$. Then
there is a constant $\delta=\delta(\|u_0\|_{L^2(\mathbb{R}^2)})$ and a unique solution $u$
to \eqref {ZK} on $[0,\delta]$,
 such that

$$\|Iu\|_{X^{0, \frac{1}{2}+}_\delta}\lesssim 1$$
where the implicit constant is independent of $\delta$.
\end{proposition}
{\bf Proof.} We use the iteration argument to show the local well-posedness.

Acting multiplier operator $I$ on both sides of \eqref {ZK},
\begin{align}
\partial_{t}Iu + (\partial^3_{x}+\partial^3_{y}) Iu + (\partial_{x}+\partial_{y}) I(u^2) =
0.\label{equationI}
\end{align}
By Duhamel's principle, one can rewrite the differential equation as an integral equation
$$Iu=e^{-t(\partial^3_{x}+\partial^3_{y})}Iu_0-\int^t_0e^{-(t-t')(\partial^3_{x}+
\partial^3_{y})}
(\partial_{x}+\partial_{y}) I(u^2)(t')dt'.$$
Estimates \eqref {estimate1}-\eqref {estimate3} give us
\begin{align}
\|Iu\|_{X^{0, \frac{1}{2}+}_\delta}
&\lesssim\|e^{-t(\partial^3_{x}+\partial^3_{y})}Iu_0\|_{X^{0,\frac{1}{2}+}_\delta}\notag \\
&\hspace{5mm}+\Big\|\int^t_0e^{-(t-t')(\partial^3_{x}+\partial^3_{y})}(\partial_{x}+
\partial_{y})I(u^2)(t')dt'\Big\|
_{X^{0, \frac{1}{2}+}_\delta}  \notag \\
& \lesssim \|Iu_0\|_{L^2(\mathbb{R}^2)} +\|(\partial_{x}+\partial_{y})I(u^2)\|_{X^{0,
-\frac{1}{2}+}_\delta}    \notag \\
&   \lesssim \|Iu_0\|_{L^2(\mathbb{R}^2)}
+\delta^{0+}\|(\partial_{x}+\partial_{y})I(u^2)\|_{X^{0, -\frac{1}{2}++}_\delta}.
\label{3.5}
\end{align}

According to the definition of the restricted norm \eqref{X1 Norm}, we can choose $\tilde u
\in X^{0, \frac{1}{2}+} $ satisfying $\tilde u|_{[0, \delta]}=u$,
\begin{align}\|Iu\|_{X^{0, \frac{1}{2}+}_\delta} \sim \|I\tilde u\|_{X^{0,
\frac{1}{2}+}}\end{align} \label{3.6}
and
\begin{align}
\|(\partial_{x}+\partial_{y})I(u^2)\|_{X^{0, -\frac{1}{2}++}_\delta}\lesssim
\|(\partial_{x}+\partial_{y})I(\tilde u^2)\|_{X^{0, -\frac{1}{2}++}}. \label{3.7}
\end{align}
We will show shortly that
\begin{align}
\|(\partial_{x}+\partial_{y})I(\tilde u^2)\|_{X^{0, -\frac{1}{2}++}}\lesssim  \|I\tilde
u\|^2_{X^{0, \frac{1}{2}+}}.\label{3.8}
\end{align}
Using the Lemma \ref{L5}, it suffices to prove
\begin{align}
\|I^s_1\partial_{x}(\tilde u^2)\|_{X^{0, -\frac{1}{2}++}}\lesssim  \|I^s_1\tilde
u\|^2_{X^{0, \frac{1}{2}+}}. \label{3.9}
\end{align}
Note that $\|I^s_1F\|_{X^{0, b}} \sim \|F\|_{X^{s, b}}$, \eqref{3.9} is an immediate
consequence of \eqref{estimate4}.

Combining \eqref{3.5}-\eqref{3.8}, we have
\begin{align}\|Iu\|_{X^{0, \frac{1}{2}+}_\delta}\lesssim \|Iu_0\|_{L^2(\mathbb{R}^2)}
+\delta^{0+}\|Iu\|^2_{X^{0, \frac{1}{2}+}_\delta},  \label{3.10} \end{align}
and
$$\|Iu-Iv\|_{X^{0, \frac{1}{2}+}_\delta}\lesssim \delta^{0+}(\|Iu\|_{X^{0,
\frac{1}{2}+}_\delta}+\|Iv\|_{X^{0, \frac{1}{2}+}_\delta})\|Iu-Iv\|_{X^{0,
\frac{1}{2}+}_\delta}.$$
Then, one can obtain the local well-posedness by means of the contraction mapping principle.

Setting $Q(\delta)\equiv \|Iu\|_{X^{0, \frac{1}{2}+}_\delta}$, the bound \eqref{3.10} yields
\begin{align}Q(\delta)\lesssim \|Iu_0\|_{L^2(\mathbb{R}^2)}+\delta^{0+}(Q(\delta))^2.
\nonumber \end{align}
As $\|Iu_0\|_{L^2(\mathbb{R}^2)}\leq 1$ and $Q(\delta)$ is continuous in the variable
$\delta$, a bootstrap argument yields $Q(\delta)\lesssim 1$, i.e.
$\|Iu\|_{X^{0,\frac{1}{2}+}_\delta}\lesssim 1.$
\begin{theorem} \label {Prop E incresement}
Let $-\frac{1}{13}<s<0, N\gg 1$. Assume $u_0$ satisfies $\|Iu_0\|_{L^2(\mathbb{R}^2)}\leq
1$. Then there is a constant $\delta=\delta(\|u_0\|_{L^2(\mathbb{R}^2)})>0$ so that there
exists a unique solution
$$u(x,y,t)\in C([0,\delta],H^s(\mathbb{R}^2))\cap X^{s,\frac{1}{2}+}_{\delta}$$
 of \eqref {ZK} satisfying
\begin{align}
\tilde{E}^1(u)(\delta)=\tilde{E}^1(u)(0)+O(N^{-\frac{1}{4}+}). \label {tildeE incresement}
 \end{align}
\end{theorem}
{\bf Proof.} From Proposition \ref{p2} , there exists a unique solution $u$ to \eqref {ZK}
on $[0,\delta]$ satisfying $\|Iu\|_{X^{0, \frac{1}{2}+}_\delta}\lesssim 1$.

Choosing $\gamma_0=N^{-\frac{1}{2}}$, it follows from \eqref {In1}, Proposition \ref{p3} and
Proposition \ref {p4} that
\begin{align}
 &|\tilde{E}^1(u)(\delta)-\tilde{E}^1(u)(0)| \notag \\
 \lesssim &\Big|\int^{\delta}_0\Lambda_3(\frac{2iM_3}{3}+\tilde{\sigma}_3h_3)dt\Big|
 +\Big|\int^{\delta}_0\Lambda_4\big(X(\tilde{\sigma}_3)\big)dt\Big|
 \notag \\
 \lesssim &
N^{-\frac{1}{4}+}\|Iu\|^3_{X^{0,\frac{1}{2}+}_{\delta}}+N^{-\frac{3}{2}+}\|Iu\|^4_{X^{0,
\frac{1}{2}+}_{\delta}}\notag \\
  \lesssim &  N^{-\frac{1}{4}+}.\nonumber
\end{align}
Thus we prove this theorem.

{\bf Proof of Theorem \ref{gwp}.} Fix $T>0$, $u_0\in L^2(\mathbb{R}^2)$. Let $\lambda\ll1$
be a scaling parameter to be chosen later. We define the rescaled solution $u_\lambda:
[0,\frac{T}{\lambda^3}]\times\mathbb{R}^2\rightarrow\mathbb{R}$
$$u_\lambda(x,y,t)=\lambda^2u(\lambda x,\lambda y,\lambda^3t).$$

The modified mass functional can be arbitrarily small by taking $\lambda$ small,
\begin{align}
E^0(u_{\lambda,0})&=\|Iu_{\lambda,0}\|^2_{L^2(\mathbb{R}^2)}\notag \\
&\leq N^{-2s}\|u_{\lambda,0}\|^2_{H^s(\mathbb{R}^2)}\notag \\
 & \leq N^{-2s}\lambda^{2(s+1)}\|u_0\|^2_{H^s(\mathbb{R}^2)}.\nonumber
\end{align}
Assuming $N\gg1$ is given ($N$ will be chose shortly), we choose
\begin{align}
\lambda=\lambda(N,\|u_0\|_{H^s(\mathbb{R}^2)})\sim N^{\frac{s}{s+1}}\nonumber
\end{align}
such that $E^0(u_{\lambda,0})\leq \frac{1}{8}$.

Now we can apply Proposition \ref {Prop E incresement} to the scaled initial data
$u_{\lambda,0}$, hence
$$\tilde{E}^1(u_\lambda)(\delta)=\tilde{E}^1(u_\lambda)(0)+O(N^{-\frac{1}{4}+}). $$
From \eqref{tildeE1(u)} and Proposition \ref{fixedtimep1}, we have
\begin{align}
E^0(u_{\lambda})(\delta)&\leq E^0(u_{\lambda,0})+|\Lambda_3(\tilde{\sigma}_3)
(\delta)|+|\Lambda_3(\tilde{\sigma}_3)(0)|\notag \\
&\hspace{4mm}+|\tilde{E}^1(u_\lambda)(\delta)-\tilde{E}^1(u_\lambda)(0)|\notag \\
&\leq\frac{1}{4}+CN^{-\frac{1}{4}+}<\frac{1}{2},\nonumber
\end{align}
where the second step holds true because
\begin{align}
|\Lambda_3(\tilde{\sigma}_3)(\delta)|&\leq C_1 \gamma_0^{-1}
N^{-1+}\|Iu_\lambda\|^2_{L^2(\mathbb{R}^2)}\notag \\
&\leq C_2 \gamma_0^{-1} N^{-1+}\|u_\lambda\|^2_{L^2(\mathbb{R}^2)}\notag \\
&\leq C_2\gamma_0^{-1} N^{-1+}\|u_{\lambda,0}\|^2_{L^2(\mathbb{R}^2)}\notag \\
&\leq C_2\gamma_0^{-1} N^{-1+} \lambda^2\|u_{0}\|^2_{L^2(\mathbb{R}^2)}\notag \\
&\leq \frac{1}{16}\nonumber
\end{align}
as long as $N$ is sufficiently large.

Thus, from Proposition \ref {p2}, the solution $u_\lambda$ can be extended to
$t\in[0,2\delta]$ .

Iterating this procedure $M$ steps, where $M\lesssim N^{\frac{1}{4}-}$, we have
$$E^0(u_\lambda)(t)\leq\frac{1}{4}+CMN^{-\frac{1}{4}+}\leq1$$
 for $t\in[0,(M+1)\delta]$. Therefore $u_\lambda$ can be extended to
$t\in[0,N^{\frac{1}{4}-}\delta].$

Taking $N(T)$ sufficiently large such that
$$N^{\frac{1}{4}-}\delta >\frac{T}{\lambda^3}\sim TN^{-\frac{3s}{s+1}},$$
which implies
$$T\sim N^{\frac{13s+1}{4(s+1)}-}.$$
Since $s>-\frac{1}{13}$, the exponent of $N$ above is positive. Thus we get the global
well-posedness for \eqref{ZK}.                                                                            \section{The fixed time estimate}

We show Proposition \ref{fixedtimep1} in this section.
\begin{lemma}\label{FTI1}
For any $(\zeta_1,\zeta_2,\zeta_3)\in\Sigma_3$, we have
\begin{align}
\big|M_3(\zeta_1,\zeta_2,\zeta_3)\big|\lesssim \max_{1\leq
j\leq3}\{m^2(\zeta_j)\}\min_{1\leq j\leq3}\{|\zeta_j|\},\nonumber
\end{align}
where $M_3=\sum^3_{j=1}m^2(\zeta_j)(\xi_j+\eta_j)$.
\end{lemma}
{\bf Proof.} Set $|\zeta_j|\sim N_j$, we know that $M_3$ vanishes when $\max_{1\leq
j\leq3}\{|\zeta_j|\}\leq N$. From symmetry, one can assume that $N_1\sim N_2\gtrsim N_3$ and
$N_1\geq N$.

Let $f(\zeta)=m^2(\zeta)(\xi+\eta)$. It suffices to prove that
$$\big|f(\zeta_1)+f(\zeta_2)+f(\zeta_3)\big|\lesssim \max_{1\leq
j\leq3}\{m^2(\zeta_j)\}\min_{1\leq j\leq3}\{|\zeta_j|\}.$$
{\bf Case 1}  \hspace{2mm}$N_1\sim N_2\sim N_3$

We can estimate all three terms on the left-hand side by $O\big(m^2(\zeta_1)|\zeta_1|\big)$.
\\{\bf Case 2} \hspace{2mm}$N_1\sim N_2\gg N_3$

Note that $\big|\nabla f(\zeta)\big|=O\big(m^2(\zeta)\big)$, so by the mean value theorem we
have
\begin{align}
\big|f(\zeta_1)+f(\zeta_2)\big|=&\big|f(\zeta_1)-f(\zeta_1+\zeta_3)\big|\notag \\
\lesssim &\big|\nabla f(\tilde{\zeta})\big|\cdot|\zeta_3|\notag \\
\lesssim &m^2(\zeta_1)\cdot |\zeta_3|\nonumber
\end{align}
where $|\tilde{\zeta}|\sim|\zeta_1|\sim|\zeta_2|$.
Hence
\begin{align}
&\big|f(\zeta_1)+f(\zeta_2)+f(\zeta_3)\big|\notag \\
\leq&\big|f(\zeta_1)+f(\zeta_2)\big|+\big|f(\zeta_3)\big|\notag \\
\lesssim &\big(m^2(\zeta_1)+m^2(\zeta_3)\big)\cdot |\zeta_3|,\nonumber
\end{align}
and the claim follows.

{\bf Proof of Proposition \ref{fixedtimep1}.} The strategy is that we first treat the dyadic
constituent $P_{N_j}u$ of which the frequency support is $|\zeta_j|\sim N_j$, then we
conclude our desired bound  by summing over all dyadic pieces $P_{N_j}u$.

For simplicity, we write $\int_{*}fgh:=\int_{\Sigma}f(\zeta_1)g(\zeta_2)h(\zeta_3)d\sigma$,
where $\mu$ is the induced measure on hypersurface $\Sigma:=\{\zeta_1+\zeta_2+\zeta_3=0\}$.
Without loss of generality we assume $\widehat {u}(\zeta)$ is non-negative. One can also
assume that $N_1\sim N_2\gtrsim N_3$ and $N_1\geq N$.

According to the notation above and Lemma \ref{FTI1}, we observe that
\begin{align}
 &\Big|\Lambda_3\big(\frac{M_3(\zeta_1,\zeta_2,\zeta_3)1_{\{|\xi_1\xi_2\xi_3+\eta_1\eta_2
\eta_3|\geq \gamma_0N_1N_2N_3\}}}{(\xi_1\xi_2\xi_3+\eta_1\eta_2\eta_3)\cdot
m(\zeta_1)m(\zeta_2)m(\zeta_3)};u_1,u_2,u_3\big)\Big|
 \notag \\
 \lesssim&
\gamma_0^{-1}N_{1}^{-1}N_{2}^{-1}N_{3}^{-1}\frac{m^2(\zeta_3)N_3}{m(\zeta_1)m(\zeta_2)
m(\zeta_3)}\int_{*}\widehat{u}^{xy}_1\cdot
 \widehat{u}^{xy}_2\cdot\widehat{u}^{xy}_3 \notag \\
 \lesssim& \gamma_0^{-1}N_{1}^{-1+}\frac{1}{\big(\frac{N_{1}}{N}\big)^{2s}}\int_{*} \langle
\zeta_1 \rangle^{-\frac{1}{2}-}  \widehat{u}^{xy}_1\cdot \langle \zeta_2
\rangle^{-\frac{1}{2}-}\widehat{u}^{xy}_2\cdot\widehat{u}^{xy}_3 \notag \\
 \lesssim& \gamma_0^{-1}N^{-1+}\|\mathscr F_{\zeta}^{-1}\langle \zeta_1
\rangle^{-\frac{1}{2}-}  \widehat{u}^{xy}_1\|_{L^4(\mathbb{R}^2)} \|\mathscr
F_{\zeta}^{-1}\langle \zeta_2
\rangle^{-\frac{1}{2}-}\widehat{u}^{xy}_2\|_{L^4(\mathbb{R}^2)} \|u_3\|_{L^2(\mathbb{R}^2)}
\notag \\
  \lesssim& \gamma_0^{-1}N^{-1+}\|\mathscr F_{\zeta}^{-1}\langle \zeta_1
\rangle^{-\frac{1}{2}-}  \widehat{u}^{xy}_1\|_{H^{\frac{1}{2}+}(\mathbb{R}^2)} \|\mathscr
F_{\zeta}^{-1}\langle \zeta_2
\rangle^{-\frac{1}{2}-}\widehat{u}^{xy}_2\|_{H^{\frac{1}{2}+}(\mathbb{R}^2)}
\|u_3\|_{L^2(\mathbb{R}^2)}     \notag \\
 \lesssim& \gamma_0^{-1}N^{-1+}
\|u_1\|_{L^2(\mathbb{R}^2)}\|u_2\|_{L^2(\mathbb{R}^2)}\|u_3\|_{L^2(\mathbb{R}^2)}.
\nonumber
\end{align}
This leads to \eqref{fixedtimep1a}.

\section{Quadrilinear estimate}

In this section, we prove Proposition \ref{p4}  which is easier than Proposition \ref{p3}.

If $\max_{1\leq j\leq4}\{|\zeta_j|\}\leq \frac{N}{2}$, then
$M_3(\zeta_1,\zeta_2,\zeta_3+\zeta_4)=0$, and thus $X(\tilde{\sigma}_3)=0$. So we shall
assume that $\max_{1\leq j\leq4}\{|\zeta_j|\}\geq \frac{N}{2}$. Moreover, one can assume
that $N_1\geq N_2$, $N_3\geq N_4$ and  $N_1\gtrsim N$ by symmetries.

Using the dyadic decomposition, it suffices to prove that
\begin{align}
  &|\int^{\delta}_0\int_{*}\frac{M_3(\zeta_1,\zeta_2,\zeta_3+\zeta_4)\cdot(\xi_3+\xi_4+
\eta_3+\eta_4)}
  {\big(\xi_1\xi_2(\xi_3+\xi_4)+\eta_1\eta_2(\eta_3+\eta_4)\big)\cdot
m(\zeta_1)m(\zeta_2)m(\zeta_3)m(\zeta_4)}\widehat{u}_1(\zeta_1)\cdots
\widehat{u}_4(\zeta_4)dt|\notag \\
  \lesssim &\gamma_0^{-1}N^{-2}\prod_{j=1}^4\|u_j\|_{X^{0,\frac{1}{2}+}_{\delta}}
\label{P4a}
\end{align}
for $|\xi_1\xi_2(\xi_3+\xi_4)+\eta_1\eta_2(\eta_3+\eta_4)|\geq \gamma_0N_1N_2N_{34}$, where
$|\zeta_3+\zeta_4|\sim N_{34}$.
\\{\bf Case 1} \hspace{2mm}$N_1\gg N_3$, $N_1\sim N_2\gtrsim N$

Using Lemma \ref{FTI1} and \eqref{bles4}, we have
\begin{align}
 LHS \  of \ \eqref {P4a}\lesssim &
\gamma_0^{-1}N_{1}^{-1}N_{2}^{-1}N_{1}^{-2s}N^{2s}\frac{m^2(\zeta_{34})N_{34}}{m(\zeta_3)
m(\zeta_4)}\|u_1u_3\|_{L^2(\mathbb{R}^3)} \|u_2u_4\|_{L^2(\mathbb{R}^3)} \notag \\
  \lesssim &\gamma_0^{-1}N^{2s}N_{1}^{-2s-2}N_{34}\frac{N_3^{\frac{1}{2}}N_4^{\frac{1}{2}}}
{N_{1}^{2}}\frac{m^2(\zeta_{34})}{m(\zeta_3)m(\zeta_4)}\prod_{j=1}^4\|u_j\|_{X^{0,\frac{1}{2}
+}_{\delta}}. \nonumber
\end{align}
In order to obtain  \eqref{P4a}, we only need to bound
$$\gamma_0^{-1}N^{2s}N_{1}^{-2s-2}N_{34}\cdot\frac{N_3^{\frac{1}{2}}N_4^{\frac{1}{2}}}
{N_{1}^{2}}\cdot\frac{m^2(\zeta_{34})}
{m(\zeta_3)m(\zeta_4)}$$
by $\gamma_0^{-1}N^{-2}$.
\\{\bf Case 1.1} \hspace{2mm} $N\gtrsim N_3\geq N_4$

As $m(\zeta_{34})=m(\zeta_{3})=m(\zeta_{4})=1$, we have
\begin{align}
 &\gamma_0^{-1}N^{2s}N_{1}^{-2s-2}N_{34}\cdot\frac{N_3^{\frac{1}{2}}N_4^{\frac{1}{2}}}
{N_{1}^{2}}\cdot\frac{m^2(\zeta_{34})}
{m(\zeta_3)m(\zeta_4)}\notag \\
\lesssim & \gamma_0^{-1}N^{2s+2}N_{1}^{-2s-4} \notag \\
\lesssim &\gamma_0^{-1}N^{-2}. \nonumber
\end{align}
{\bf Case 1.2} \hspace{2mm} $N_3\gg N\gtrsim N_4$

As $m(\zeta_{34})\sim m(\zeta_{3})\sim(\frac{N_3}{N})^s$ and $m(\zeta_{4})\sim 1$, we have
\begin{align}
 &\gamma_0^{-1}N^{2s}N_{1}^{-2s-2}N_{34}\cdot\frac{N_3^{\frac{1}{2}}N_4^{\frac{1}{2}}}
{N_{1}^{2}}\cdot\frac{m^2(\zeta_{34})}
{m(\zeta_3)m(\zeta_4)}\notag \\
\lesssim &
\gamma_0^{-1}N^{2s}N_{1}^{-2s-4}N_{3}^{\frac{3}{2}}N_{4}^{\frac{1}{2}}\big(\frac{N_3}
{N}\big)^s \notag \\
\lesssim &\gamma_0^{-1}N^{-2}. \nonumber
\end{align}
{\bf Case 1.3} \hspace{2mm} $N_3\gg N_4\gg N$
\begin{align}
 &\gamma_0^{-1}N^{2s}N_{1}^{-2s-2}N_{34}\cdot\frac{N_3^{\frac{1}{2}}N_4^{\frac{1}{2}}}
{N_{1}^{2}}\cdot\frac{m^2(\zeta_{34})}
{m(\zeta_3)m(\zeta_4)}\notag \\
\lesssim &
\gamma_0^{-1}N^{2s}N_{1}^{-2s-4}N_{3}^{\frac{3}{2}}N_{4}^{\frac{1}{2}}\cdot\frac{m(\zeta_{3})}
{m(\zeta_4)}\notag \\
\lesssim &
\gamma_0^{-1}N^{2s}N_{1}^{-2s-4}N_{3}^{\frac{3}{2}}N_{4}^{\frac{1}{2}}\big(\frac{N_3}
{N_4}\big)^s \notag \\
\lesssim &\gamma_0^{-1}N^{-2}. \nonumber
\end{align}
{\bf Case 1.4} \hspace{2mm} $N_3\sim N_4\gg N$
\begin{align}
 &\gamma_0^{-1}N^{2s}N_{1}^{-2s-2}N_{34}\cdot\frac{N_3^{\frac{1}{2}}N_4^{\frac{1}{2}}}
{N_{1}^{2}}\cdot\frac{m^2(\zeta_{34})}
{m(\zeta_3)m(\zeta_4)}\notag \\
\lesssim & \gamma_0^{-1}N^{2s}N_{1}^{-2s-4}N_{3}^{2}\cdot\frac{1}
{m^2(\zeta_3)}\notag \\
\lesssim &\gamma_0^{-1}N^{-2}. \nonumber
\end{align}
{\bf Case 2} \hspace{2mm} $N_1\sim N_3$

We shall divide this case into four subcases depending on relationships between $N_1$ and
$N_2, N_4$.
\\{\bf Case 2.1} \hspace{2mm} $N_1\gg N_2$ and $N_1\sim N_3\gg N_4$

According to Lemma \ref{FTI1} and \eqref{bles4}, we have
\begin{align}
 LHS \  of \ \eqref {P4a}
 \lesssim &
\big(\gamma_0N_{1}N_{2}N_{34}\big)^{-1}\frac{m^2(\zeta_2)N_{2}N_{34}}{m^2(\zeta_1)m(\zeta_2)
m(\zeta_4)}\|u_1u_2\|_{L^2(\mathbb{R}^3)} \|u_3u_4\|_{L^2(\mathbb{R}^3)} \notag \\
  \lesssim
&\big(\gamma_0N^2_{1}N_{2}\big)^{-1}N_{1}^{-2s}N^{2s}\frac{N_2^{\frac{1}{2}}N_4^{\frac{1}{2}}}
{N_{1}^{2}}\frac{1}{m(\zeta_4)}
  \prod_{j=1}^4\|u_j\|_{X^{0,\frac{1}{2}+}_{\delta}}. \notag \\
  \lesssim &\gamma_0^{-1}N^{-2}\prod_{j=1}^4\|u_j\|_{X^{0,\frac{1}{2}+}_{\delta}}.\nonumber
\end{align}
{\bf Case 2.2} \hspace{2mm} $N_1\sim N_3\sim N_2\gg N_4$
\begin{align}
 LHS \  of \ \eqref {P4a}\lesssim &
\theta^{-1}_0N^{-3}_{1}\frac{m^2(\zeta_1)N_{1}N_{34}}{m^3(\zeta_1)m(\zeta_4)}\|u_1u_2\|
_{L^2(\mathbb{R}^3)} \|u_3u_4\|_{L^2(\mathbb{R}^3)} \notag \\
  \lesssim
&\theta^{-1}_0N^{-1}_{1}N_{1}^{-s}N^{s}\frac{N_2^{\frac{1}{2}}N_4^{\frac{1}{2}}}{N_{1}^{2}}
\frac{1}{m(\zeta_4)}
  \prod_{j=1}^4\|u_j\|_{X^{0,\frac{1}{2}+}_{\delta}}. \notag \\
  \lesssim &\gamma_0^{-1}N^{-2}\prod_{j=1}^4\|u_j\|_{X^{0,\frac{1}{2}+}_{\delta}}.\nonumber
\end{align}
{\bf Case 2.3} \hspace{2mm} $N_1\sim N_3\sim N_4\gg N_2$
\begin{align}
 LHS \  of \ \eqref {P4a}\lesssim &
\big(\gamma_0N_{1}N_{2}N_{34}\big)^{-1}\frac{m^2(\zeta_2)N_{2}N_{34}}{m^3(\zeta_1)m(\zeta_2)}
\|u_1u_2\|_{L^2(\mathbb{R}^3)} \|u_3u_4\|_{L^2(\mathbb{R}^3)} \notag \\
  \lesssim
&\theta^{-1}_0N^{-1}_{1}N_{1}^{-3s}N^{3s}\frac{N_2^{\frac{1}{2}}N_4^{\frac{1}{2}}}{N_{1}^{2}}
  \prod_{j=1}^4\|u_j\|_{X^{0,\frac{1}{2}+}_{\delta}}. \notag \\
  \lesssim &\gamma_0^{-1}N^{-2}\prod_{j=1}^4\|u_j\|_{X^{0,\frac{1}{2}+}_{\delta}}.\nonumber
\end{align}
{\bf Case 2.4} \hspace{2mm} $N_1\sim N_2\sim N_3\sim N_4$
\begin{align}
 LHS \  of \ \eqref {P4a}\lesssim &
\big(\gamma_0N_{1}N_{2}N_{34}\big)^{-1}\frac{m^2(\zeta_{34})N^2_{34}}{m^4(\zeta_1)}\|u_1u_2\|
_{L^2(\mathbb{R}^3)} \|u_3u_4\|_{L^2(\mathbb{R}^3)} \notag \\
  \lesssim
&\theta^{-1}_0N^{-2}_{1}N_{34}N_{1}^{-4s}N^{4s}m^2(\zeta_{34})\frac{N_2^{\frac{1}{2}}
N_4^{\frac{1}{2}}}{N_{1}^{2}}
  \prod_{j=1}^4\|u_j\|_{X^{0,\frac{1}{2}+}_{\delta}}. \notag \\
  \lesssim &\gamma_0^{-1}N^{-2}\prod_{j=1}^4\|u_j\|_{X^{0,\frac{1}{2}+}_{\delta}}.\nonumber
\end{align}

Next we consider the last case $N_3\gg N_1$ which implies $N_3\sim N_4$.
\\{\bf Case 3} \hspace{2mm} $N_3\sim N_4\gg N_1\geq N_2$

We should investigate the size relationship between $N_1$ and $ N_2$.
\\{\bf Case 3.1} \hspace{2mm} $N_1\sim N_2\sim N_{34}$
\begin{align}
 LHS \  of \ \eqref {P4a}\lesssim &
\big(\gamma_0N_{1}N_{2}N_{34}\big)^{-1}\frac{m^2(\zeta_{1})N^2_{1}}{m^2(\zeta_1)m^2(\zeta_3)}
\|u_1u_3\|_{L^2(\mathbb{R}^3)} \|u_2u_4\|_{L^2(\mathbb{R}^3)} \notag \\
  \lesssim
&\theta^{-1}_0N^{-1}_{1}N_{3}^{-2s}N^{2s}\frac{N_1^{\frac{1}{2}}N_2^{\frac{1}{2}}}{N_{3}^{2}}
  \prod_{j=1}^4\|u_j\|_{X^{0,\frac{1}{2}+}_{\delta}}. \notag \\
  \lesssim &\gamma_0^{-1}N^{-2}\prod_{j=1}^4\|u_j\|_{X^{0,\frac{1}{2}+}_{\delta}}.\nonumber
\end{align}
{\bf Case 3.2} \hspace{2mm} $N_1\sim N_2\gg N_{34}$
\begin{align}
 LHS \  of \ \eqref {P4a}\lesssim &
\big(\gamma_0N_{1}N_{2}N_{34}\big)^{-1}\frac{m^2(\zeta_{34})N^2_{34}}{m^2(\zeta_1)m^2
(\zeta_3)}\|u_1u_3\|_{L^2(\mathbb{R}^3)} \|u_2u_4\|_{L^2(\mathbb{R}^3)} \notag \\
  \lesssim
&\theta^{-1}_0N^{-2}_{1}N_{34}\frac{m^2(\zeta_{34})}{m^2(\zeta_{1})}\big(\frac{N_3}
{N}\big)^{-2s}\frac{N_1^{\frac{1}{2}}N_2^{\frac{1}{2}}}{N_{3}^{2}}
  \prod_{j=1}^4\|u_j\|_{X^{0,\frac{1}{2}+}_{\delta}}. \notag \\
  \lesssim &\gamma_0^{-1}N^{-2}\prod_{j=1}^4\|u_j\|_{X^{0,\frac{1}{2}+}_{\delta}}.\nonumber
\end{align}
{\bf Case 3.3} \hspace{2mm} $N_1\sim N_{34}\gg N_{2}$
\begin{align}
 LHS \  of \ \eqref {P4a}\lesssim &
\big(\gamma_0N_{1}N_{2}N_{34}\big)^{-1}\frac{m^2(\zeta_{2})N_2N_{34}}{m(\zeta_1)m(\zeta_2)
m^2(\zeta_3)}\|u_1u_3\|_{L^2(\mathbb{R}^3)} \|u_2u_4\|_{L^2(\mathbb{R}^3)} \notag \\
  \lesssim
&\theta^{-1}_0N^{-1}_{1}\frac{m(\zeta_{2})}{m(\zeta_{1})}\big(\frac{N_3}{N}\big)^{-2s}
\frac{N_1^{\frac{1}{2}}N_2^{\frac{1}{2}}}{N_{3}^{2}}
  \prod_{j=1}^4\|u_j\|_{X^{0,\frac{1}{2}+}_{\delta}}. \notag \\
  \lesssim &\gamma_0^{-1}N^{-2}\prod_{j=1}^4\|u_j\|_{X^{0,\frac{1}{2}+}_{\delta}}.\nonumber
\end{align}

This concludes the proof of Proposition \ref {p4}.                                                              \section{Trilinear estimate}

This section is devoted to the proof of Proposition \ref {p3}.

As we can see from the symbol $\frac{2iM_3}{3}+\tilde{\sigma}_3h_3$ that
$$\frac{2iM_3}{3}+\tilde{\sigma}_3h_3=0$$
if $|\xi_1\xi_2\xi_3+\eta_1\eta_2\eta_3|\geq \gamma_0N_1N_2N_{3}$, we need a better control
on the resonant interactions.

By Parseval formula, the left-hand side of \eqref{Lambda(3)} can be expanded as
\begin{align}
\int_{\mathbb{R}}\int_{\sum^3_{{\iota}=1}\tau_{\iota}=\tau_0,\sum^3_{{\iota}=1}
\zeta_{\iota}=0}\mathcal{M}_3(\zeta_1,\zeta_2,\zeta_3)\widehat{1}_{[0,
\delta]}(\tau_0)\widehat{u}_{1}(\lambda_1)\widehat{u}_{2}(\lambda_2)\widehat{u}_{3}
(\lambda_3)d\tau_0,\label{M1zeta}
\end{align}
where
$$\mathcal{M}_3(\zeta_1,\zeta_2,\zeta_3)=\frac{\sum_{{\iota}=1}^3m^2(\zeta_{\iota})
(\xi_{\iota}+\eta_{\iota})}{m(\zeta_1)m(\zeta_2)m(\zeta_3)}1_{\{|\xi_1\xi_2\xi_3+
\eta_1\eta_2\eta_3|\leq \gamma_0N_1N_2N_{3}\}}.$$
It's easy to check that
$$\widehat{1}_{[0,\delta]}(\tau_0)=  \langle\tau_0 \rangle^{-1} .$$

Using spacetime convolution we rewrite
$$\int_{\mathbb{R}}\int_{\sum^3_{{\iota}=1}\tau_{\iota}=\tau_0,\sum^3_{{\iota}=1}
\zeta_{\iota}=0}\widehat{1}_{[0,
\delta]}(\tau_0)\widehat{u}_{1}(\lambda_1)\widehat{u}_{2}(\lambda_2)\widehat{u}_{3}
(\lambda_3)d\tau_0$$ as
\begin{align}
\int_{\mathbb{R}}\langle\tau_0\rangle^{-1}\widehat{u}_{1}*\widehat{u}_{2}*\widehat{u}_{3}
(\tau_0,0)d\tau_0.\nonumber
\end{align}
It is unfortunate that $\langle\tau_0\rangle^{-1}$ fails to be integrable. However, we can
use the logarithmic weight introduced in \cite{CKSTT08}.

Set $\omega(\tau)=1+\log^2\langle\tau\rangle$, then
$\langle\tau_0\rangle^{-1}\omega^{-1}(\tau)$ is integrable. From elementary estimate
$$\omega(\tau_1+\tau_2+\tau_3)\lesssim\omega(\tau_1)\omega(\tau_2)\omega(\tau_3),$$
we have
\begin{align}
\omega(\tau_1+\tau_2+\tau_3)\widehat{u}_{1}*\widehat{u}_{2}*\widehat{u}_{3}\lesssim
\big(\omega(\tau_1)\widehat{u}_{1}\big)*\big(\omega(\tau_2)\widehat{u}_{2}\big)*\big(
\omega(\tau_3)\widehat{u}_{3}\big).\label{omega1}
\end{align}
If $\widehat{v}_{\iota}$ denotes $\omega(\tau_{\iota})\widehat{u}_{\iota}$, one easily
verifies that
\begin{align}
\|v_{\iota}\|_{X^{0,b-}_{\delta}}\lesssim
\log^2(1+\overline{N})\|u_{\iota}\|_{X^{0,b}_{\delta}}\label{vjlog1}
\end{align}
for $b>\frac{1}{2}$, where $\overline {N}=\max_{1\leq {\iota}\leq3}\{|\zeta_{\iota}|\}$.

By \eqref{omega1}, \eqref{vjlog1} and Parseval formula, we have
\begin{align}
&\int_{\mathbb{R}}\langle\tau_0\rangle^{-1}\widehat{u}_{1}*\widehat{u}_{2}*\widehat{u}_{3}
(\tau_0,0)d\tau_0\notag \\
\lesssim &\int_{\mathbb{R}}\langle\tau_0\rangle^{-1}\omega^{-1}(\tau_0)\cdot\big
(\omega\widehat{u}_{1}\big)*
\big(\omega\widehat{u}_{2}\big)*\big(\omega\widehat{u}_{3}\big)(\tau_0,0)d\tau_0\notag \\
\lesssim &\big\|\big(\omega\widehat{u}_{1}\big)*
\big(\omega\widehat{u}_{2}\big)*\big(\omega\widehat{u}_{3}\big)\big\|_{L^{\infty}
_{\tau,\zeta}}\notag \\
\lesssim &\big\|\mathscr F^{-1}\big(\omega\widehat{u}_{1}\big)\cdot
\mathscr F^{-1}\big(\omega\widehat{u}_{2}\big)\cdot
\mathscr F^{-1}\big(\omega\widehat{u}_{3}\big)\big\|_{L^{1}(\mathbb{R}^3)}\notag \\
\lesssim &\int\int_{*}\big(\omega\widehat{u}_{1}\big)\cdot
\big(\omega\widehat{u}_{2}\big)\cdot
\big(\omega\widehat{u}_{3}\big).\nonumber
\end{align}
Thus, in order to show \eqref{Lambda(3)}, partition up into Littlewood-Paley pieces, it
suffices to show that
\begin{align}
&\bigg|\iint_{\Omega_r^{*}}\frac{\sum_{{\iota}=1}^3m^2(\zeta_{\iota})(\xi_{\iota}+
\eta_{\iota})}{m(\zeta_1)m(\zeta_2)m(\zeta_3)}
\widehat{v}_{1}(\lambda_1)\widehat{v}_{2}(\lambda_2)\widehat{v}_{3}(\lambda_3)\bigg|\notag
\\
\lesssim &\big(\gamma_0^{-\frac{1}{2}+}N^{-\frac{1}{2}+}+
N^{-\frac{1}{4}+}\big)(L_1L_2L_3)^{\frac{1}{2}}\|v_1\|_{L^{2}(\mathbb{R}^3)}
\|v_2\|_{L^{2}(\mathbb{R}^3)}
\|v_3\|_{L^{2}(\mathbb{R}^3)},\label{mv123}
\end{align}
where $v_{{\iota}}=\mathscr{F}^{-1}\omega(\tau_{\iota})\widehat{u}_{\iota}\in
L^{2}(\mathbb{R}^3)$ satisfying
$$supp \ \widehat{v}_{\iota}\subset
G_{N_{\iota},L_{\iota}}=\big\{(\tau_{\iota},\zeta_{\iota})\in \mathbb{R}\times \mathbb{R}^2
\big| \  |\zeta_{\iota}|\sim N_{\iota}, \ |\tau_{\iota}-\xi_{\iota}^3-\eta_{\iota}^3|\sim
L_{\iota} \big\},$$
and
\begin{align}
\Omega_r^{*}=\big\{&(\lambda_1,\lambda_2,\lambda_3)\in \mathbb{R}^{3\times3} \big| \
\lambda_1+\lambda_2+\lambda_3=0, \  \notag \\
&\hspace{2mm} \max_{1\leq {\iota}\leq 3}\{|\zeta_{\iota}|\}\ \gtrsim N, \ \
|\xi_1\xi_2\xi_3+\eta_1\eta_2\eta_3|<\gamma_0N_1N_2N_3\big\}.\nonumber
\end{align}

From symmetries, one can assume that $N_1\sim N_2\gtrsim N_3$, $L_1\gtrsim L_2$ and
$|\xi_1|\geq |\eta_1|$. By definition of the norms it is enough to consider functions with
non-negative Fourier transform.

\subsection{High modulation}

We first treat the high modulation case:
$$\max\{L_1,L_2,L_3\}\gtrsim\gamma_0N_1N_2N_3.$$
If $L_1=\max\{L_1,L_2,L_3\}$, using bilinear Strichartz estimate \eqref{bles4}, one has
\begin{align}
&\bigg|\iint_{\Omega_r^{*}}\frac{\sum_{{\iota}=1}^3m^2(\zeta_{\iota})(\xi_{\iota}+
\eta_{\iota})}{m(\zeta_1)m(\zeta_2)m(\zeta_3)}
\widehat{v}_{1}(\lambda_1)\widehat{v}_{2}(\lambda_2)\widehat{v}_{3}(\lambda_3)\bigg|\notag
\\
\lesssim
&\frac{m^2(\zeta_3)N_3}{m^2(\zeta_1)m(\zeta_3)}\|v_1\|_{L^{2}(\mathbb{R}^3)}\|v_2v_3\|_{L^{2}
(\mathbb{R}^3)}\notag \\
\lesssim &\frac{m(\zeta_3)N_3}{m^2(\zeta_1)}\frac{N_3^{\frac{1}{2}}}{N_1}L_1^{-\frac{1}{2}}
(L_1L_2L_3)^{\frac{1}{2}}\|v_1\|_{L^{2}}\|v_2\|_{L^{2}}
\|v_3\|_{L^{2}}\notag \\
\lesssim
&\frac{N_3^{\frac{3}{2}}}{N_1}\big(\gamma_0N_1N_2N_3\big)^{-\frac{1}{2}}(L_1L_2L_3)^
{\frac{1}{2}}\|v_1\|_{L^{2}}\|v_2\|_{L^{2}}
\|v_3\|_{L^{2}}\notag \\
\lesssim
&\gamma_0^{-\frac{1}{2}}N^{-1}(L_1L_2L_3)^{\frac{1}{2}}\|v_1\|_{L^{2}}\|v_2\|_{L^{2}}
\|v_3\|_{L^{2}}.\nonumber
\end{align}
If $L_3=\max\{L_1,L_2,L_3\}$, using \eqref{estimate16}, we get
\begin{align}
&\bigg|\iint_{\Omega_r^{*}}\frac{\sum_{{\iota}=1}^3m^2(\zeta_{\iota})(\xi_{\iota}+
\eta_{\iota})}{m(\zeta_1)m(\zeta_2)m(\zeta_3)}
\widehat{v}_{1}(\lambda_1)\widehat{v}_{2}(\lambda_2)\widehat{v}_{3}(\lambda_3)\bigg|\notag
\\
\lesssim &\frac{m^2(\zeta_3)N_3}{m^2(\zeta_1)m(\zeta_3)}\|v_1\|_{L^{4}(\mathbb{R}^3)}
\|v_2\|_{L^{4}(\mathbb{R}^3)}\|v_3\|_{L^{2}(\mathbb{R}^3)}\notag \\
\lesssim
&\frac{m(\zeta_3)N_3}{m^2(\zeta_1)}(L_1L_2)^{\frac{5}{12}}\|v_1\|_{L^{2}}\|v_2\|_{L^{2}}
\|v_3\|_{L^{2}}\notag \\
\lesssim
&N_3\big(\gamma_0N_1N_2N_3\big)^{-\frac{1}{2}}(L_1L_2L_3)^{\frac{1}{2}}\|v_1\|_{L^{2}}
\|v_2\|_{L^{2}}
\|v_3\|_{L^{2}}\notag \\
\lesssim
&\gamma_0^{-\frac{1}{2}}N^{-\frac{1}{2}}(L_1L_2L_3)^{\frac{1}{2}}\|v_1\|_{L^{2}}\|v_2\|
_{L^{2}}
\|v_3\|_{L^{2}}.\nonumber
\end{align}
This concludes \eqref{mv123}.

\subsection{Low modulation}

Next we focus our attention on the low modulation case which will be divided into high-low
interactions interactions and high-high interactions. In each subcase, we need to consider
non-parallel interactions and parallel interactions respectively.

Low modulation means that
$$\max\{L_1,L_2,L_3\}\ll\gamma_0N_1N_2N_3.$$
For the high modulation case, we use Strichartz estimates to obtain some decay
$L^{-\frac{1}{2}}\sim N^{-\frac{3}{2}}$ which offsets the increment from $\mathcal{M}_3$.
However, for the low modulation case we can not obtain decay from $L$ any more. So we should
thoroughly consider the quantity $\int_{*}fgh$. In fact, the quantity $\int_{*}fgh$ can be
treated by applying the convolution estimate on hypersufaces, which is called the nonlinear
Loomis-Whitney inequality. The nonlinear Loomis-Whitney inequality is an important tool for
low regularity problems. When using this inequality, we will find that the upper bound is
related to a transversality condition
$$\big|\xi_1\eta_2-\xi_2\eta_1\big|\cdot\big|\xi_1\eta_2+\xi_2\eta_1+2(\xi_1\eta_1+
\xi_2\eta_2)
\big|.$$

As the former term $\big|\xi_1\eta_2-\xi_2\eta_1\big|$ is comparable to
$N_1^2\big|\sin\angle\big((\xi_1,\eta_1),(\xi_2,\eta_2)\big)\big|$, we call them parallel
interactions when $\big|\xi_1\eta_2-\xi_2\eta_1\big|\ll 1$. On one hand, we use angular
dyadic decomposition for parallel interactions. On the other hand, for non-parallel
interactions we utilize Whitney type decomposition and almost orthogonality. These ideas are
originally from Kinoshita \cite{Kinoshita19}. But there are also great difference in technique. For
instance, we find that $\big|\xi_1\eta_2-\xi_2\eta_1\big|$ and
$\big|\xi_1\eta_2+\xi_2\eta_1+2(\xi_1\eta_1+\xi_2\eta_2)\big|$ can not be very small at the
same time in our problem. Moreover, we can easily get almost orthogonality after making
decomposition.

In order to take advantage of the nonlinear Loomis--Whitney inequality, we need to decompose
$\mathbb{R}^2$ into square tiles.
\begin{definition}[see \cite{Kinoshita19} Def. 1]Let $A\gg1$ be a dyadic number and
$k=(k(1),k(2))\in \mathbb{Z}^2$. We define square-tiles $\{\mathcal{T}^A_k\}_{k\in
\mathbb{Z}^2}$ whose side length is $A^{-1}N_1$ and prisms
$\{\widetilde{\mathcal{T}}^A_k\}_{k\in \mathbb{Z}^2}$ as follows:
\begin{align}
&\mathcal{T}^A_k:=\big\{\zeta\in\mathbb{R}^2 \ \big| \ \zeta\in
A^{-1}N_1\big([k_{(1)},k_{(1)}+1)\times[k_{(2)},k_{(2)}+1)\big)\big\}, \notag \\
&\widetilde{\mathcal{T}}^A_k:=\mathbb{R}\times\mathcal{T}^A_k. \nonumber
\end{align}
\end{definition}
Let's consider the high-low interactions case.
\\{\bf Case 1} \hspace{2mm} $N_1\sim N_2\gg N_3$ and $|\xi_1|\geq |\eta_1|$

If $|\xi_1|\sim |\eta_1|$, as $|\xi_3|\lesssim N_3\ll N_1\sim |\xi_1|$ , $|\eta_3|\lesssim
N_3\ll N_1\sim |\eta_1|$ and $\zeta_1+\zeta_2+\zeta_3=0$, we have $|\xi_1|\sim |\xi_2|$ and
$|\eta_1|\sim |\eta_2|$. We claim that $|\xi_3|\sim|\eta_3|$. Otherwise,
$|\xi_3|\gg|\eta_3|$ or $|\xi_3|\ll |\eta_3|$, then
$|\xi_1\xi_2\xi_3+\eta_1\eta_2\eta_3|\sim N_1^2N_3\gg\gamma_0N_1^2N_3$, which contradicts to
the assumption $|\xi_1\xi_2\xi_3+\eta_1\eta_2\eta_3|\leq\gamma_0N_1^2N_3$.

$\zeta_1+\zeta_2+\zeta_3=0$ and $|\xi_3|\sim|\eta_3|\ll|\xi_1|\sim|\eta_1|$ help imply
$$
\big|\xi_1\eta_2+\xi_2\eta_1+2(\xi_1\eta_1+\xi_2\eta_2)\big|=\big|\xi_1\eta_3
+\xi_3\eta_1+2(\xi_1\eta_1+\xi_3\eta_3)\big|,$$
and
$$\big|\xi_1\eta_3+\xi_3\eta_1+2(\xi_1\eta_1+\xi_3\eta_3)\big|\sim |\xi_1\eta_1|\sim
N_1^2,$$
hence
$$\big|\xi_1\eta_2+\xi_2\eta_1+2(\xi_1\eta_1+\xi_2\eta_2)\big|\sim N_1^2.$$

So, under the high-low frequencies case, we consider parallel interactions.
\\{\bf Case 1.1} \hspace{2mm} $|\xi_1|\sim |\eta_1|\sim |\xi_2|\sim |\eta_2|\gg |\xi_3|\sim
|\eta_3|$

As $\xi_3+\eta_3$ can be very small, we will use angular decomposition (see
\eqref{angulardecomp} for ${\mathcal{D}}_j^A$).

Let $\angle(\zeta_1,\zeta_2)\in (0,\frac{\pi}{2}]$ denote the angle between the line spanned
by $\zeta_1,\zeta_2\in \mathbb{R}^2$. For dyadic numbers $64\leq A\leq M$ we consider the
following angular decomposition:
\begin{align}
\mathbb{R}^2\times\mathbb{R}^2=&\big\{ \angle(\zeta_1,\zeta_2)\leq\frac{16\pi}{M} \big\} \
\cup\bigcup\limits_{64\leq A\leq M}\big\{ \frac{16\pi}{A}\leq
\angle(\zeta_1,\zeta_2)\leq\frac{32\pi}{A} \big\} \notag \\
=&\bigcup\limits_{\substack{0\leq j_1,j_2\leq M-1 \\ |j_1-j_2|\leq16} }{\mathcal
D}^{M}_{j_1}\times {\mathcal D}^{M}_{j_2}\ \cup
\bigcup\limits_{64\leq A\leq M}\bigcup\limits_{\substack{0\leq j_1,j_2\leq A-1 \\
16\leq|j_1-j_2|\leq32} }{\mathcal D}^{A}_{j_1}\times {\mathcal D}^{A}_{j_2}.\nonumber
\end{align}

Let $\mathcal{I}_1,\mathcal{I}_2\in\mathbb{R}^2\times\mathbb{R}^2$ be defined as
$$\mathcal{I}_1={\mathcal D}^{2^{11}}_{2^{9}\times3}\times{\mathcal
D}^{2^{11}}_{2^{9}\times3},\hspace{5mm} \widetilde{\mathcal{I}}_1={\widetilde{\mathcal
D}}^{2^{11}}_{2^{9}\times3}\times{\widetilde{\mathcal D}}^{2^{11}}_{2^{9}\times3},$$
$$\mathcal{I}_2=(\mathbb{R}^2\times\mathbb{R}^2)\setminus
\mathcal{I}_1,\hspace{5mm}\widetilde{\mathcal{I}}_2=(\mathbb{R}^3\times\mathbb{R}^3)\setminus
\widetilde{\mathcal{I}}_1.$$
Note that
$${\mathcal
D}^{2^{11}}_{2^{9}\times3}=\Big\{\big(|(\xi,\eta)|\cos\theta,|(\xi,\eta)|\sin\theta)\big)
\in\mathbb{R}^2
\ \big| \
\min\big(|\theta-\frac{3\pi}{4}|,|\theta+\frac{\pi}{4}|\big)\leq\frac{\pi}{2^{10}}\Big\} .$$
Next we treat $(\xi_1,\eta_1)\times(\xi_2,\eta_2)\notin \mathcal{I}_1$ and
$(\xi_1,\eta_1)\times(\xi_2,\eta_2)\in \mathcal{I}_1$ respectively.
\begin{lemma}\label{mathcalI1a0}
Let $A\geq 2^{25}$ be dyadic, $\lambda_{\iota}=(\tau_{\iota},\xi_{\iota},\eta_{\iota})\in
G_{N_{\iota},L_{\iota}}\cap {\widetilde{\mathcal D}}^{A}_{j_{\iota}}$ $({\iota}=1,2)$ and
$\lambda_1+\lambda_2\in G_{N_3,L_3}$. Assume that $\max_{1\leq
{\iota}\leq3}\{L_{\iota}\}\lesssim A^{-1}N_1^3$, $|j_1-j_2|\leq 32$ and
$(\xi_1,\eta_1)\times(\xi_2,\eta_2)\notin \mathcal{I}_1$.
In the assumption of Case 1.1, we have $N_3\lesssim A^{-1}N_1$.
\end{lemma}
{\bf Proof.} We write
$(\xi_{\iota},\eta_{\iota})=r_{\iota}(\cos\theta_{\iota},\sin\theta_{\iota})$ for
${\iota}=1,2$. Since $(\xi_1,\eta_1)\times(\xi_2,\eta_2)\notin \mathcal{I}_1$, without loss
of generality, we shall assume that $(\xi_1,\eta_1)\notin {\mathcal
D}^{2^{11}}_{2^{9}\times3}$ which implies
$$|\cos\theta_1+\sin\theta_1|=\sqrt2|\sin(\theta_1+\frac{\pi}{4})|>2^{-11}\pi.$$

From the assumption of Case 1.1, we have $\xi_1\xi_2<0$ and $\eta_1\eta_2<0$. Thus we deduce
from $|j_1-j_2|\leq 32$ that
$$|(\cos\theta_1,\sin\theta_1)+(\cos\theta_2,\sin\theta_2)|\leq2^8A^{-1}.$$

As we observe that
\begin{align}
&|\xi_1\xi_2(\xi_1+\xi_2)+\eta_1\eta_2(\eta_1+\eta_2)|\notag \\
=&\big|r_1r_2(r_1-r_2)(\cos^3\theta_1+\sin^3\theta_1)-r_1r_2(r_1-r_2)\cos^2\theta_1
(\cos\theta_1+\cos\theta_2)\notag \\
&-r_1r_2(r_1-r_2)\sin^2\theta_1(\sin\theta_1+\sin\theta_2)-r_1r^2_2\cos\theta_1\cos\theta_2
(\cos\theta_1+\cos\theta_2)\notag \\
&-r_1r^2_2\sin\theta_1\sin\theta_2(\sin\theta_1+\sin\theta_2)\big|\notag \\
\geq&2^{-13}r_1r_2|r_1-r_2|-2^{9}A^{-1}r_1r^2_2, \nonumber
\end{align}
and
\begin{align}
 A^{-1}N_1^3\gtrsim \max_{1\leq {\iota}\leq3}\{L_{\iota}\}\gtrsim
\bigg|\sum^3_{{\iota}=1}\tau_{\iota}-\xi^3_{\iota}-\eta^3_{\iota}\bigg|\gtrsim
|\xi_1\xi_2(\xi_1+\xi_2)+\eta_1\eta_2(\eta_1+\eta_2)|,
 \nonumber
\end{align}
hence we get
$$|r_1-r_2|\lesssim A^{-1}N_1.$$
Thus
\begin{align}
 N_3=&|(\xi_1+\xi_2,\eta_1+\eta_2)|\notag \\
 \leq&|(r_1\cos\theta_1-r_2\cos\theta_1,r_1\sin\theta_1-r_2\sin\theta_1)|\notag \\
 & \ +r_2|(\cos\theta_1,\sin\theta_1)+(\cos\theta_2,\sin\theta_2)|\notag \\
 \leq& |r_1-r_2|+2^8A^{-1}N_1\notag \\
 \lesssim& A^{-1}N_1. \nonumber
\end{align}
\begin{lemma}\label{mathcalI1a}
Let $A\geq 2^{25}$ be dyadic, $|j_1-j_2|\leq 32$ and
$(\xi_1,\eta_1)\times(\xi_2,\eta_2)\notin \mathcal{I}_1$.
In the assumption of Case 1.1, we have
\begin{align}
 &\big\|\chi_{G_{N_3,L_3}}(\lambda_3)\int \widehat{v}_1\big|_{{\widetilde{\mathcal
D}}^{A}_{j_1}}(\lambda_1) \widehat{v}_2\big|_{{\widetilde{\mathcal
D}}^{A}_{j_2}}(\lambda_1+\lambda_3)d\lambda_1\big\|_{L^2_{\lambda_3}}\notag \\
 \lesssim &
(AN_1)^{-\frac{1}{2}}(L_1L_2)^{\frac{1}{2}}\big\|\widehat{v}_{1}|_{{\widetilde{\mathcal
D}}^{A}_{j_1}}\big\|_{L^{2}}\big\|\widehat{v}_{2}|_{{\widetilde{\mathcal
D}}^{A}_{j_2}}\big\|_{L^{2}},\label{case11bles1}
\end{align}
\begin{align}
 &\big\|\chi_{G_{N_2,L_2}}(\lambda_2)\int \widehat{v}_1\big|_{{\widetilde{\mathcal
D}}^{A}_{j_1}}(\lambda_1)
\widehat{v}_3(\lambda_1+\lambda_2)d\lambda_1\big\|_{L^2_{\lambda_2}}\notag \\
 \lesssim &
(AN_1)^{-\frac{1}{2}}(L_1L_3)^{\frac{1}{2}}\big\|\widehat{v}_{1}|_{{\widetilde{\mathcal
D}}^{A}_{j_1}}\big\|_{L^{2}}\big\|\widehat{v}_{3}\big\|_{L^{2}},\label{case11bles2}
\end{align}
and
\begin{align}
 &\big\|\chi_{G_{N_1,L_1}}(\lambda_1)\int \widehat{v}_2\big|_{{\widetilde{\mathcal
D}}^{A}_{j_2}}(\lambda_2)
\widehat{v}_3(\lambda_1+\lambda_2)d\lambda_2\big\|_{L^2_{\lambda_1}}\notag \\
 \lesssim &
(AN_1)^{-\frac{1}{2}}(L_2L_3)^{\frac{1}{2}}\big\|\widehat{v}_{2}|_{{\widetilde{\mathcal
D}}^{A}_{j_2}}\big\|_{L^{2}}\big\|\widehat{v}_{3}\big\|_{L^{2}},\label{case11bles3}
\end{align}
where $supp \ \widehat{v}_k\subset G_{N_k,L_k}$ for $k=1,2,3$.
\end{lemma}
{\bf Proof.} By the assumption of Case 1.1,
$$|\xi_1|\sim N_1,\ \ \ |\xi_2|\sim N_2.$$
From the definition of ${\widetilde{\mathcal D}}^{A}_{j}$,
$$\Big|\left\{\eta_{\iota}\big|\lambda_{\iota}=(\tau_{\iota},\xi_{\iota},\eta_{\iota})\in
{\widetilde{\mathcal D}}^{A}_{j_{\iota}} \right\}\Big|\lesssim A^{-1}N_1$$
for ${\iota}=1,2$.

Hence, from Lemma \ref{X Lema3}, we know this lemma holds true.
\begin{lemma}\label{mathcalI1b}
Let $ 2^{25}\leq A\lesssim \frac{N_1}{N_3}$ be dyadic, $16\leq|j_1-j_2|\leq 32$ and
$(\xi_1,\eta_1)\times(\xi_2,\eta_2)\notin \mathcal{I}_1$.
In the assumption of Case 1.1, we have
\begin{align}
&\bigg|\iint_{\Omega_r^{*}}\widehat{v}_{1}|_{{\widetilde{\mathcal
D}}^{A}_{j_1}}(\lambda_1)\widehat{v}_{2}|_{{\widetilde{\mathcal
D}}^{A}_{j_2}}(\lambda_2)\widehat{v}_{3}(\lambda_3)\bigg| \notag \\
 \lesssim &
A^{\frac{1}{2}}N_1^{-2}(L_1L_2L_3)^{\frac{1}{2}}\big\|\widehat{v}_{1}|_{{\widetilde{\mathcal
D}}^{A}_{j_1}}\big\|_{L^{2}}\big\|\widehat{v}_{2}|_{{\widetilde{\mathcal
D}}^{A}_{j_2}}\big\|_{L^{2}}\big\|\widehat{v}_{3}\big\|_{L^{2}},\label{case11b1}
\end{align}
where $supp \ \widehat{v}_{\iota}\subset G_{N_{\iota},L_{\iota}}$ for ${\iota}=1,2,3$.
\end{lemma}
{\bf Proof.} To prove \eqref{case11b1}, the strategy is using the nonlinear Loomis-Whitney
inequality (see \cite{BHHT09} and \cite{Kinoshita19}). We take advantage of square prisms
decomposition so that the nonlinear Loomis-Whitney inequality can be applied.

Let $A'=2^{30}A$. Assume that $f,g$ and $h$ are functions satisfying
$$supp \ f\subset G_{N_1,L_1}\cap\widetilde{\mathcal{T}}^{A'}_{k_1}, \hspace{2mm}supp \
g\subset G_{N_2,L_2}\cap\widetilde{\mathcal{T}}^{A'}_{k_2},\hspace{2mm} supp \ h\subset
G_{N_3,L_3}\cap\widetilde{\mathcal{T}}^{A'}_{k_3}, $$
where $\big(\mathcal{T}^{A'}_{k_1}\times\mathcal{T}^{A'}_{k_2}\big)\cup\big(\mathcal
{D}^{A}_{j_1}\times\mathcal {D}^{A}_{j_2}\big)\not=\emptyset$.

Since $N_3\lesssim A^{-1}N_1$, we reduce to showing that
\begin{align}
\bigg|\int_{*}f(\lambda_1)g(\lambda_2)h(\lambda_3)\bigg|
 \lesssim
A^{\frac{1}{2}}N_1^{-2}(L_1L_2L_3)^{\frac{1}{2}}\|f\|_{L^{2}}\|g\|_{L^{2}}\|h\|_{L^{2}}.
\nonumber
\end{align}

First we change variables $c_1=\tau_1-\xi_1^3-\eta_1^3$ and $c_2=\tau_2-\xi_2^3-\eta_2^3$.
Then, by decomposing $h$ into $L_3$ pieces and applying the Cauchy-Schwarz inequality, it
suffices to prove that
\begin{align}
&\bigg|\int
f\big(\phi_{c_1}(\xi_1,\eta_1)\big)g\big(\phi_{c_2}(\xi_2,\eta_2)\big)h\big(\phi_{c_1}
(\xi_1,\eta_1)
+\phi_{c_2}(\xi_2,\eta_2)\big)d\xi_1d\eta_1d\xi_2d\eta_2\bigg|  \notag \\
 \lesssim &
A^{\frac{1}{2}}N_1^{-2}\|f\circ\phi_{c_1}\|_{L^{2}_{\xi\eta}}\|g\circ\phi_{c_2}\|_{L^{2}
_{\xi\eta}}\|h\|_{L^{2}_{\tau\xi\eta}},\label{case11b2}
\end{align}
where $h(\tau,\xi,\eta)$ is supported in $c_0\leq\tau-\xi^3-\eta^3\leq c_0+1$ and
$$\phi_{c_{\iota}}(\xi_{\iota},\eta_{\iota})=(\xi_{\iota}^3+\eta_{\iota}^3+c_{\iota},
\xi_{\iota},\eta_{\iota})\  \  \ \  \  for \  \ {\iota}=1,2.$$

We use the scaling $(\tau,\xi,\eta)\rightarrow(N_1^3\tau,N_1\xi,N_1\eta)$ to define
 \begin{align}
 \tilde f(\tau_1,\xi_1,\eta_1)=f(N_1^3\tau_1,N_1\xi_1,N_1\eta_1),\notag \\
 \tilde g(\tau_2,\xi_2,\eta_2)=g(N_1^3\tau_2,N_1\xi_2,N_1\eta_2),\notag \\
  \tilde h(\tau_3,\xi_3,\eta_3)=h(N_1^3\tau_3,N_1\xi_3,N_1\eta_3).\nonumber
\end{align}

Setting $\tilde c_{\iota}=N_1^{-3}c_{\iota}$, then \eqref{case11b2} reduces to
\begin{align}
&\bigg|\int \tilde f\big(\phi_{\tilde c_1}(\xi_1,\eta_1)\big)\tilde g\big(\phi_{\tilde
c_2}(\xi_2,\eta_2)\big)\tilde h\big(\phi_{\tilde c_1}(\xi_1,\eta_1)
+\phi_{\tilde c_2}(\xi_2,\eta_2)\big)d\xi_1d\eta_1d\xi_2d\eta_2\bigg|  \notag \\
 \lesssim & A^{\frac{1}{2}}N_1^{-\frac{3}{2}}\|\tilde f\circ\phi_{\tilde
c_1}\|_{L^{2}_{\xi\eta}}\|\tilde g\circ\phi_{\tilde c_2}\|_{L^{2}_{\xi\eta}}\|\tilde
h\|_{L^{2}_{\tau\xi\eta}},\label{case11b2666}
\end{align}
where $
supp \ \tilde  f\subset\widetilde{\mathcal{T}}^{N_1A'}_{\widetilde  k_1}, \hspace{2mm}supp \
\tilde  g\subset\widetilde{\mathcal{T}}^{N_1A'}_{\widetilde  k_2}$, and $\tilde  h$ is
supported in a neighbourhood of size $N_1^{-3}$ of the surface
$S_3:=\big\{(\tau_3,\xi_3,\eta_3)\in
\mathbb{R}^3\big|\tau_3-\xi^3_3-\eta^3_3=c_3N_1^{-3}\big\}$, precisely, $supp \ \tilde
h\subset S_3(N_1^{-3})$,
where
$$S_3(N_1^{-3}):=\bigg\{(\tau_3,\xi_3,\eta_3)\in \widetilde{\mathcal{T}}^{N_1A'}_{\widetilde
k_3}\ \bigg|\
\xi^3_3+\eta^3_3+\frac{c_0}{N_1^3}\leq\tau_3\leq\xi^3_3+\eta^3_3+\frac{c_0+1}{N_1^3}\bigg\},$$
$c_3\in [c_0,c_0+1]$ and $\widetilde  k_{\iota}=k_{\iota}/N_1$ for $\iota=1,2,3$.

By density and duality it is enough to estimate
\begin{align}
\big\|\tilde  f|_{S_1}*\tilde  g|_{S_2}\big\|_{L^2(S_3(N_1^{-3}))}
 \lesssim  A^{\frac{1}{2}}N_1^{-\frac{3}{2}}\|\tilde  f\|
 _{L^2(S_1)}\|\tilde  g\|_{L^2(S_2)}\label{case11b3}
\end{align}
for continuous $\tilde  f, \tilde  g$, where $S_{\iota}$ are parametrized by $\phi_{\tilde
c_{\iota}}$ for ${\iota}=1,2$.
Moreover, it will suffice to show
\begin{align}
\big\|\tilde  f|_{S_1}*\tilde  g|_{S_2}\big\|_{L^2(S_3)}
 \lesssim  A^{\frac{1}{2}}\|\tilde  f\|_{L^2(S_1)}\|\tilde  g\|_{L^2(S_2)}.\label{case11b4}
\end{align}

As $S_{\iota}$ are restricted in $\widetilde{\mathcal{T}}^{N_1A'}_{\widetilde  k_{\iota}}$,
hence
\begin{align}
diam(S_{\iota})\leq 2^{-25}A^{-1}\label{case11bdiam}
\end{align}
for ${\iota}=1,2,3$.
For any $\lambda_{\iota}=\phi_{\tilde c_{\iota}}(\xi_{\iota},\eta_{\iota})\in S_{\iota}$,
the unit normals $\mathfrak {n}_{\iota}$ on $\lambda_{\iota}$ are
$$\mathfrak
{n}_j(\lambda_{\iota})=\frac{1}{1+9\xi^4_{\iota}+9\eta^4_{\iota}}(-1,3\xi^2_{\iota},
3\eta^2_{\iota})$$
for ${\iota}=1,2,3$. Thus, the surfaces $S_1,S_2,S_3$ satisfy the H\"older condition
\begin{align}
\sup_{\lambda_{\iota},\widetilde{\lambda}_{\iota}\in S_{\iota}^*}\frac{|\mathfrak
{n}_{\iota}(\lambda_{\iota})-\mathfrak
{n}_{\iota}(\widetilde{\lambda}_{\iota})|}{|\lambda_{\iota}-\widetilde{\lambda}_{\iota}|}
+\frac{|\mathfrak{n}_{\iota}(\lambda_{\iota})(\lambda_{\iota}-\widetilde{\lambda}_{\iota})|}
{|\lambda_{\iota}
-\widetilde{\lambda}_{\iota}|^2}\lesssim 1.\nonumber
\end{align}
As the region we need to consider is $\big(S_1+S_2\big)\cap S_3 $, we can assume that there
exist $(\xi'_{\iota},\eta'_{\iota})$ such that
$$(\xi'_1,\eta'_1)+(\xi'_2,\eta'_2)=(\xi'_3,\eta'_3),$$
$$\lambda'_{\iota}=\phi_{\tilde c_{\iota}}(\xi'_{\iota},\eta'_{\iota})\in S_{\iota}
\hspace{5mm}{\iota}=1,2,3.$$
From \eqref{case11bdiam}, we know that
$$|\mathfrak {n}_{\iota}(\lambda_{\iota})-\mathfrak
{n}_{\iota}(\lambda'_{\iota})|\leq2^{-20}A^{-1}$$
for ${\iota}=1,2,3$.

As $16\leq|j_1-j_2|\leq 32$, hence $\big|\xi'_1\eta'_2-\xi'_2\eta'_1\big|\sim
\sin\angle\big((\xi_1,\eta_1),(\xi_2,\eta_2)\big)\sim A^{-1}.$ And we can bound
\begin{align}
&\big|\xi'_1\eta'_2+\xi'_2\eta'_1+2(\xi'_1\eta'_1+\xi'_2\eta'_2)\big|\notag \\
\sim&N_1^{-2}\big|\xi_1\eta_2+\xi_2\eta_1+2(\xi_1\eta_1+\xi_2\eta_2)\big|\sim 1.\nonumber
\end{align}
Therefore,
\begin{align}
|det\mathcal {N}(\lambda_1,\lambda_2,\lambda_3)|&\gtrsim |det\mathcal
{N}(\lambda'_1,\lambda'_2,\lambda'_3)|\notag \\
&\gtrsim \frac{1}{\prod^3_{{\iota}=1}\langle (\xi'_{\iota},\eta'_{\iota})\rangle
^2}\Bigg|det
\left(
  \begin{array}{ccc}
    -1&3\xi'^2_1&3\eta'^2_1\\
    -1&3\xi'^2_2&3\eta'^2_2\\
    -1&3\xi'^2_3&3\eta'^2_3\\
  \end{array}
\right)
\Bigg|\notag \\
&\gtrsim
\big|\xi'_1\eta'_2-\xi'_2\eta'_1\big|\cdot\big|\xi'_1\eta'_2+\xi'_2\eta'_1+2(\xi'_1\eta'_1
+\xi'_2\eta'_2)\big|\notag \\
&\gtrsim A^{-1},\nonumber
\end{align}
which implies \eqref{case11b4} by Lemma \ref{LW1}. This concludes the proof of Lemma
\ref{mathcalI1b}.
\begin{proposition}\label{mathcalI1}
Assume that $(\xi_1,\eta_1)\times(\xi_2,\eta_2)\notin \mathcal{I}_1$ , then we have
\begin{align}
\bigg|\iint_{\Omega_r^{*}}\widehat{v}_{1}(\lambda_1)\widehat{v}_{2}(\lambda_2)
\widehat{v}_{3}(\lambda_3)\bigg|
\lesssim
N_1^{-\frac{3}{2}}N_3^{-\frac{1}{2}}(L_1L_2L_3)^{\frac{1}{2}}\|\widehat{v}_{1}\|_{L^{2}}\|
\widehat{v}_{2}\|_{L^{2}}
\|\widehat{v}_{3}\|_{L^{2}}.\nonumber
\end{align}
\end{proposition}
{\bf Proof.} First we claim that $A\lesssim\frac{N_1}{N_3}$ or
$\gamma_0^{-1}\frac{N_1}{N_3}\ll A$. Otherwise, $\frac{N_1}{N_3}\ll
A\lesssim\gamma_0^{-1}\frac{N_1}{N_3}$, thus we have $A^{-1}N_1\geq \gamma_0N_3$ and $N_3\gg
A^{-1}N_1$ which imply
$$\max_{1\leq {\iota}\leq3}\{L_{\iota}\}\gg A^{-1}N_1^3\gtrsim\gamma_0N^2_1N_3$$
by Lemma \ref{mathcalI1a0}. However, this contradicts to the low modulation assumption.

If $\gamma_0^{-1}\frac{N_1}{N_3}\ll A$, then $A^{-1}N_1\ll\gamma_0N_3$ which implies
\begin{align}
\max_{1\leq {\iota}\leq3}\{L_{\iota}\}\gg A^{-1}N_1^3.  \label{AAAA}
\end{align}
In fact, by Lemma \ref{mathcalI1a0}, $\max_{1\leq {\iota}\leq3}\{L_{\iota}\}\lesssim
A^{-1}N_1^3$ will give
$$N_3\lesssim A^{-1}N_1\ll\gamma_0N_3\ll N_3.$$

According to angular decomposition, Lemma \ref{mathcalI1b} , Lemma \ref{mathcalI1a} and
\eqref{AAAA} , one has
\begin{align}
&\bigg|\iint_{\Omega_r^{*}}\widehat{v}_{1}(\lambda_1)\widehat{v}_{2}(\lambda_2)
\widehat{v}_{3}(\lambda_3)\bigg|\notag \\
\lesssim&\sum_{2^{25}\leq A\lesssim\frac{N_1}{N_3}}\sum\limits_{\substack{0\leq j_1,j_2\leq
A-1 \\ 16\leq|j_1-j_2|\leq32} }\bigg|\int_{*}\widehat{v}_{1}|_{{\widetilde{\mathcal
D}}^{A}_{j_1}}(\lambda_1)\widehat{v}_{2}|_{{\widetilde{\mathcal
D}}^{A}_{j_2}}(\lambda_2)\widehat{v}_{3}(\lambda_3)\bigg|\notag \\
&\ \ +\sum_{\gamma_0^{-1}\frac{N_1}{N_3}\ll A}\sum\limits_{\substack{0\leq j_1,j_2\leq A-1
\\ 16\leq|j_1-j_2|\leq32} }\bigg|\int_{*}\widehat{v}_{1}|_{{\widetilde{\mathcal
D}}^{A}_{j_1}}(\lambda_1)\widehat{v}_{2}|_{{\widetilde{\mathcal
D}}^{A}_{j_2}}(\lambda_2)\widehat{v}_{3}(\lambda_3)\bigg|\notag \\
\lesssim&\sum_{2^{25}\leq A\lesssim\frac{N_1}{N_3}}\sum_{j_1\sim
j_2}A^{\frac{1}{2}}N_1^{-2}(L_1L_2L_3)^{\frac{1}{2}}\big\|\widehat{v}_{1}|_{{\widetilde
{\mathcal D}}^{A}_{j_1}}\big\|_{L^{2}}\big\|\widehat{v}_{2}|_{{\widetilde{\mathcal
D}}^{A}_{j_2}}\big\|_{L^{2}}\big\|\widehat{v}_{3}\big\|_{L^{2}}\notag \\
&\ \ +\sum_{\gamma_0^{-1}\frac{N_1}{N_3}\ll A}\sum_{j_1\sim
j_2}A^{-\frac{1}{2}}N_1^{-\frac{1}{2}}A^{\frac{1}{2}-}N_1^{-\frac{3}{2}+}(L_1L_2L_3)
^{\frac{1}{2}}\big\|\widehat{v}_{1}|_{{\widetilde{\mathcal
D}}^{A}_{j_1}}\big\|_{L^{2}}\big\|\widehat{v}_{2}|_{{\widetilde{\mathcal
D}}^{A}_{j_2}}\big\|_{L^{2}}\big\|\widehat{v}_{3}\big\|_{L^{2}}\notag \\
\lesssim&\big(N_1^{-\frac{3}{2}}N_3^{-\frac{1}{2}}+\gamma_0^{0+}N_1^{-2}N_3^{0+}\big)
(L_1L_2L_3)^{\frac{1}{2}}\|\widehat{v}_{1}\|_{L^{2}}\|\widehat{v}_{2}\|_{L^{2}}\|
\widehat{v}_{3}\|_{L^{2}}.\nonumber
\end{align}
We finish the proof.

Now let's consider the situation when $(\xi_1,\eta_1)\times(\xi_2,\eta_2)\in \mathcal{I}_1$.
The low frequency effect ( see Lemma \ref{mathcalI1a0} ) is invalid near the line
$\xi+\eta=0$, however the derivative loss $\xi_3+\eta_3$ is small. To be specific, we divide
$\mathcal{I}_1$ into tiny dyadic pieces in which angular decomposition are applied again.
For low angular frequency, we use the nonlinear Loomis-Whitney inequality. For high angular
frequency, we take advantage of bilinear Strichartz estimates. Further more, we also obtain
smallness on the line $\xi+\eta=0$.

We make a decomposition of $\mathcal{I}_1$ so that $\xi+\eta$ is controllable in each area.
\begin{definition}[see \cite{Kinoshita19} Def. 6]Let $M\geq 2^{11}$ be a dyadic number.
Define
\begin{align}
&\mathcal{I}^M_1:=\big({\mathcal D}^{M}_{\frac{3}{4}M}\times{\mathcal
D}^{M}_{\frac{3}{4}M}\big)\setminus \big({\mathcal D}^{2M}_{\frac{3}{2}M}\times{\mathcal
D}^{2M}_{\frac{3}{2}M}\big), \notag \\
&\widetilde{\mathcal{I}}^M_1:=\big(\widetilde{{\mathcal
D}}^{M}_{\frac{3}{4}M}\times\widetilde{{\mathcal D}}^{M}_{\frac{3}{4}M}\big)\setminus
\big(\widetilde{{\mathcal D}}^{2M}_{\frac{3}{2}M}\times\widetilde{{\mathcal
D}}^{2M}_{\frac{3}{2}M}\big). \nonumber
\end{align}
\end{definition}
It is easy to see that
\begin{align}
&\hspace{5mm}(r\cos\theta,r\sin\theta)\in{\mathcal D}^{M}_{\frac{3}{4}M}\setminus{\mathcal
D}^{2M}_{\frac{3}{2}M}\ \notag \\
\Longleftrightarrow \
&M^{-1}\pi\leq\min\big\{\big|\theta-\frac{3}{4}\pi\big|,\big|\theta+\frac{1}{4}\pi\big|\big\}
\leq 2M^{-1}\pi,\nonumber
\end{align}
and
\begin{align}
\mathcal{I}_1=\bigcup\limits_{2^{11}\leq M\leq M_0}\mathcal{I}^M_1\ \cup\bigg({\mathcal
D}^{M_0}_{\frac{3}{4}M_0}\times{\mathcal D}^{M_0}_{\frac{3}{4}M_0}\bigg).\nonumber
\end{align}
Moreover, $|\xi_3+\eta_3|\lesssim M^{-1}N_1$ if $(\xi_1,\eta_1)\times(\xi_2,\eta_2)\in
\mathcal{I}^M_1$.

As $(\xi_1,\eta_1)\times(\xi_2,\eta_2)\in \mathcal{I}^M_1$, without loss of generality, we
assume that $(\xi_1,\eta_1)\in {\mathcal D}^{M}_{\frac{3}{4}M}\setminus {\mathcal
D}^{2M}_{\frac{3}{2}M}$ and $(\xi_2,\eta_2)\in {\mathcal D}^{M}_{\frac{3}{4}M}$. Hence,
$$M^{-1}\pi\leq\min\big\{\big|\theta_1-\frac{3}{4}\pi\big|,\big|\theta_1+\frac{1}{4}\pi\big|
\big\}\leq 2M^{-1}\pi,$$
and
$$\min\big\{\big|\theta_2-\frac{3}{4}\pi\big|,\big|\theta_2+\frac{1}{4}\pi\big|\big\}\leq
2M^{-1}\pi.$$
Therefore,
$$|\cos\theta_1+\sin\theta_1|\sim M^{-1},$$
and
$$|\cos\theta_2+\sin\theta_2|\lesssim M^{-1}.$$

For $M\geq 2^{11}$, we use angular decomposition in each $\mathcal{I}^M_1$ :
\begin{align}
\mathcal{I}^M_1=
\bigcup\limits_{M\leq A}\bigcup\limits_{\substack{0\leq j_1,j_2\leq A-1 \\
16\leq|j_1-j_2|\leq32} }\big({\mathcal D}^{A}_{j_1}\times {\mathcal
D}^{A}_{j_2}\big)\cap\mathcal{I}^M_1 .\nonumber
\end{align}
We also write
$$\mathcal{J}^{\mathcal{I}^M_1}_A:=\big\{(j_1,j_2)\ | \ 0\leq j_1,j_2\leq A-1, \  {\mathcal
D}^{A}_{j_1}\times {\mathcal D}^{A}_{j_2}\subset \mathcal{I}^M_1\big\}$$
for simplicity.
\begin{lemma}\label{mathcalI1a0in}
Let $A\geq 2^{5}M$ be dyadic, $\lambda_{\iota}=(\tau_{\iota},\xi_{\iota},\eta_{\iota})\in
G_{N_{\iota},L_{\iota}}\cap {\widetilde{\mathcal D}}^{A}_{j_{\iota}}$ $({\iota}=1,2)$ and
$\lambda_1+\lambda_2\in G_{N_3,L_3}$. Assume that $(j_1,j_2)\in
\mathcal{J}^{\mathcal{I}^M_1}_A$ and $|j_1-j_2|\leq 32$.
In the assumption of Case 1.1,
\begin{itemize}
 \item[\rm (i)] if $M\ll\gamma_0^{-1}$, then $\gamma_0N_3\ll A^{-1}N_1$;
\item[\rm (ii)] if $\gamma_0^{-1}\lesssim M$ and $\frac{MN_1}{N_3}\ll A$, then
    $\max_{1\leq {\iota}\leq3}\{L_{\iota}\}\ll A^{-1}N_1^3$.
\end{itemize}
\end{lemma}
{\bf Proof.} We write
$(\xi_{\iota},\eta_{\iota})=r_{\iota}(\cos\theta_{\iota},\sin\theta_{\iota})$ for
${\iota}=1,2$. Since
$$|\cos\theta_1+\sin\theta_1|\sim M^{-1}$$
and
$$|(\cos\theta_1,\sin\theta_1)+(\cos\theta_2,\sin\theta_2)|\leq2^8A^{-1},$$
one has
\begin{align}
&|\xi_1\xi_2(\xi_1+\xi_2)+\eta_1\eta_2(\eta_1+\eta_2)|\notag \\
=&\big|r_1r_2(r_1-r_2)(\cos^3\theta_1+\sin^3\theta_1)-r_1r_2(r_1-r_2)\cos^2\theta_1
(\cos\theta_1+\cos\theta_2)\notag \\
&-r_1r_2(r_1-r_2)\sin^2\theta_1(\sin\theta_1+\sin\theta_2)-r_1r^2_2\cos\theta_1\cos\theta_2
(\cos\theta_1+\cos\theta_2)\notag \\
&-r_1r^2_2\sin\theta_1\sin\theta_2(\sin\theta_1+\sin\theta_2)\big|\notag \\
\gtrsim &M^{-1}r_1r_2|r_1-r_2|-A^{-1}r_1r_2|r_1-r_2|-A^{-1}r_1r^2_2. \label{LN3item12}
\end{align}

If $M\ll\gamma_0^{-1}$ and $ A^{-1}N_1\lesssim \gamma_0N_3$, from \eqref{LN3item12} we get
\begin{align}
\gamma_0N_1^2N_3\gg &\max_{1\leq {\iota}\leq3}\{L_{\iota}\}\notag \\
\gtrsim& |\xi_1\xi_2(\xi_1+\xi_2)+\eta_1\eta_2(\eta_1+\eta_2)|\notag \\
\gtrsim&M^{-1}N_1^2|r_1-r_2|-A^{-1}N_1^3. \nonumber
\end{align}
Thus $|r_1-r_2|\ll M\gamma_0N_3$.
However
\begin{align}
 N_3=&|(\xi_1+\xi_2,\eta_1+\eta_2)|\notag \\
 \leq&|(r_1\cos\theta_1-r_2\cos\theta_1,r_1\sin\theta_1-r_2\sin\theta_1)|\notag \\
 & \ +r_2|(\cos\theta_1,\sin\theta_1)+(\cos\theta_2,\sin\theta_2)|\notag \\
 \lesssim & |r_1-r_2|+A^{-1}N_1\notag \\
 \lesssim& M\gamma_0N_3+\gamma_0N_3 \ll N_3. \nonumber
\end{align}
This is impossible. Hence $\gamma_0N_3\ll A^{-1}N_1$ if $M\ll\gamma_0^{-1}$.

If $\gamma_0^{-1}\lesssim M\ll AN_3/N_1$ and $\max_{1\leq {\iota}\leq3}\{L_{\iota}\}\lesssim
A^{-1}N_1^3$, from \eqref{LN3item12} we have
$$M^{-1}N_1^2|r_1-r_2|-A^{-1}N_1^3\lesssim |\xi_1\xi_2\xi_3+\eta_1\eta_2\eta_3|\lesssim
\max_{1\leq {\iota}\leq3}\{L_{\iota}\}\lesssim A^{-1}N_1^3,$$
which implies $|r_1-r_2|\lesssim MA^{-1}N_1$.
But
\begin{align}
 N_3=&|(\xi_1+\xi_2,\eta_1+\eta_2)|\notag \\
 \lesssim & |r_1-r_2|+A^{-1}N_1\notag \\
 \lesssim& MA^{-1}N_1+A^{-1}N_1\notag \\
 \lesssim& N_3+\gamma_0N_3 \ll N_3, \nonumber
\end{align}
which gives contradiction.

Thus,
$$\max_{1\leq {\iota}\leq3}\{L_{\iota}\}\ll A^{-1}N_1^3$$
if $\gamma_0^{-1}\lesssim M$ and $\frac{MN_1}{N_3}\gg A$.
\begin{lemma}\label{mathcalI1ain}
Let $A\geq 2^{5}M$ be dyadic. Assume that $(j_1,j_2)\in \mathcal{J}^{\mathcal{I}^M_1}_A$ and
$|j_1-j_2|\leq 32$.
In the assumption of Case 1.1, we have
\begin{align}
 &\big\|\chi_{G_{N_3,L_3}}(\lambda_3)\int \widehat{v}_1\big|_{{\widetilde{\mathcal
D}}^{A}_{j_1}}(\lambda_1) \widehat{v}_2\big|_{{\widetilde{\mathcal
D}}^{A}_{j_2}}(\lambda_1+\lambda_3)d\lambda_1\big\|_{L^2_{\lambda_3}}\notag \\
 \lesssim &
(AN_1)^{-\frac{1}{2}}(L_1L_2)^{\frac{1}{2}}\big\|\widehat{v}_{1}|_{{\widetilde{\mathcal
D}}^{A}_{j_1}}\big\|_{L^{2}}\big\|\widehat{v}_{2}|_{{\widetilde{\mathcal
D}}^{A}_{j_2}}\big\|_{L^{2}},\label{case11bles1in}
\end{align}
\begin{align}
 &\big\|\chi_{G_{N_2,L_2}}(\lambda_2)\int \widehat{v}_1\big|_{{\widetilde{\mathcal
D}}^{A}_{j_1}}(\lambda_1)
\widehat{v}_3(\lambda_1+\lambda_2)d\lambda_1\big\|_{L^2_{\lambda_2}}\notag \\
 \lesssim &
(AN_1)^{-\frac{1}{2}}(L_1L_3)^{\frac{1}{2}}\big\|\widehat{v}_{1}|_{{\widetilde{\mathcal
D}}^{A}_{j_1}}\big\|_{L^{2}}\big\|\widehat{v}_{3}\big\|_{L^{2}},\label{case11bles2in}
\end{align}
and
\begin{align}
 &\big\|\chi_{G_{N_1,L_1}}(\lambda_1)\int \widehat{v}_2\big|_{{\widetilde{\mathcal
D}}^{A}_{j_2}}(\lambda_2)
\widehat{v}_3(\lambda_1+\lambda_2)d\lambda_2\big\|_{L^2_{\lambda_1}}\notag \\
 \lesssim &
(AN_1)^{-\frac{1}{2}}(L_2L_3)^{\frac{1}{2}}\big\|\widehat{v}_{2}|_{{\widetilde{\mathcal
D}}^{A}_{j_2}}\big\|_{L^{2}}\big\|\widehat{v}_{3}\big\|_{L^{2}},\label{case11bles3in}
\end{align}
where $supp \ \widehat{v}_{\iota}\subset G_{N_{\iota},L_{\iota}}$ for ${\iota}=1,2,3$.
\end{lemma}
{\bf Proof.} In fact, the proof is exactly the same as that of Lemma \ref{mathcalI1a}.

\begin{lemma}\label{mathcalI1bin}
Let $ \max\{2^{25},M\}\leq A$ be dyadic. Suppose that $(j_1,j_2)\in
\mathcal{J}^{\mathcal{I}^M_1}_A$ , $16\leq|j_1-j_2|\leq 32$ and $N_3\lesssim A^{-1}MN_1$.
In the assumption of Case 1.1, we have
\begin{align}
&\bigg|\iint_{\Omega_r^{*}}\widehat{v}_{1}|_{{\widetilde{\mathcal
D}}^{A}_{j_1}}(\lambda_1)\widehat{v}_{2}|_{{\widetilde{\mathcal
D}}^{A}_{j_2}}(\lambda_2)\widehat{v}_{3}(\lambda_3)\bigg| \notag \\
 \lesssim &
A^{\frac{1}{2}}MN_1^{-2}(L_1L_2L_3)^{\frac{1}{2}}\big\|\widehat{v}_{1}|_{{\widetilde{\mathcal
D}}^{A}_{j_1}}\big\|_{L^{2}}\big\|\widehat{v}_{2}|_{{\widetilde{\mathcal
D}}^{A}_{j_2}}\big\|_{L^{2}}\big\|\widehat{v}_{3}\big\|_{L^{2}},\label{case11b1in}
\end{align}
where $supp \ \widehat{v}_{\iota}\subset G_{N_{\iota},L_{\iota}}$ for ${\iota}=1,2,3$.
\end{lemma}
{\bf Proof.} We imitate the proof of Lemma \ref{mathcalI1b} to show \eqref{case11b1in}. Let
$A'=2^{30}A$. Since $N_3\lesssim A^{-1}MN_1$, by the almost orthogonality, we may assume
that the integral region are restricted to sets of length $A^{-1}MN_1$. That's to say,
$$supp \ \widehat{v}_{\iota}\subset \widetilde{\mathcal{S}}^{A^{-1}MN_1}_{l_{\iota}} \ \ \
for \ \ \iota=1,2 $$
where $|l_1-l_2|\leq 2^{10}$ and
$$\widetilde{\mathcal{S}}^{A^{-1}MN_1}_{l}:=\big\{(\tau,\xi,\eta)\in \mathbb{R}^3\ \big| \
N_1+l A^{-1}MN_1\leq \langle\xi\rangle\leq N_1+(l+1) A^{-1}MN_1\big\}.$$

Note that $\widetilde{\mathcal{S}}^{A^{-1}MN_1}_{l_{\iota}}\cap {\widetilde{\mathcal
D}}^{A}_{j_{\iota}}$ is contained in a rectangle whose long side length is $\sim A^{-1}MN_1$
and short side length is $\sim A^{-1}N_1$. Thus we shall divide the rectangle further. One
can easily verify that the number of tiles $\widetilde{\mathcal{T}}^{A'}_{k_{\iota}}$
satisfying
$$\widetilde{\mathcal{S}}^{A^{-1}MN_1}_{l_{\iota}}\cap {\widetilde{\mathcal
D}}^{A}_{j_{\iota}}\cap \widetilde{\mathcal{T}}^{A'}_{k_{\iota}}\not=\emptyset,\ \ \
\widetilde{\mathcal{S}}^{A^{-1}MN_1}_{l_{\iota}}\cap {\widetilde{\mathcal
D}}^{A}_{j_{\iota}}\subset \bigcup\limits_{\#k_{\iota}\sim
M}\widetilde{\mathcal{T}}^{A'}_{k_{\iota}}.$$
is approximately $M$.

Assume that $f,g$ and $h$ are functions satisfying
$$supp \ f\subset G_{N_1,L_1}\cap\widetilde{\mathcal{T}}^{A'}_{k_1}, \hspace{2mm}supp \
g\subset G_{N_2,L_2}\cap\widetilde{\mathcal{T}}^{A'}_{k_2},\hspace{2mm} supp \ h\subset
G_{N_3,L_3}, $$
where $\widetilde{\mathcal{T}}^{A'}_{k_{\iota}}\cap
\widetilde{\mathcal{S}}^{A^{-1}MN_1}_{l_{\iota}}\cap {\widetilde{\mathcal
D}}^{A}_{j_{\iota}}\not=\emptyset$ for $\iota=1,2$ , we reduce to showing that
\begin{align}
\bigg|\int_{*}f(\lambda_1)g(\lambda_2)h(\lambda_3)\bigg|
 \lesssim
A^{\frac{1}{2}}N_1^{-2}(L_1L_2L_3)^{\frac{1}{2}}\|f\|_{L^{2}}\|g\|_{L^{2}}\|h\|_{L^{2}}
\label{mathcalI1bin111}
\end{align}
which can be proved by the same way as Lemma \ref{mathcalI1b}.

From \eqref{mathcalI1bin111}, we have
\begin{align}
&\bigg|\iint_{\Omega_r^{*}}\widehat{v}_{1}|_{{\widetilde{\mathcal
D}}^{A}_{j_1}}(\lambda_1)\widehat{v}_{2}|_{{\widetilde{\mathcal
D}}^{A}_{j_2}}(\lambda_2)\widehat{v}_{3}(\lambda_3)\bigg| \notag \\
 \lesssim & \bigg|\iint_{\Omega_r^{*}}\widehat{v}_{1}|_{{\widetilde{\mathcal
D}}^{A}_{j_1}\cap\widetilde{\mathcal{S}}^{A^{-1}MN_1}_{l_1}}(\lambda_1)\cdot\widehat{v}_{2}|
_{{\widetilde{\mathcal
D}}^{A}_{j_2}\cap\widetilde{\mathcal{S}}^{A^{-1}MN_1}_{l_2}}(\lambda_2)\cdot\widehat{v}_{3}
(\lambda_3)\bigg| \notag \\
  \lesssim & \sum\limits_{\#k_1\sim M}\sum\limits_{\#k_2\sim
M}\bigg|\iint_{\Omega_r^{*}}\widehat{v}_{1}|_{\widetilde{\mathcal{T}}^{A'}_{k_1}}(\lambda_1)
\cdot\widehat{v}_{2}|_{\widetilde{\mathcal{T}}^{A'}_{k_2}}(\lambda_2)\cdot\widehat{v}_{3}
(\lambda_3)\bigg| \notag \\
 \lesssim & A^{\frac{1}{2}} N_1^{-2}(L_1L_2L_3)^{\frac{1}{2}}\sum\limits_{\#k_1\sim
M}\sum\limits_{\#k_2\sim M}\big\|\widehat{v}_{1}|_{{\widetilde{\mathcal
T}}^{A'}_{k_1}}\big\|_{L^{2}}\big\|\widehat{v}_{2}|_{{\widetilde{\mathcal
T}}^{A'}_{k_2}}\big\|_{L^{2}}\big\|\widehat{v}_{3}\big\|_{L^{2}}\notag \\
 \lesssim & A^{\frac{1}{2}}M
N_1^{-2}(L_1L_2L_3)^{\frac{1}{2}}\big\|\widehat{v}_{1}|_{{\widetilde{\mathcal
D}}^{A}_{j_1}}\big\|_{L^{2}}\big\|\widehat{v}_{2}|_{{\widetilde{\mathcal
D}}^{A}_{j_2}}\big\|_{L^{2}}\big\|\widehat{v}_{3}\big\|_{L^{2}}\nonumber
\end{align}
which concludes the proof of Lemma \ref{mathcalI1bin}.
\begin{remark}
In Lemma \ref{mathcalI1bin}, the upper bound is $A^{\frac{1}{2}}M N_1^{-2}$ ( see RHS of
\eqref{case11b1in} ) which is different from $(3.60)$ in Proposition 3.17 of
\cite{Kinoshita19} where the upper bound is $A^{\frac{1}{2}}M^{\frac{1}{2}} N_1^{-2}$.
\end{remark}

\begin{proposition}\label{mathcalI1in}
Assume that $(\xi_1,\eta_1)\times(\xi_2,\eta_2)\in \mathcal{I}_1$ , then we have
\begin{align}
&\bigg|\iint_{\Omega_r^{*}}(\xi_3+\eta_3)\widehat{v}_{1}(\lambda_1)\widehat{v}_{2}(\lambda_2)
\widehat{v}_{3}(\lambda_3)\bigg|\notag \\
\lesssim &
C(\gamma_0,N_1,N_3)(L_1L_2L_3)^{\frac{1}{2}}\|\widehat{v}_{1}\|_{L^{2}}\|\widehat{v}_{2}\|
_{L^{2}}
\|\widehat{v}_{3}\|_{L^{2}}\nonumber
\end{align}
where
$$C(\gamma_0,N_1,N_3)=\gamma_0^{-\frac{1}{2}-}N_1^{-\frac{1}{2}}N_3^{-\frac{1}{2}}+
\gamma_0^{1+}N_1^{-1+}N_3^{0+}+N_1^{-\frac{1}{12}}N_3^{-\frac{1}{3}}.$$
\end{proposition}
{\bf Proof.} Since $(\xi_1,\eta_1)\times(\xi_2,\eta_2)\in \mathcal{I}_1^M$,
$$|\xi_3+\eta_3|\lesssim |\xi_1+\eta_1|+|\xi_2+\eta_2|\lesssim M^{-1}N_1.$$

According to the twice decomposition of $\mathcal{I}_1$ , we have
\begin{align}
&\bigg|\iint_{\Omega_r^{*}}(\xi_3+\eta_3)\widehat{v}_{1}(\lambda_1)\widehat{v}_{2}(\lambda_2)
\widehat{v}_{3}(\lambda_3)\bigg|\notag \\
\lesssim&\sum_{M\ll \gamma_0^{-1}}
\sum_{M\leq A}\sum\limits_{\substack{(j_1,j_2)\in\mathcal{J}^{\mathcal{I}^M_1}_A \\
16\leq|j_1-j_2|\leq32} }\bigg|\int_{*}(\xi_3+\eta_3)\widehat{v}_{1}|_{{\widetilde{\mathcal
D}}^{A}_{j_1}}(\lambda_1)\widehat{v}_{2}|_{{\widetilde{\mathcal
D}}^{A}_{j_2}}(\lambda_2)\widehat{v}_{3}(\lambda_3)\bigg|\notag \\
& +\sum_{\gamma_0^{-1}\lesssim M \lesssim M_0}
\sum_{M\leq A}\sum\limits_{\substack{(j_1,j_2)\in\mathcal{J}^{\mathcal{I}^M_1}_A \\
16\leq|j_1-j_2|\leq32} }\bigg|\int_{*}(\xi_3+\eta_3)\widehat{v}_{1}|_{{\widetilde{\mathcal
D}}^{A}_{j_1}}(\lambda_1)\widehat{v}_{2}|_{{\widetilde{\mathcal
D}}^{A}_{j_2}}(\lambda_2)\widehat{v}_{3}(\lambda_3)\bigg|\notag \\
& +\bigg|\int_{*}(\xi_3+\eta_3)\widehat{v}_{1}\big|_{{\widetilde{\mathcal
D}}^{M_0}_{\frac{3}{4}M_0}}(\lambda_1)\widehat{v}_{2}\big|_{{\widetilde{\mathcal
D}}^{M_0}_{\frac{3}{4}M_0}}(\lambda_2)\widehat{v}_{3}(\lambda_3)\bigg|.\label{I1term123}
\end{align}
We denote the right-hand side of \eqref{I1term123} by $T_1,T_2,T_3$ and estimate these three
terms respectively.

First, we consider the contribution of $T_1$.
From Lemma \ref{mathcalI1a0in} (i), we know that $A\ll\gamma_0^{-1}N_1N_3^{-1}$ and
$$|r_1-r_2|\lesssim A^{-1}MN_1, \ \ N_3\lesssim |r_1-r_2|+A^{-1}N_1\lesssim A^{-1}MN_1 $$
when $M\ll \gamma_0^{-1}$.

Applying Lemma \ref{mathcalI1bin}, we have
\begin{align}
T_1\lesssim&\sum_{M\ll \gamma_0^{-1}}
\sum_{A\ll\gamma_0^{-1}N_1N_3^{-1}}\sum\limits_{\substack{(j_1,j_2)\in\mathcal{J}^
{\mathcal{I}^M_1}_A \\ 16\leq|j_1-j_2|\leq32}
}\bigg|\int_{*}(\xi_3+\eta_3)\widehat{v}_{1}|_{{\widetilde{\mathcal
D}}^{A}_{j_1}}(\lambda_1)\widehat{v}_{2}|_{{\widetilde{\mathcal
D}}^{A}_{j_2}}(\lambda_2)\widehat{v}_{3}(\lambda_3)\bigg|\notag \\
\lesssim&\sum_{M\ll \gamma_0^{-1}}
\sum_{A\ll\gamma_0^{-1}N_1N_3^{-1}}\sum\limits_{\substack{(j_1,j_2)\in\mathcal{J}^
{\mathcal{I}^M_1}_A \\ 16\leq|j_1-j_2|\leq32} }M^{-1}N_1A^{\frac{1}{2}}M
N_1^{-2}(L_1L_2L_3)^{\frac{1}{2}}\times\notag \\
&\hspace{55mm}\big\|\widehat{v}_{1}|_{{\widetilde{\mathcal
D}}^{A}_{j_1}}\big\|_{L^{2}}\big\|\widehat{v}_{2}|_{{\widetilde{\mathcal
D}}^{A}_{j_2}}\big\|_{L^{2}}\big\|\widehat{v}_{3}\big\|_{L^{2}}\notag \\
\lesssim&\gamma_0^{-\frac{1}{2}-}N_1^{\frac{1}{2}}N_3^{-\frac{1}{2}}N_1^{-1}(L_1L_2L_3)^
{\frac{1}{2}}
\|\widehat{v}_{1}\|_{L^{2}}\|\widehat{v}_{2}\|_{L^{2}}\|\widehat{v}_{3}\|_{L^{2}}\notag \\
\lesssim&\gamma_0^{-\frac{1}{2}-}(N_1N_3)^{-\frac{1}{2}}(L_1L_2L_3)^{\frac{1}{2}}
\|\widehat{v}_{1}\|_{L^{2}}\|\widehat{v}_{2}\|_{L^{2}}\|\widehat{v}_{3}\|_{L^{2}}.
\label{pp6T1}
\end{align}

Secondly, we consider $T_2$ which will be split into two terms based on the size of $A$. On
one hand, we use the nonlinear Loomis-Whitney inequality if $A\lesssim MN_1/N_3$. On the
other hand, we take advantage of bilinear Strichartz estimates if $A\gg MN_1/N_3$.

Let us write $T_2=T_{2,1}+T_{2,2}$, where
\begin{align}T_{2,1}:=\sum_{\gamma_0^{-1}\lesssim M \lesssim M_0}
\sum_{A\lesssim\frac{MN_1}{N_3}}\sum\limits_{\substack{(j_1,j_2)\in\mathcal{J}^
{\mathcal{I}^M_1}_A \\ 16\leq|j_1-j_2|\leq32}
}\bigg|\int_{*}(\xi_3+\eta_3)\widehat{v}_{1}|_{{\widetilde{\mathcal
D}}^{A}_{j_1}}(\lambda_1)\widehat{v}_{2}|_{{\widetilde{\mathcal
D}}^{A}_{j_2}}(\lambda_2)\widehat{v}_{3}(\lambda_3)\bigg|\nonumber
\end{align}
and
\begin{align}T_{2,2}:=\sum_{\gamma_0^{-1}\lesssim M \lesssim M_0}
\sum_{\frac{MN_1}{N_3}\ll A}\sum\limits_{\substack{(j_1,j_2)\in\mathcal{J}
^{\mathcal{I}^M_1}_A \\ 16\leq|j_1-j_2|\leq32}
}\bigg|\int_{*}(\xi_3+\eta_3)\widehat{v}_{1}|_{{\widetilde{\mathcal
D}}^{A}_{j_1}}(\lambda_1)\widehat{v}_{2}|_{{\widetilde{\mathcal
D}}^{A}_{j_2}}(\lambda_2)\widehat{v}_{3}(\lambda_3)\bigg|.\nonumber
\end{align}

For $T_{2,1}$, since $N_3\lesssim A^{-1}MN_1$, by Lemma \ref{mathcalI1bin}, we have
\begin{align}T_{2,1}\lesssim&\sum_{\gamma_0^{-1}\lesssim M \lesssim M_0}
\sum_{A\lesssim\frac{MN_1}{N_3}}\sum\limits_{\substack{(j_1,j_2)\in\mathcal{J}^
{\mathcal{I}^M_1}_A \\ 16\leq|j_1-j_2|\leq32} }M^{-1}N_1A^{\frac{1}{2}}M
N_1^{-2}(L_1L_2L_3)^{\frac{1}{2}}\times\notag \\
&\hspace{55mm}\big\|\widehat{v}_{1}|_{{\widetilde{\mathcal
D}}^{A}_{j_1}}\big\|_{L^{2}}\big\|\widehat{v}_{2}|_{{\widetilde{\mathcal
D}}^{A}_{j_2}}\big\|_{L^{2}}\big\|\widehat{v}_{3}\big\|_{L^{2}}\notag \\
\lesssim&\sum_{\gamma_0^{-1}\lesssim M \lesssim M_0}\big(\frac{M
N_1}{N_3}\big)^{\frac{1}{2}}N_1^{-1}(L_1L_2L_3)^{\frac{1}{2}}
\|\widehat{v}_{1}\|_{L^{2}}\|\widehat{v}_{2}\|_{L^{2}}\|\widehat{v}_{3}\|_{L^{2}}\notag \\
\lesssim&M_0^{\frac{1}{2}}(N_1N_3)^{-\frac{1}{2}}(L_1L_2L_3)^{\frac{1}{2}}
\|\widehat{v}_{1}\|_{L^{2}}\|\widehat{v}_{2}\|_{L^{2}}\|\widehat{v}_{3}\|_{L^{2}}.
\label{pp6T21}
\end{align}

For $T_{2,2}$, by Lemma \ref{mathcalI1a0in} (i), we know that $\max_{1\leq \iota\leq
3}\{L_{\iota}\}\gg A^{-1}N_1^3$. Using Lemma \ref{mathcalI1ain},
\begin{align}T_{2,2}\lesssim&\sum_{\gamma_0^{-1}\lesssim M \lesssim M_0}
\sum_{\frac{MN_1}{N_3}\ll A}\sum\limits_{\substack{(j_1,j_2)\in\mathcal{J}
^{\mathcal{I}^M_1}_A \\ 16\leq|j_1-j_2|\leq32}
}M^{-1}N_1(AN_1)^{-\frac{1}{2}}A^{\frac{1}{2}-}N_1^{-\frac{3}{2}+}\times\notag \\
&\hspace{40mm}(L_1L_2L_3)^{\frac{1}{2}}\big\|\widehat{v}_{1}|_{{\widetilde{\mathcal
D}}^{A}_{j_1}}\big\|_{L^{2}}\big\|\widehat{v}_{2}|_{{\widetilde{\mathcal
D}}^{A}_{j_2}}\big\|_{L^{2}}\big\|\widehat{v}_{3}\big\|_{L^{2}}\notag \\
\lesssim&\sum_{\gamma_0^{-1}\lesssim M \lesssim
M_0}M^{-1-}N_1^{-1+}N_3^{0+}(L_1L_2L_3)^{\frac{1}{2}}
\|\widehat{v}_{1}\|_{L^{2}}\|\widehat{v}_{2}\|_{L^{2}}\|\widehat{v}_{3}\|_{L^{2}}\notag \\
\lesssim&\gamma_0^{1+}N_1^{-1+}N_3^{0+}(L_1L_2L_3)^{\frac{1}{2}}
\|\widehat{v}_{1}\|_{L^{2}}\|\widehat{v}_{2}\|_{L^{2}}\|\widehat{v}_{3}\|_{L^{2}}.
\label{pp6T22}
\end{align}

Finally, let's consider $T_3$. By using \eqref{estimate15}, we get
\begin{align}T_{3}\lesssim&M_0^{-1}N_1N_1^{-\frac{1}{4}}\big\||\nabla_x|^{\frac{1}{8}}|\nabla_y|
^{\frac{1}{8}}v_1\big\|
_{L^4}\cdot\big\||\nabla_x|^{\frac{1}{8}}|\nabla_y|^{\frac{1}{8}}v_2\big\|_{L^4}\cdot
\|v_{3}\|_{L^{2}}\notag \\
\lesssim&M_0^{-1}N_1^{\frac{3}{4}}(L_1L_2L_3)^{\frac{1}{2}}
\|\widehat{v}_{1}\|_{L^{2}}\|\widehat{v}_{2}\|_{L^{2}}\|\widehat{v}_{3}\|_{L^{2}}.
\label{pp6T3}
\end{align}

Collecting \eqref{pp6T1}--\eqref{pp6T3} and choosing
$M_0=N_1^{\frac{5}{6}}N_3^{\frac{1}{3}}$, we have
\begin{align}
&\bigg|\iint_{\Omega_r^{*}}(\xi_3+\eta_3)\widehat{v}_{1}(\lambda_1)\widehat{v}_{2}
(\lambda_2)\widehat{v}_{3}(\lambda_3)\bigg|\notag \\
\lesssim&
\big(\gamma_0^{-\frac{1}{2}-}N_1^{-\frac{1}{2}}N_3^{-\frac{1}{2}}+\gamma_0^{1+}N_1^{-1+}N_3^
{0+}+N_1^{-\frac{1}{12}}N_3^{-\frac{1}{3}}\big)(L_1L_2L_3)^{\frac{1}{2}}\|\widehat{v}_{1}\|_
{L^{2}}\|\widehat{v}_{2}\|_{L^{2}}
\|\widehat{v}_{3}\|_{L^{2}}\nonumber
\end{align}
which concludes the proof of Proposition \ref{mathcalI1in}.

Up to now, we have treated $(\xi_1,\eta_1)\times(\xi_2,\eta_2)\notin \mathcal{I}_1$ and
$(\xi_1,\eta_1)\times(\xi_2,\eta_2)\in \mathcal{I}_1$ respectively. Therefore, under the
assumption of Case 1.1, from Proposition \ref{mathcalI1} and Proposition \ref{mathcalI1in},
one obtains that
\begin{align}
&\bigg|\iint_{\Omega_r^{*}}\frac{\sum^3_{{\iota}=1}m^2(\zeta_{\iota})(\xi_{\iota}+\eta_
{\iota})}
{m(\zeta_1)m(\zeta_2)m(\zeta_3)}\widehat{v}_{1}(\lambda_1)\widehat{v}_{2}(\lambda_2)
\widehat{v}_{3}(\lambda_3)\bigg|\notag \\
\lesssim
&\big(\gamma_0^{-\frac{1}{2}-}N^{-1+}+N^{-\frac{5}{12}}\big)(L_1L_2L_3)^{\frac{1}{2}}\|
\widehat{v}_{1}\|_{L^{2}}\|\widehat{v}_{2}\|_{L^{2}}
\|\widehat{v}_{3}\|_{L^{2}}.\label{maintricase11}
\end{align}

In the next place, we consider non-parallel interactions under the high-low frequencies
case. The non-parallel interactions are caused by the size relationship between $|\xi_1|$
and $|\eta_1|$. Actually, if $|\xi_1|\gg |\eta_1|$, as $|\eta_2|=|\eta_1+\eta_3|\ll N_1$ ,
thus $|\xi_1|\sim |\xi_2|\gg |\eta_2|$. We claim that $|\xi_3|\ll|\eta_3|$. Otherwise,
$|\xi_3|\gtrsim|\eta_3|$, then $|\xi_1\xi_2\xi_3+\eta_1\eta_2\eta_3|\sim
|\xi_1\xi_2\xi_3|\sim N_1^2N_3\gg\gamma_0N_1^2N_3$, which contradicts to the assumption
$|\xi_1\xi_2\xi_3+\eta_1\eta_2\eta_3|\ll\gamma_0N_1^2N_3$.

Moreover, for $\iota=1,2$,
$$\big|\xi_{\iota}\eta_3+\xi_3\eta_{\iota}+2(\xi_{\iota}\eta_{\iota}+\xi_3\eta_3)\big|\sim
 |\xi_{\iota}\eta_{\iota}|\gtrsim N_1N_3 \hspace{10mm}  if \  |\eta_{\iota}|\gg|\eta_3|,$$
and
$$\big|\xi_{\iota}\eta_3+\xi_3\eta_{\iota}+2(\xi_{\iota}\eta_{\iota}+\xi_3\eta_3)\big|\sim
|\xi_{\iota}\eta_3|\sim N_1N_3 \hspace{10mm} if \ |\eta_{\iota}|\ll|\eta_3|,$$
hence we can directly use the nonlinear Loomis-Whitney inequality and bilinear Strichartz
estimates.

So, it is reduced to:
\\{\bf Case 1.2} \hspace{2mm} $|\xi_1|\sim |\xi_2|\gg|\eta_1|\sim |\eta_2|\sim |\eta_3|\gg
|\xi_3|$

From the assumption above, it is easy to know that
$$|\xi_1\eta_3-\xi_3\eta_1|\sim|\xi_1\eta_3|\sim N_1N_3.$$
This is the reason why we call them non-parallel interactions.

To bound the high frequency part, we need the following lemma which is knew as bilinear
Strichartz estimate.
\begin{lemma}\label{case12a}
Let $A\geq 2^{25}$ be dyadic, $k_1,k_2\in \mathbb{Z}^2$.
In the assumption of Case 1.2, we have
\begin{align}
 &\big\|\chi_{G_{N_3,L_3}}(\lambda_3)\int \widehat{v}_1\big|_{{\widetilde{\mathcal
T}}^{A}_{k_1}}(\lambda_1) \widehat{v}_2\big|_{{\widetilde{\mathcal
T}}^{A}_{k_2}}(\lambda_1+\lambda_3)d\lambda_1\big\|_{L^2_{\lambda_3}}\notag \\
 \lesssim &
(AN_1)^{-\frac{1}{2}}(L_1L_2)^{\frac{1}{2}}\big\|\widehat{v}_{1}|_{{\widetilde{\mathcal
T}}^{A}_{k_1}}\big\|_{L^{2}}\big\|\widehat{v}_{2}|_{{\widetilde{\mathcal
T}}^{A}_{k_2}}\big\|_{L^{2}},\label{case12ales1}
\end{align}
\begin{align}
 &\big\|\chi_{G_{N_2,L_2}}(\lambda_2)\int \widehat{v}_1\big|_{{\widetilde{\mathcal
T}}^{A}_{k_1}}(\lambda_1)
\widehat{v}_3(\lambda_1+\lambda_2)d\lambda_1\big\|_{L^2_{\lambda_2}}\notag \\
 \lesssim &
(AN_1)^{-\frac{1}{2}}(L_1L_3)^{\frac{1}{2}}\big\|\widehat{v}_{1}|_{{\widetilde{\mathcal
T}}^{A}_{k_1}}\big\|_{L^{2}}\big\|\widehat{v}_{3}\big\|_{L^{2}},\label{case12ales2}
\end{align}
and
\begin{align}
 &\big\|\chi_{G_{N_1,L_1}}(\lambda_1)\int \widehat{v}_2\big|_{{\widetilde{\mathcal
T}}^{A}_{k_2}}(\lambda_2)
\widehat{v}_3(\lambda_1+\lambda_2)d\lambda_2\big\|_{L^2_{\lambda_1}}\notag \\
 \lesssim &
(AN_1)^{-\frac{1}{2}}(L_2L_3)^{\frac{1}{2}}\big\|\widehat{v}_{2}|_{{\widetilde{\mathcal
T}}^{A}_{k_2}}\big\|_{L^{2}}\big\|\widehat{v}_{3}\big\|_{L^{2}},\label{case12ales3}
\end{align}
where $supp \ \widehat{v}_{\iota}\subset G_{N_{\iota},L_{\iota}}$ for ${\iota}=1,2,3$.
\end{lemma}
{\bf Proof.} By the assumption of Case 1.2,
$$|\xi_1|\sim N_1,\ \ \ |\xi_2|\sim N_2.$$
According to the definition of ${\widetilde{\mathcal T}}^{A}_{k}$, we have
$$\Big|\left\{\eta_{\iota}\big|\lambda_{\iota}=(\tau_{\iota},\xi_{\iota},\eta_{\iota})\in
{\widetilde{\mathcal T}}^{A}_{k_{\iota}} \right\}\Big|\lesssim A^{-1}N_1$$
for ${\iota}=1,2$.

Thus, from Lemma \ref{X Lema3}, we know this lemma holds true.
\begin{lemma}\label{case12b}
Let $A\geq 2^{25}$ be dyadic. Suppose that
$$|\xi_1\xi_2(\xi_1+\xi_2)+\eta_1\eta_2(\eta_1+\eta_2)|\geq A^{-1}N_1^3 \ \ \ for  \  \
(\xi_{\iota},\eta_{\iota})\in{\mathcal T}^{A}_{k_{\iota}} $$
where ${\iota}=1,2$. In the assumption of Case 1.2, then we have
\begin{align}
&\bigg|\iint_{*}\widehat{v}_{1}|_{{\widetilde{\mathcal
T}}^{A}_{k_1}}(\lambda_1)\widehat{v}_{2}|_{{\widetilde{\mathcal
T}}^{A}_{k_2}}(\lambda_2)\widehat{v}_{3}(\lambda_3)\bigg| \notag \\
 \lesssim & N_1^{-2}(L_1L_2L_3)^{\frac{1}{2}}\big\|\widehat{v}_{1}|_{{\widetilde{\mathcal
T}}^{A}_{k_1}}\big\|_{L^{2}}\big\|\widehat{v}_{2}|_{{\widetilde{\mathcal
T}}^{A}_{k_2}}\big\|_{L^{2}}\big\|\widehat{v}_{3}\big\|_{L^{2}},\label{case12b1}
\end{align}
where $supp \ \widehat{v}_{\iota}\subset G_{N_{\iota},L_{\iota}}$ for ${\iota}=1,2,3$.
\end{lemma}
{\bf Proof.} Without loss of generality, one can assume that
\begin{align}L_1=\max_{1\leq {\iota}\leq3}\{L_{\iota}\}\gtrsim
|\xi_1\xi_2(\xi_1+\xi_2)+\eta_1\eta_2(\eta_1+\eta_2)|\geq
A^{-1}N_1^3.\label{case12b1ina}\end{align}

Using \eqref{case12ales3} and \eqref{case12b1ina}, we have
\begin{align}
&\bigg|\iint_{*}\widehat{v}_{1}|_{{\widetilde{\mathcal
T}}^{A}_{k_1}}(\lambda_1)\widehat{v}_{2}|_{{\widetilde{\mathcal
T}}^{A}_{k_2}}(\lambda_2)\widehat{v}_{3}(\lambda_3)\bigg| \notag \\
 \lesssim &
(AN_1)^{-\frac{1}{2}}L_1^{-\frac{1}{2}}(L_1L_2L_3)^{\frac{1}{2}}\big\|\widehat{v}_{1}|_
{{\widetilde{\mathcal
T}}^{A}_{k_1}}\big\|_{L^{2}}\big\|\widehat{v}_{2}|_{{\widetilde{\mathcal
T}}^{A}_{k_2}}\big\|_{L^{2}}\big\|\widehat{v}_{3}\big\|_{L^{2}}\notag \\
\lesssim & N_1^{-2}(L_1L_2L_3)^{\frac{1}{2}}\big\|\widehat{v}_{1}|_{{\widetilde{\mathcal
T}}^{A}_{k_1}}\big\|_{L^{2}}\big\|\widehat{v}_{2}|_{{\widetilde{\mathcal
T}}^{A}_{k_2}}\big\|_{L^{2}}\big\|\widehat{v}_{3}\big\|_{L^{2}}.\nonumber
\end{align}
This concludes the lemma.

To estimate the low frequency part, we need the nonlinear Loomis-Whitney inequality.
\begin{lemma}\label{case12ca1}
Let $ 2^{25}\leq A$ be dyadic. Suppose that
$$\big|\xi_1\eta_2+\xi_2\eta_1+2(\xi_1\eta_1+\xi_2\eta_2)\big|\gtrsim A^{-1}N_1^{2}\ \ \ for
\  \ (\xi_{\iota},\eta_{\iota})\in{\mathcal T}^{A}_{k_{\iota}} $$
where ${\iota}=1,2$.  In the assumption of Case 1.2, we have
\begin{align}
&\bigg|\iint_{\Omega_r^{*}}\widehat{v}_{1}|_{{\widetilde{\mathcal
T}}^{A}_{k_1}}(\lambda_1)\widehat{v}_{2}|_{{\widetilde{\mathcal
T}}^{A}_{k_2}}(\lambda_2)\widehat{v}_{3}(\lambda_3)\bigg| \notag \\
 \lesssim &
A^{\frac{1}{2}}N_1^{-\frac{3}{2}}N_3^{-\frac{1}{2}}(L_1L_2L_3)^{\frac{1}{2}}\big\|\widehat{v}
_{1}|_{{\widetilde{\mathcal
T}}^{A}_{k_1}}\big\|_{L^{2}}\big\|\widehat{v}_{2}|_{{\widetilde{\mathcal
T}}^{A}_{k_2}}\big\|_{L^{2}}\big\|\widehat{v}_{3}\big\|_{L^{2}},\label{case11b1}
\end{align}
where $supp \ \widehat{v}_{\iota}\subset G_{N_{\iota},L_{\iota}}$ for ${\iota}=1,2,3$.
\end{lemma}
{\bf Proof.} The proof is very similar to that of Lemma \ref{mathcalI1b}. We just sketch it
.

Let $A'=2^{30}A$. Assume that $f,g$ and $h$ are functions satisfying
$$supp \ f\subset G_{N_1,L_1}\cap\widetilde{\mathcal{T}}^{A'}_{k_1}, \hspace{2mm}supp \
g\subset G_{N_2,L_2}\cap\widetilde{\mathcal{T}}^{A'}_{k_2},\hspace{2mm} supp \ h\subset
G_{N_3,L_3}\cap\widetilde{\mathcal{T}}^{A'}_{k_3}, $$
we reduce to showing that
\begin{align}
\bigg|\int_{*}f(\lambda_1)g(\lambda_2)h(\lambda_3)\bigg|
 \lesssim
A^{\frac{1}{2}}N_1^{-\frac{3}{2}}N_3^{-\frac{1}{2}}(L_1L_2L_3)^{\frac{1}{2}}\|f\|_{L^{2}}
\|g\|_{L^{2}}\|h\|_{L^{2}}.\nonumber
\end{align}

By changing variables $c_1=\tau_1-\xi_1^3-\eta_1^3$, $c_2=\tau_2-\xi_2^3-\eta_2^3$, then
decomposing $h$ into $L_3$ pieces and applying the Cauchy-Schwarz inequality, it suffices to
prove that
\begin{align}
&\bigg|\int
f\big(\phi_{c_1}(\xi_1,\eta_1)\big)g\big(\phi_{c_2}(\xi_2,\eta_2)\big)h\big(\phi_{c_1}
(\xi_1,\eta_1)
+\phi_{c_2}(\xi_2,\eta_2)\big)d\xi_1d\eta_1d\xi_2d\eta_2\bigg|  \notag \\
 \lesssim &
A^{\frac{1}{2}}N_1^{-\frac{3}{2}}N_3^{-\frac{1}{2}}\|f\circ\phi_{c_1}\|_{L^{2}_{\xi\eta}}
\|g\circ\phi_{c_2}\|_{L^{2}_{\xi\eta}}\|h\|_{L^{2}_{\tau\xi\eta}},\nonumber
\end{align}
where $h(\tau,\xi,\eta)$ is supported in $c_0\leq\tau-\xi^3-\eta^3\leq c_0+1$ and
$$\phi_{c_{\iota}}(\xi_{\iota},\eta_{\iota})=(\xi_{\iota}^3+\eta_{\iota}^3+c_{\iota},
\xi_{\iota},\eta_{\iota})\  \  \ \  \  for \  \ {\iota}=1,2.$$

Using the scaling $(\tau,\xi,\eta)\rightarrow(N_1^3\tau,N_1\xi,N_1\eta)$ to define
 \begin{align}
 \tilde f(\tau_1,\xi_1,\eta_1)=f(N_1^3\tau_1,N_1\xi_1,N_1\eta_1),\notag \\
 \tilde g(\tau_2,\xi_2,\eta_2)=g(N_1^3\tau_2,N_1\xi_2,N_1\eta_2),\notag \\
  \tilde h(\tau_3,\xi_3,\eta_3)=h(N_1^3\tau_3,N_1\xi_3,N_1\eta_3),\nonumber
\end{align}
then we need to show
\begin{align}
\big\|\tilde  f|_{S_1}*\tilde  g|_{S_2}\big\|_{L^2(S_3)}
 \lesssim  A^{\frac{1}{2}}N_1^{\frac{1}{2}}N_3^{-\frac{1}{2}}\|\tilde
f\|_{L^2(S_1)}\|\tilde  g\|_{L^2(S_2)}.\label{case12b4}
\end{align}

We can verify that transversal conditions are satisfied, thus from the nonlinear
Loomis-Whitney inequality it suffices to estimate the determinant consisting of the unit
normal vectors.

Since
\begin{align}
\big|\xi'_1\eta'_2-\xi'_2\eta'_1\big|=\big|\xi'_1\eta'_3-\xi'_3\eta'_1\big|
\sim N_1^{-2}\big|\xi_1\eta_3-\xi_3\eta_1\big|\sim N_3/N_1\nonumber
\end{align}
and
\begin{align}
\big|\xi'_1\eta'_2+\xi'_2\eta'_1+2(\xi'_1\eta'_1+\xi'_2\eta'_2)\big|
\sim N_1^{-2}\big|\xi_1\eta_2+\xi_2\eta_1+2(\xi_1\eta_1+\xi_2\eta_2)\big|\gtrsim
A^{-1},\nonumber
\end{align}
we have
\begin{align}
&|det\mathcal {N}(\lambda_1,\lambda_2,\lambda_3)|\notag \\
\gtrsim&
\big|\xi'_1\eta'_2-\xi'_2\eta'_1\big|\cdot\big|\xi'_1\eta'_2+\xi'_2\eta'_1+2(\xi'_1\eta'_1
+\xi'_2\eta'_2)\big|\notag \\
\gtrsim& A^{-1}N_1^{-1}N_3,\nonumber
\end{align}
which implies \eqref{case12b4} by Lemma \ref{LW1}. We finish the proof of Lemma
\ref{case12ca1}.

In order to utilize the nonlinear Loomis-Whitney inequality and thus take advantage of Lemma
\ref{case12ca1} for low frequency part, we need to get the lower bound of
$\big|\xi_1\eta_3+\xi_3\eta_1+2(\xi_1\eta_1+\xi_3\eta_3)\big|$. To this end, Whitney type
decomposition is used.

\begin{definition}[Whitney type decomposition] (see Def 2 in \cite{Kinoshita19})
Let $ 2^{25}\leq A$ be dyadic and
$$H_1(\zeta_1,\zeta_2):=\big|\xi_1\xi_2(\xi_1+\xi_2)+\eta_1\eta_2(\eta_1+\eta_2)\big|,$$
$$H_2(\zeta_1,\zeta_2):=\big|\xi_1\eta_2+\xi_2\eta_1+2(\xi_1\eta_1+\xi_2\eta_2)\big|.$$
Define
$$Z^1_{A}:=\big\{(k_1,k_2)\in \mathbb{Z}^2\times\mathbb{Z}^2\ \big| \
H_1(\zeta_1,\zeta_2)\geq A^{-1}N_1^3, \ \ \   \forall \ (\xi_{\iota},\eta_{\iota})\in
{\mathcal T}^{A}_{k_{\iota}},\ \ {\iota}=1,2\big\},$$
$$Z^2_{A}:=\big\{(k_1,k_2)\in \mathbb{Z}^2\times\mathbb{Z}^2\ \big| \
H_2(\zeta_1,\zeta_2)\geq A^{-1}N_1^2, \ \ \   \forall \
(\xi_{\iota},\eta_{\iota})\in{\mathcal T}^{A}_{k_{\iota}},\ \ {\iota}=1,2\big\},$$
and
$$Z_{A}:=Z^1_{A}\cup
Z^2_{A}\subset\mathbb{Z}^2\times\mathbb{Z}^2,\hspace{5mm}R_A:=\bigcup\limits_{(k_1,k_2)\in
Z_{A}}{\mathcal T}^{A}_{k_1}\times{\mathcal
T}^{A}_{k_2}\subset\mathbb{R}^2\times\mathbb{R}^2.$$
\end{definition}
It is easy to see that $A_1\leq A_2 \Rightarrow  R_{A_1}\subset R_{A_2}.$
We further write
 \begin{equation}
   Q_A=\left\{
    \begin{aligned}
    & R_A\setminus R_{\frac{A}{2}}      \hspace{10mm}  for \ \ A\geq 2^{25} ,   \quad \\
    &R_{2^{25}} \hspace{16.9mm}   for \ \ A= 2^{25}, \nonumber
    \end{aligned}
    \right.
   \end{equation}
and
$$\widetilde{Z}_{A}:=\big\{(k_1,k_2)\in Z_A\ \big| \
{\mathcal T}^{A}_{k_1}\times{\mathcal T}^{A}_{k_2}\cap Q_A\not=\emptyset
\big\},$$
$$\widehat{Z}_{A_0}:=\big\{(k_1,k_2)\in \mathbb{Z}^2\times\mathbb{Z}^2\ \big| \
{\mathcal T}^{A_0}_{k_1}\times{\mathcal T}^{A_0}_{k_2}\cap Q_A=\emptyset  \ \ \   for \ \ \
A\leq A_0\big\}.$$
Then
$$Q_A=\bigcup\limits_{(k_1,k_2)\in\widetilde{Z}_{A}}{\mathcal T}^{A_0}_{k_1}\times{\mathcal
T}^{A_0}_{k_2},$$
and
$$\mathbb{R}^2\times\mathbb{R}^2=\bigcup\limits_{2^{25}\leq A\leq A_0}Q_A \ \cup \
(R_{A_0})^c.$$

Above we choose two functions $H_1$ and $H_2$ to define Whitney decomposition for the
purpose of almost orthogonality. In other words, since the intersection of $H_1$ and $H_2$
is comparable, we can sum up those tiles without loss.
\begin{lemma} \label{orthogonality}
Let $A,A_0$ be dyadic and $ 2^{25}\leq A\leq A_0$ . In the assumption of Case 1.2, for fixed
$k_1\in \mathbb{Z}^2$ , we have
$$\#\big\{k_2\in \mathbb{Z}^2 \ \big| \  (k_1,k_2)\in\widetilde{Z}_{A}\big \}\lesssim 1,$$
and
$$\#\big\{k_2\in \mathbb{Z}^2 \ \big| \  (k_1,k_2)\in\widehat{Z}_{A_0}\big \}\lesssim 1.$$
\end{lemma}
{\bf Proof.} We just show the first conclusion because the second one can be dealt with in a
similar way.

As $k_1$ is fixed, for every $k_2$  satisfying $(k_1,k_2)\in\widetilde{Z}_{A}$, we can
always find  $\tilde k_1$ and $\tilde k_2$ such that
$${\mathcal T}^{A}_{k_1}\subset{\mathcal T}^{A/2}_{\tilde k_1},  \ \
{\mathcal T}^{A}_{k_2}\subset{\mathcal T}^{A/2}_{\tilde k_2} \ \ \ and \ \ \ (\tilde
k_1,\tilde k_2)\notin\widetilde{Z}_{A/2}
$$
which imply that $\exists (\tilde \xi_1, \tilde \eta_1), ({\tilde \xi}'_1, {\tilde
\eta}'_1)\in {\mathcal T}^{A/2}_{\tilde k_1}$
and
$\exists (\tilde \xi_2, \tilde \eta_2), ({\tilde \xi}'_2, {\tilde \eta}'_2)\in {\mathcal
T}^{A/2}_{\tilde k_2}$  satisfying
$$H_1(\tilde \zeta_1,\tilde \zeta_2)\leq 2A^{-1}N_1^3 \ \ \ \ and \ \ \ \
H_2({\tilde \zeta}'_1,\tilde \zeta'_2)\leq 2A^{-1}N_1^2.$$

Let $(\xi_1,\eta_1)$ be the center of ${\mathcal T}^{A}_{k_1}$. It suffices to show that
there exist $k_{2,\iota}\in \mathbb{Z}^2 \ \ (\iota=1,2,3,4)$ such that
\begin{align}
 \big\{(\xi'_2,\eta'_2) \ \big| \ H_1(\zeta_1, \zeta'_2)\leq 2^4A^{-1}N_1^3  \ \
H_2(\zeta_1, \zeta'_2)\leq 2^4A^{-1}N_1^2
 \big\}\subset \bigcup\limits_{1\leq \iota\leq 4}{\mathcal T}^{c_0A}_{k_{2,\iota}},
\label{almostorto}
\end{align}
where $c_0$ is constant, for example we can choose $c_0=2^{-20}$.

$\ H_1(\zeta_1, \zeta'_2)\leq2^4A^{-1}N_1^3$ and $H_2(\zeta_1, \zeta'_2)\leq 2^4A^{-1}N_1^2
$ mean that
$$H_1(\zeta_1,\zeta'_2):=\big|\xi_1\xi'_2(\xi_1+\xi'_2)+\eta_1\eta'_2(\eta_1+\eta'_2)\big|\leq
2^4A^{-1}N_1^3,$$
and
$$H_2(\zeta_1,\zeta'_2):=\big|\xi_1\eta'_2+\xi'_2\eta_1+2(\xi_1\eta_1+\xi'_2\eta'_2)\big|\leq
2^4A^{-1}N_1^2.$$
By setting $\xi_2=\xi'_2+\frac{\xi_1}{2}$ and $\eta_2=\eta'_2+\frac{\eta_1}{2}$, those two
inequalities are equivalent to
\begin{align}
\widetilde
H_1(\zeta_2):&=\bigg|\xi_1\xi^2_2+\eta_1\eta^2_2-\frac{\xi^3_1+\eta^3_1}{4}\bigg|\leq
2^4A^{-1}N_1^3,       \label{widetilde H1} \\
  \widetilde
H_2(\zeta_2):&=\bigg|\frac{3}{2}\xi_1\eta_1+2\xi_2\eta_2\bigg|\leq2^4A^{-1}N_1^2.
\label{widetilde H2}
\end{align}
Since $|(\xi_1,\eta_1)|\sim N_1$ and $|(\xi'_2,\eta'_2)|\sim N_2\sim N_1$, we denote
$\zeta_3{''}=-\zeta_1-\zeta'_2$, then $|\zeta_3{''}|\sim N_3\ll N_1$  and
$$\big|\xi_2\big|=\bigg|\frac{\xi_1+2\xi'_2}{2}\bigg|=\bigg|\frac{\xi'_2-\xi_3{''}}{2}\bigg|
\sim N_1, \ \ \ |\eta_2|\lesssim N_3.$$
Thus, from \eqref{widetilde H2} we conclude
 \begin{align}
\bigg|\eta_2+\frac{3\xi_1\eta_1}{4\xi_2}\bigg|\leq2^8A^{-1}N_1.  \label{widetilde H2a}
\end{align}

\eqref{widetilde H1} and  \eqref{widetilde H2a} yield
 \begin{align}
 &\bigg|\xi_1\xi^2_2+\eta_1\frac{9\xi^2_1\eta^2_1}{16\xi^2_2}-\frac{\xi^3_1+\eta^3_1}{4}\bigg|
\notag \\
 \leq& \bigg|\xi_1\xi^2_2+\eta_1\eta^2_2-\frac{\xi^3_1+\eta^3_1}{4}\bigg|+\big|\eta_1\big|
\big|\eta_2+\frac{3\xi_1\eta_1}{4\xi_2}\big| \big|\eta_2-\frac{3\xi_1\eta_1}{4\xi_2}\big|
\notag \\
\leq&2^5A^{-1}N_1^3.\label{widetilde H2b}
\end{align}
We observe that
\begin{align}
 &\bigg|\frac{\partial}{\partial\xi_2}\bigg(\xi_1\xi^2_2+\eta_1\frac{9\xi^2_1\eta^2_1}
{16\xi^2_2}-\frac{\xi^3_1+\eta^3_1}{4}\bigg)\bigg| \notag \\
 =&\bigg|\frac{2\xi_1}{\xi_2^3}(\xi_2^4-\frac{9\xi_1\eta^3_1}{16})\bigg|\notag \\
 \geq& 2^5N_1^2. \nonumber
\end{align}
Since
$G(\xi_2)=\xi_1\xi^2_2+\eta_1\frac{9\xi^2_1\eta^2_1}{16\xi^2_2}-\frac{\xi^3_1+\eta^3_1}{4}=C$
has at most four roots, there exist at most four constants $c_{\iota}\ (\iota=1,2,3,4)$ such
that
\begin{align}
|\xi_2-c_{\iota}|\leq 2^{15}A^{-1}N_1  \label{widetilde H2b2x}
\end{align}
for any $\xi_2$ satisfying \eqref{widetilde H2b}.

Moreover, from \eqref{widetilde H2a} we know that
\begin{align}
\bigg|\frac{\partial}{\partial\xi_2}\bigg(\eta_2+\frac{3\xi_1\eta_1}{4\xi_2}\bigg)\bigg|=
\bigg|\frac{3\xi_1\eta_1}{4\xi^2_2} \bigg|\sim \frac{N_3}{N_1}\ll 1. \label{widetilde H2b2y}
\end{align}
Thus \eqref{widetilde H2b2x} and \eqref{widetilde H2b2y} conclude \eqref{almostorto}. This
finish the proof.
\begin{proposition}\label{case12dprop1}In the assumption of Case 1.2, we have
\begin{align}
\bigg|\iint_{\Omega_r^{*}}\widehat{v}_{1}(\lambda_1)\widehat{v}_{2}(\lambda_2)\widehat{v}_{3}
(\lambda_3)\bigg|
\lesssim
N_1^{-1}N_3^{-\frac{1}{4}}(L_1L_2L_3)^{\frac{5}{12}}\|\widehat{v}_{1}\|_{L^{2}}\|\widehat{v}
_{2}\|_{L^{2}}
\|\widehat{v}_{3}\|_{L^{2}}.\nonumber
\end{align}
\end{proposition}
{\bf Proof.} By Whitney decomposition, we have
 \begin{align}
&\bigg|\iint_{\Omega_r^{*}}\widehat{v}_{1}(\lambda_1)\widehat{v}_{2}(\lambda_2)\widehat{v}
_{3}(\lambda_3)\bigg|\notag \\
\lesssim&\sum_{ A\leq A_0}\sum\limits_{(k_1,k_2)\in\widetilde{Z}_{A}}
\bigg|\int_{*}\widehat{v}_{1}|_{{\widetilde{\mathcal
T}}^{A}_{k_1}}(\lambda_1)\widehat{v}_{2}|_{{\widetilde{\mathcal
T}}^{A}_{k_2}}(\lambda_2)\widehat{v}_{3}(\lambda_3)\bigg|\notag \\
&\ \ +\sum\limits_{(k_1,k_2)\in\widehat{Z}_{A_0}}
\bigg|\int_{*}\widehat{v}_{1}|_{{\widetilde{\mathcal
T}}^{A_0}_{k_1}}(\lambda_1)\widehat{v}_{2}|_{{\widetilde{\mathcal
T}}^{A_0}_{k_2}}(\lambda_2)\widehat{v}_{3}(\lambda_3)\bigg|\notag \\
:=&Y_1+Y_2.\nonumber
\end{align}

Using Lemma \ref{case12b}, Lemma \ref{case12ca1} and Lemma \ref{orthogonality} to estimate
the first term
\begin{align}
Y_1= &\sum_{ A\leq A_0}\sum\limits_{(k_1,k_2)\in\widetilde{Z}_{A}}
\bigg|\int_{*}\widehat{v}_{1}|_{{\widetilde{\mathcal
T}}^{A}_{k_1}}(\lambda_1)\widehat{v}_{2}|_{{\widetilde{\mathcal
T}}^{A}_{k_2}}(\lambda_2)\widehat{v}_{3}(\lambda_3)\bigg|\notag  \\
 \lesssim & \sum_{ A\leq
A_0}A^{\frac{1}{2}}N_1^{-\frac{3}{2}}N_3^{-\frac{1}{2}}(L_1L_2L_3)^{\frac{1}{2}}
 \sum\limits_{k_1\sim k_2}\big\|\widehat{v}_{1}|_{{\widetilde{\mathcal T}}^{A}_{k_1}}
\big\|_{L^{2}}\big\|\widehat{v}_{2}|_{{\widetilde{\mathcal
T}}^{A}_{k_2}}\big\|_{L^{2}}\big\|\widehat{v}_{3}\big\|_{L^{2}} \notag  \\
  \lesssim &A_0^{\frac{1}{2}}N_1^{-\frac{3}{2}}N_3^{-\frac{1}{2}}(L_1L_2L_3)^{\frac{1}{2}}
\|\widehat{v}_{1}\|_{L^{2}} \|\widehat{v}_{2}\|_{L^{2}}
\|\widehat{v}_{3}\|_{L^{2}}.\label{case12dprop1xc1}
\end{align}

By Lemma \ref{case12a} and Lemma \ref{orthogonality}, we can control the second term
\begin{align}
Y_2= &\sum\limits_{(k_1,k_2)\in\widehat{Z}_{A_0}}
 \bigg|\int_{*}\widehat{v}_{1}|_{{\widetilde{\mathcal
 T}}^{A_0}_{k_1}}(\lambda_1)\widehat{v}_{2}|_{{\widetilde{\mathcal
 T}}^{A_0}_{k_2}}(\lambda_2)\widehat{v}_{3}(\lambda_3)\bigg|\notag  \\
  \lesssim & (A_0N_1L_0)^{-\frac{1}{2}}(L_1L_2L_3)^{\frac{1}{2}}
  \sum\limits_{k_1\sim k_2}\big\|\widehat{v}_{1}|_{{\widetilde{\mathcal T}}^{A_0}_{k_1}}
 \big\|_{L^{2}}\big\|\widehat{v}_{2}|_{{\widetilde{\mathcal
 T}}^{A_0}_{k_2}}\big\|_{L^{2}}\big\|\widehat{v}_{3}\big\|_{L^{2}} \notag  \\
 \lesssim &(A_0N_1L_0)^{-\frac{1}{2}}(L_1L_2L_3)^{\frac{1}{2}}
 \|\widehat{v}_{1}\|_{L^{2}} \|\widehat{v}_{2}\|_{L^{2}}
\|\widehat{v}_{3}\|_{L^{2}},\label{case12dprop1xc2}
 \end{align}
where $L_0=\max_{1\leq {\iota}\leq3}\{L_{\iota}\}$.

Choosing $A_0=N_1N_3^{-\frac{1}{2}}L_0^{-\frac{1}{2}}$, thus from \eqref{case12dprop1xc1}
and \eqref{case12dprop1xc2} we have
\begin{align}
&\bigg|\iint_{\Omega_r^{*}}\widehat{v}_{1}(\lambda_1)\widehat{v}_{2}(\lambda_2)\widehat{v}
_{3}(\lambda_3)\bigg|\notag
\\
\lesssim&L_0^{-\frac{1}{4}}N_1^{-1}N_3^{-\frac{1}{4}}(L_1L_2L_3)^{\frac{1}{2}}
 \|\widehat{v}_{1}\|_{L^{2}} \|\widehat{v}_{2}\|_{L^{2}} \|\widehat{v}_{3}\|_{L^{2}} \notag
\\
 \lesssim&N_1^{-1}N_3^{-\frac{1}{4}}(L_1L_2L_3)^{\frac{5}{12}}
 \|\widehat{v}_{1}\|_{L^{2}} \|\widehat{v}_{2}\|_{L^{2}} \|\widehat{v}_{3}\|_{L^{2}},
\nonumber
\end{align}
which concludes the proof.

Hence, under the assumption of Case 1.2, from Proposition \ref{case12dprop1},
we obtain that
\begin{align}
&\bigg|\iint_{\Omega_r^{*}}\frac{\sum^3_{{\iota}=1}m^2(\zeta_{\iota})(\xi_{\iota}+\eta
_{\iota})}
{m(\zeta_1)m(\zeta_2)m(\zeta_3)}\widehat{v}_{1}(\lambda_1)\widehat{v}_{2}(\lambda_2)
\widehat{v}_{3}(\lambda_3)\bigg|\notag \\
\lesssim
&N^{-\frac{1}{4}}(L_1L_2L_3)^{\frac{5}{12}}\|\widehat{v}_{1}\|_{L^{2}}
\|\widehat{v}_{2}\|_{L^{2}}
\|\widehat{v}_{3}\|_{L^{2}}.\label{maintricase12}
\end{align}

Now let's consider the high-high interactions case.  By symmetry, we can assume that
$|\xi_2|\geq |\eta_2|$.
\\{\bf Case 2} \hspace{2mm} $N_1\sim N_2\sim N_3$ and $|\xi_2|\geq |\eta_2|$

We shall investigate these subcases of Case 2 relying on the relationship between $|\xi_2|$
and $|\xi_1-\xi_3|$.

First, if $|\xi_2|\gg |\xi_1-\xi_3|$, then $|\xi_1|\sim|\xi_2|\sim|\xi_3|$, thus
$|\xi_{\iota}|\sim|\eta_{\iota}|\sim N_1$  for ${\iota}=1,2,3$. Otherwise,
$$|\xi_1\xi_2\xi_3+\eta_1\eta_2\eta_3|\sim|\xi_1\xi_2\xi_3|\sim N_1^3$$
which contradicts to the low modulation assumption $|\xi_1\xi_2\xi_3+\eta_1\eta_2\eta_3|\ll
\gamma_0N_1N_2N_3$.
\\{\bf Case 2.1} \hspace{2mm} $|\xi_2|\gg |\xi_1-\xi_3|$ \ and \
$|\xi_{\iota}|\sim|\eta_{\iota}|\sim N_1$  for $\iota=1,2,3$

If $\xi_1\xi_3<0$, then $|\xi_1-\xi_3|=|\xi_1|+|\xi_3|\sim N_1$ which contradicts to
$|\xi_1-\xi_3|\ll |\xi_2|$, thus $\xi_1\xi_3>0$. Without loss of generality, we assume that
$\xi_1>0$ and $\xi_3>0$. From $\xi_1+\xi_2+\xi_3=0$, we know that $\xi_2<0$, thereby
$\eta_1\eta_2\eta_3>0$. Otherwise, $\eta_1\eta_2\eta_3<0$ gives
$$|\xi_1\xi_2\xi_3+\eta_1\eta_2\eta_3|=|\xi_1\xi_2\xi_3|+|\eta_1\eta_2\eta_3|\sim N_1^3.$$
This also contradicts to the  low modulation assumption.
\begin{table}[]
\centering
 \caption{Subcases of 2.1 and the corresponding controllable quantities}
 \begin{tabular}{|l|c|c|c|c|c|c|c|c|c|}
 \hline
 \hspace{4mm}Case                                          & \multicolumn{3}{c|}{Case 2.1.1}
& \multicolumn{3}{c|}{Case 2.1.2}
& \multicolumn{3}{c|}{Case 2.1.3}
\\ \hline
 $\hspace{7mm}\iota$                                       & \multicolumn{3}{c|}{1
\hspace{9mm}   2   \hspace{8mm}    3}
& \multicolumn{3}{c|}{1  \hspace{8mm}   2   \hspace{8mm}    3}
& \multicolumn{3}{c|}{1  \hspace{8mm}   2   \hspace{8mm}    3}
\\ \hline
 $\hspace{7mm}\xi_{\iota}$                                 & \hspace{2mm}$+$    \hspace{2mm}
&\hspace{2mm} $-$     \hspace{2mm}                         &\hspace{2mm} $+$   \hspace{2mm}
& \hspace{2mm}$+$    \hspace{2mm}                      &\hspace{2mm} $-$     \hspace{2mm}
&\hspace{2mm} $+$   \hspace{2mm}                                 & \hspace{2mm}$+$
\hspace{2mm}                      &\hspace{2mm} $-$     \hspace{2mm}
&\hspace{2mm} $+$   \hspace{2mm}                      \\ \hline
 $\hspace{7mm}\eta_{\iota}$                                & \hspace{2mm}$+$    \hspace{2mm}
&\hspace{2mm} $-$     \hspace{2mm}                         &\hspace{2mm} $-$   \hspace{2mm}
& \hspace{2mm}$-$    \hspace{2mm}                      &\hspace{2mm} $+$     \hspace{2mm}
&\hspace{2mm} $-$   \hspace{2mm}                                 & \hspace{2mm}$-$
\hspace{2mm}                      &\hspace{2mm} $-$     \hspace{2mm}
&\hspace{2mm} $+$   \hspace{2mm}                         \\ \hline
 $H_2$ or  $H_3$                                 &
\multicolumn{3}{c|}{$H_3(\zeta_1,\zeta_3)\sim N_1^2$} &
\multicolumn{3}{l|}{\hspace{3mm}$H_2(\zeta_1,\zeta_3)\sim N_1^2$} &
\multicolumn{3}{l|}{$\hspace{3mm}H_3(\zeta_1,\zeta_3)\sim N_1^2$} \\ \hline
\end{tabular}
 \end{table}
From the discussion above, we can divide Case 2.1 into the following three subcases:
\\{ Case 2.1.1} \hspace{2mm} $\xi_1>0 \ \ \xi_2<0 \ \ \xi_3>0$ \ and \ $\eta_1>0 \ \
\eta_2<0 \ \ \eta_3<0$,
\\{ Case 2.1.2} \hspace{2mm} $\xi_1>0 \ \ \xi_2<0 \ \ \xi_3>0$ \ and \ $\eta_1<0 \ \
\eta_2>0 \ \ \eta_3<0$,
\\{ Case 2.1.3} \hspace{2mm} $\xi_1>0 \ \ \xi_2<0 \ \ \xi_3>0$ \ and \ $\eta_1<0 \ \
\eta_2<0 \ \ \eta_3>0$.

The most important thing is that in each subcase we find one of $|\xi_1\eta_2-\xi_2\eta_1|$
and $|\xi_1\eta_2+\xi_2\eta_1+2(\xi_1\eta_1+\xi_2\eta_2)|$ is controllable ( see Table 1 where we write $H_3(\zeta_1,\zeta_3)=|\xi_1\eta_3-\xi_3\eta_1|$ for simplicity).
Thus it will be not too complicated to use the nonlinear Loomis--Whitney inequality.

As $\xi_1+\xi_2+\xi_3=0$ and $\eta_1+\eta_2+\eta_3=0$, simple computation implies
$$|\xi_1\eta_2-\xi_2\eta_1|=|\xi_1\eta_3-\xi_3\eta_1|=|\xi_2\eta_3-\xi_3\eta_2|$$
and
\begin{align}
&|\xi_1\eta_2+\xi_2\eta_1+2(\xi_1\eta_1+\xi_2\eta_2)|  \notag \\
=&|\xi_1\eta_3+\xi_3\eta_1+2(\xi_1\eta_1+\xi_3\eta_3)|  \notag \\
=&|\xi_2\eta_3+\xi_3\eta_2+2(\xi_2\eta_2+\xi_3\eta_3)|.  \nonumber
\end{align}

To avoid being lengthy and tedious, we just take Case 2.1.3 for instance to explain how
$|\xi_1\eta_3-\xi_3\eta_1|$ can be controlled.
When $\xi_1>0, \xi_2<0, \xi_3>0$  and  $\eta_1<0 , \eta_2<0, \eta_3>0$,
$$|\xi_1\eta_3-\xi_3\eta_1|=|\xi_1\eta_3|+|\xi_3\eta_1|\sim N_1^2.$$

Now we will deal with these three subcases respectively.
 \\{\bf Case 2.1.1} \hspace{2mm} $\xi_1>0 \ \ \xi_2<0 \ \ \xi_3>0$ \ and \ $\eta_1>0 \ \
\eta_2<0 \ \ \eta_3<0$
\begin{proposition}\label{case211dprop1}In the assumption of Case 2.1.1, we have
 \begin{align}
 \bigg|\iint_{\Omega_r^{*}}\widehat{v}_{1}(\lambda_1)\widehat{v}_{2}(\lambda_2)\widehat{v}
_{3}(\lambda_3)\bigg|
 \lesssim
N_1^{-\frac{5}{4}}(L_1L_2L_3)^{\frac{5}{12}}\|\widehat{v}_{1}\|_{L^{2}}\|\widehat{v}_{2}\|
_{L^{2}}
 \|\widehat{v}_{3}\|_{L^{2}}.\nonumber
 \end{align}
 \end{proposition}
 {\bf Proof.} In fact, the proof is very similar to that of Proposition \ref{case12dprop1}.
As $|\xi_1\eta_3-\xi_3\eta_1|\sim N_1^2$, it is easy to verify that Lemma
\ref{case12a}--\ref{case12ca1} also hold true under the assumption. Hence, it suffices to
verify Lemma \ref{orthogonality}.

In order to avoid confusion and misunderstanding, we use new variable when we check those
condition provided in Lemma \ref{orthogonality}. For example, here we use
$\overline\xi_{\iota}$ to replace $\xi_{\iota}$ in Lemma \ref{orthogonality}.

Set 
$${\overline\xi}'_2=\xi_1, \ \ \overline\xi_1=\xi_3,\ \
\overline\xi_2={\overline\xi}'_2+\frac{\overline\xi_1}{2},$$ 
and
$${\overline\eta}'_2=\eta_1,\ \ \overline\eta_1=\eta_3,\ \
\overline\eta_2={\overline\eta}'_2+\frac{\overline\eta_1}{2},$$
then

$$\overline\xi_2=\xi_1+\frac{\xi_3}{2}=\frac{\xi_1-\xi_2}{2}>0, \ \ |\overline\xi_2|\sim
N_1,$$
$$\overline\eta_2=\eta_1+\frac{\eta_3}{2}=\frac{\eta_1-\eta_2}{2}>0, \ \
|\overline\eta_2|\sim N_1.$$
Note that
\begin{align}
 &\bigg|\frac{\partial}{\partial\overline\xi_2}\bigg(\overline\xi_1\overline\xi^2_2
+\overline\eta_1\frac{9\overline\xi^2_1\overline\eta^2_1}{16\overline\xi^2_2}-
\frac{\overline\xi^3_1+\overline\eta^3_1}{4}\bigg)\bigg|
\notag \\
 =&\bigg|\frac{2\overline\xi_1}{\overline\xi_2^3}\bigg(\overline\xi_2^4-
\frac{9\overline\xi_1\overline\eta^3_1}{16}\bigg)\bigg|\notag \\
 =&\bigg|\frac{2\xi_3}{\overline\xi_2^3}\bigg|\cdot\bigg(\bigg|\overline\xi_2^4\bigg|+
 \bigg|\frac{9\xi_3\eta^3_3}{16}\bigg|\bigg)\notag \\
\gtrsim& N_1^2, \nonumber
\end{align}
and
\begin{align}
\bigg|\frac{\partial}{\partial\overline\xi_2}\bigg(\overline\eta_2+\frac{3\overline\xi_1
\overline\eta_1}{4\overline\xi_2}\bigg)\bigg|=\bigg|\frac{3\overline\xi_1\overline\eta_1}
{4\overline\xi^2_2}\bigg| =\bigg|\frac{3\xi_3\eta_3} {4\overline\xi^2_2}\bigg|\sim 1.
\nonumber
\end{align}
Thus Lemma \ref{orthogonality} holds true under the assumption of Case 2.1.1. We finish the
proof.
\\{\bf Case 2.1.2} \hspace{2mm} $\xi_1>0 \ \ \xi_2<0 \ \ \xi_3>0$ \ and \ $\eta_1<0 \ \
\eta_2>0 \ \ \eta_3<0$

At this subcase, $|\xi_1\eta_2+\xi_2\eta_1+2(\xi_1\eta_1+\xi_2\eta_2)|\sim N_1^2$, we shall
take the same steps as in Case 1.1.

In order to control $|\xi_1\eta_2-\xi_2\eta_1|$ locally, we need the angular decomposition:
\begin{align}
\mathbb{R}^2\times\mathbb{R}^2=&\bigcup\limits_{\substack{0\leq j_1,j_2\leq M-1 \\
|j_1-j_2|\leq16} }
{\mathcal D}^{M}_{j_1}\times {\mathcal D}^{M}_{j_2}\ \cup
\bigcup\limits_{64\leq A\leq M}\bigcup\limits_{\substack{0\leq j_1,j_2\leq A-1 \\ 16
\leq|j_1-j_2|\leq32} }{\mathcal D}^{A}_{j_1}\times {\mathcal D}^{A}_{j_2}.\nonumber
\end{align}
Note that $|\xi_1\eta_2-\xi_2\eta_1|\sim A^{-1}N_1^2$ in ${\mathcal D}^{A}_{j_1}\times
{\mathcal D}^{A}_{j_2}$ when $|j_1-j_2|\leq 2^5$.

Let's first treat $(\xi_1,\eta_1)\times(\xi_2,\eta_2)\notin \mathcal{I}_1$.
 \begin{proposition}\label{case212mathcalI1}
 Suppose that $(\xi_1,\eta_1)\times(\xi_2,\eta_2)\notin \mathcal{I}_1$ , in the assumption
of Case 2.1.2£¬ we have
 \begin{align}
 \bigg|\iint_{\Omega_r^{*}}\widehat{v}_{1}(\lambda_1)\widehat{v}_{2}(\lambda_2)
\widehat{v}_{3}(\lambda_3)\bigg|
 \lesssim \gamma^{0+}_0N_1^{-2+}(L_1L_2L_3)^{\frac{1}{2}}\|\widehat{v}_{1}\|_{L^{2}}\|
\widehat{v}_{2}\|_{L^{2}}
 \|\widehat{v}_{3}\|_{L^{2}}.\nonumber
 \end{align}
 \end{proposition}
 {\bf Proof.} If $(\xi_1,\eta_1)\times(\xi_2,\eta_2)\in {\mathcal D}^{A}_{j_1}\times
{\mathcal D}^{A}_{j_2} $ for $|j_1-j_2|\leq 2^{5}$ and
$(\xi_1,\eta_1)\times(\xi_2,\eta_2)\notin \mathcal{I}_1$,
we claim that
 $$\max\{L_1,L_2,L_3\}\gg A^{-1}N_1^3, \ \ \  \ \ \ A\gg\gamma_0^{-1}.$$

Since $|\xi_{\iota}|\sim|\eta_{\iota}|\sim N_1 \  (\iota=1,2)$ and
$\xi_1>0,\eta_1<0,\xi_2<0 ,\eta_2>0$, we deduce
$$|(\cos\theta_1,\sin\theta_1)+(\cos\theta_2,\sin\theta_2)|\leq2^8A^{-1}$$
from $|j_1-j_2|\leq 32$.

If $\max\{L_1,L_2,L_3\}\lesssim A^{-1}N_1^3$, hence from Lemma \ref{mathcalI1a0}
 we have $N_3\lesssim A^{-1}N_1\ll N_1$. This contradicts to $N_3\sim N_1$.
 Thus
$$\gamma_0N_1^3\gg\max\{L_1,L_2,L_3\}\gg A^{-1}N_1^3,$$
which leads to $A\gg\gamma_0^{-1}$. The claim follows.

According to the angular decomposition, Lemma \ref{mathcalI1a} and the claim, we have
 \begin{align}
 &\bigg|\iint_{\Omega_r^{*}}\widehat{v}_{1}(\lambda_1)\widehat{v}_{2}(\lambda_2)
\widehat{v}_{3}(\lambda_3)\bigg|\notag \\
 \lesssim&\sum_{\gamma_0^{-1}\ll A}
\sum\limits_{\substack{0\leq j_1,j_2\leq A-1 \\ 16\leq|j_1-j_2|\leq32}
}\bigg|\int_{*}\widehat{v}_{1}|_{{\widetilde{\mathcal
D}}^{A}_{j_1}}(\lambda_1)\widehat{v}_{2}|_{{\widetilde{\mathcal
D}}^{A}_{j_2}}(\lambda_2)\widehat{v}_{3}(\lambda_3)\bigg|\notag \\
 \lesssim& \sum_{\gamma_0^{-1}\ll A}
\sum_{j_1\sim j_2}A^{-\frac{1}{2}}N_1^{-\frac{1}{2}}A^{\frac{1}{2}-}N_1^{-\frac{3}{2}+}
(L_1L_2L_3)^{\frac{1}{2}}\big\|\widehat{v}_{1}|_{{\widetilde{\mathcal
D}}^{A}_{j_1}}\big\|_{L^{2}}\big\|\widehat{v}_{2}|_{{\widetilde{\mathcal
D}}^{A}_{j_2}}\big\|_{L^{2}}\big\|\widehat{v}_{3}\big\|_{L^{2}}\notag \\
 \lesssim&\gamma_0^{0+}N_1^{-2+}(L_1L_2L_3)^{\frac{1}{2}}\|\widehat{v}_{1}\|_{L^{2}}
\|\widehat{v}_{2}\|_{L^{2}}\|\widehat{v}_{3}\|_{L^{2}},\nonumber
 \end{align}
which concludes the proof.

Next we consider $(\xi_1,\eta_1)\times(\xi_2,\eta_2)\in \mathcal{I}_1$. We make twice
decomposition for $\mathcal{I}_1$:
 \begin{align}
 \mathcal{I}_1=\bigcup\limits_{2^{11}\leq M}\mathcal{I}^M_1\nonumber
 \end{align}
where
\begin{align}
&\mathcal{I}^M_1:=\big({\mathcal D}^{M}_{\frac{3}{4}M}\times
{\mathcal D}^{M}_{\frac{3}{4}M}\big)\setminus \big({\mathcal
D}^{2M}_{\frac{3}{2}M}\times{\mathcal D}^{2M}_{\frac{3}{2}M}\big), \nonumber
\end{align}
and
\begin{align}
\mathcal{I}^M_1=
\bigcup\limits_{2^{30}M\leq A}\bigcup\limits_{\substack{0\leq j_1,j_2\leq A-1 \\
16\leq|j_1-j_2|\leq32} }\big({\mathcal D}^{A}_{j_1}\times {\mathcal
D}^{A}_{j_2}\big)\cap\mathcal{I}^M_1 .\nonumber
\end{align}
\begin{proposition}\label{case212mathcalI1in}
Assume that $(\xi_1,\eta_1)\times(\xi_2,\eta_2)\in \mathcal{I}_1$ , then we have
\begin{align}
&\bigg|\iint_{\Omega_r^{*}}(\xi_3+\eta_3)\widehat{v}_{1}(\lambda_1)\widehat{v}_{2}(\lambda_2)
\widehat{v}_{3}(\lambda_3)\bigg|\notag \\
\lesssim & N_1^{-1+}(L_1L_2L_3)^{\frac{1}{2}}\|\widehat{v}_{1}\|_{L^{2}}
\|\widehat{v}_{2}\|_{L^{2}}
\|\widehat{v}_{3}\|_{L^{2}}.\nonumber
\end{align} \end{proposition}
{\bf Proof.} Since $(\xi_1,\eta_1)\times(\xi_2,\eta_2)\in \mathcal{I}_1^M$,
$$|\xi_3+\eta_3|\lesssim |\xi_1+\eta_1|+|\xi_2+\eta_2|\lesssim M^{-1}N_1.$$
We also observe that
$$\max\{L_1,L_2,L_3\}\gg A^{-1}N_1^3.$$
Because $(\xi_1, \eta_1)$ and $(\xi_2, \eta_2)$ locate in the second and the forth quadrant
respectively, one has
$$|(\cos\theta_1,\sin\theta_1)+(\cos\theta_2,\sin\theta_2)|\leq2^8A^{-1}$$
for $|j_1-j_2|\leq 32$.

As
\begin{align}
&|\xi_1\xi_2(\xi_1+\xi_2)+\eta_1\eta_2(\eta_1+\eta_2)|\notag \\
=&\big|r_1r_2(r_1-r_2)(\cos^3\theta_1+\sin^3\theta_1)-r_1r_2(r_1-r_2)\cos^2\theta_1
(\cos\theta_1+\cos\theta_2)\notag \\
&-r_1r_2(r_1-r_2)\sin^2\theta_1(\sin\theta_1+\sin\theta_2)-r_1r^2_2\cos\theta_1
\cos\theta_2(\cos\theta_1+\cos\theta_2)\notag \\
&-r_1r^2_2\sin\theta_1\sin\theta_2(\sin\theta_1+\sin\theta_2)\big|\notag \\
\gtrsim &M^{-1}r_1r_2|r_1-r_2|-A^{-1}r_1r^2_2, \nonumber
\end{align}
 if $\max\{L_1,L_2,L_3\}\lesssim A^{-1}N_1^3$, then $|r_1-r_2|\lesssim MA^{-1}N_1$.

However
\begin{align}
  N_3\lesssim  |r_1-r_2|+A^{-1}N_1
  \lesssim MA^{-1}N_1+A^{-1}N_1 \ll N_1, \nonumber
 \end{align}
 which contradicts to $N_3\sim N_1$. Hence we verify the observation.

According to the twice decomposition of $\mathcal{I}_1$
, Lemma \ref{mathcalI1a} and the observation, we have
\begin{align}
&\bigg|\iint_{\Omega_r^{*}}(\xi_3+\eta_3)\widehat{v}_{1}(\lambda_1)\widehat{v}_{2}
(\lambda_2)\widehat{v}_{3}(\lambda_3)\bigg|\notag \\
\lesssim&\sum_{2^{11}\lesssim M }
\sum_{2^{30}M\leq A}\sum\limits_{\substack{(j_1,j_2)\in\mathcal{J}^{\mathcal{I}^M_1}_A \\
 16\leq|j_1-j_2|\leq32} }\bigg|\int_{*}(\xi_3+\eta_3)\widehat{v}_{1}|_{{\widetilde{\mathcal
D}}^{A}_{j_1}}(\lambda_1)\widehat{v}_{2}|_{{\widetilde{\mathcal
D}}^{A}_{j_2}}(\lambda_2)\widehat{v}_{3}(\lambda_3)\bigg|\notag \\
\lesssim&\sum_{2^{11}\lesssim M }\sum_{2^{30}M\leq A}
\sum_{j_1\sim j_2}M^{-1}N_1A^{0-}N_1^{-2+}
(L_1L_2L_3)^{\frac{1}{2}}\big\|\widehat{v}_{1}|_{{\widetilde{\mathcal
D}}^{A}_{j_1}}\big\|_{L^{2}}\big\|\widehat{v}_{2}|_{{\widetilde{\mathcal
D}}^{A}_{j_2}}\big\|_{L^{2}}\big\|\widehat{v}_{3}\big\|_{L^{2}}\notag \\
 \lesssim&N_1^{-1+}(L_1L_2L_3)^{\frac{1}{2}}\|\widehat{v}_{1}\|_{L^{2}}
\|\widehat{v}_{2}\|_{L^{2}}\|\widehat{v}_{3}\|_{L^{2}}.\nonumber
\end{align}
We finish the proof.

Therefore, under the assumption of Case 2.1.2, we obtain that
\begin{align}
&\bigg|\iint_{\Omega_r^{*}}(\xi_3+\eta_3)\widehat{v}_{1}(\lambda_1)\widehat{v}_{2}(\lambda_2)
\widehat{v}_{3}(\lambda_3)\bigg|\notag \\
\lesssim & N_1^{-1+}(L_1L_2L_3)^{\frac{1}{2}}\|\widehat{v}_{1}\|_{L^{2}}
\|\widehat{v}_{2}\|_{L^{2}}
\|\widehat{v}_{3}\|_{L^{2}}\label{maincase212}
\end{align}
from Proposition \ref{case212mathcalI1} and Proposition \ref{case212mathcalI1in}.
\\{\bf Case 2.1.3} \hspace{2mm} $\xi_1>0 \ \ \xi_2<0 \ \ \xi_3>0$ \ and \ $\eta_1<0 \ \
\eta_2<0 \ \ \eta_3>0$
 \begin{proposition}\label{case213dprop1}In the assumption of Case 2.1.3, we have
 \begin{align}
 \bigg|\iint_{\Omega_r^{*}}\widehat{v}_{1}(\lambda_1)\widehat{v}_{2}(\lambda_2)
\widehat{v}_{3}(\lambda_3)\bigg|
 \lesssim
N_1^{-\frac{5}{4}}(L_1L_2L_3)^{\frac{5}{12}}\|\widehat{v}_{1}\|_{L^{2}}
\|\widehat{v}_{2}\|_{L^{2}}
 \|\widehat{v}_{3}\|_{L^{2}}.\nonumber
 \end{align}
 \end{proposition}
 {\bf Proof.} As $|\xi_1\eta_3-\xi_3\eta_1|\sim N_1^2$, like Proposition
\ref{case211dprop1}, it suffices to
verify Lemma \ref{orthogonality}.

Let's denote
$${\overline\xi}'_2=\xi_2, \ \ \overline\xi_1=\xi_1,\ \
\overline\xi_2={\overline\xi}'_2+\frac{\overline\xi_1}{2},$$ and
$${\overline\eta}'_2=\eta_2,\ \ \overline\eta_1=\eta_1,\ \
\overline\eta_2={\overline\eta}'_2+\frac{\overline\eta_1}{2},$$
then
$$\overline\xi_2=\xi_2+\frac{\xi_1}{2}=\frac{\xi_2-\xi_3}{2}<0, \ \ |\overline\xi_2|\sim
N_1,$$
$$\overline\eta_2=\eta_2+\frac{\eta_1}{2}=\frac{\eta_2-\eta_3}{2}<0, \ \
|\overline\eta_2|\sim N_1.$$
Note that
\begin{align}
 &\bigg|\frac{\partial}{\partial\overline\xi_2}\bigg(\overline\xi_1\overline\xi^2_2
+\overline\eta_1\frac{9\overline\xi^2_1\overline\eta^2_1}{16\overline\xi^2_2}-
\frac{\overline\xi^3_1+\overline\eta^3_1}{4}\bigg)\bigg|
\notag \\
 =&\bigg|\frac{2\overline\xi_1}{\overline\xi_2^3}\bigg(\overline\xi_2^4-
\frac{9\overline\xi_1\overline\eta^3_1}{16}\bigg)\bigg|\notag \\
 =&\bigg|\frac{2\xi_1}{\overline\xi_2^3}\bigg|\cdot\bigg(\bigg|\overline\xi_2^4\bigg|+
 \bigg|\frac{9\xi_1\eta^3_1}{16}\bigg|\bigg)\notag \\
\gtrsim& N_1^2, \nonumber
\end{align}
and
\begin{align}
\bigg|\frac{\partial}{\partial\overline\xi_2}\bigg(\overline\eta_2+\frac{3\overline\xi_1
\overline\eta_1}{4\overline\xi_2}\bigg)\bigg|=\bigg|\frac{3\overline\xi_1\overline\eta_1}
{4\overline\xi^2_2}\bigg| =\bigg|\frac{3\xi_1\eta_1} {4\overline\xi^2_2}\bigg|\sim 1.
\nonumber
\end{align}
This concludes Lemma \ref{orthogonality}. Thus Proposition \ref{case213dprop1} follows.

Hence, in Case 2.1, from Proposition \ref{case211dprop1}, \eqref{maincase212} and
Proposition \ref{case213dprop1}, we have
 \begin{align}
 &\bigg|\iint_{\Omega_r^{*}}\frac{\sum^3_{{\iota}=1}m^2(\zeta_{\iota})(\xi_{\iota}+
\eta_{\iota})}
 {m(\zeta_1)m(\zeta_2)m(\zeta_3)}\widehat{v}_{1}(\lambda_1)\widehat{v}_{2}(\lambda_2)
 \widehat{v}_{3}(\lambda_3)\bigg|\notag \\
 \lesssim
 &N^{-\frac{1}{4}}(L_1L_2L_3)^{\frac{1}{2}}\|\widehat{v}_{1}\|_{L^{2}}
 \|\widehat{v}_{2}\|_{L^{2}}
 \|\widehat{v}_{3}\|_{L^{2}}.\label{maintricase21}
 \end{align}

Next, if $|\xi_2|\lesssim |\xi_1-\xi_3|$ and $|\xi_1|\sim|\xi_2|\sim|\xi_3|$,
 it is enough to consider $|\xi_{\iota}|\sim|\eta_{\iota}|\sim N_1$  for $\iota=1,2,3$.
Moreover, if $\xi_1\xi_3>0$, the strategy we use is almost the same as that employed in Case
2.1. \\{\bf Case 2.2} \hspace{2mm} $|\xi_2|\lesssim |\xi_1-\xi_3|$, $\xi_1\xi_3<0$ \ and \
 $|\xi_{\iota}|\sim|\eta_{\iota}|\sim N_1$  for $\iota=1,2,3$

Without loss of generality, we assume that
 $\xi_1>0,\xi_2>0$ and $\xi_3<0$.  If $\eta_1\eta_2\eta_3<0$ , then $|\xi_1\xi_2\xi_3+\eta_1\eta_2\eta_3|=|\xi_1\xi_2\xi_3|+|\eta_1\eta_2\eta_3|\sim N_1^3$,
 which contradicts to the low modulation assumption.

Hence $\eta_1\eta_2\eta_3>0$ and we can divide Case 2.2 into three subcases:
 \\{ Case 2.2.1} \hspace{2mm} $\xi_1>0 \ \ \xi_2>0 \ \ \xi_3<0$ \ and \ $\eta_1>0 \ \
 \eta_2<0 \ \ \eta_3<0$,
 \\{ Case 2.2.2} \hspace{2mm} $\xi_1>0 \ \ \xi_2>0 \ \ \xi_3<0$ \ and \ $\eta_1<0 \ \
 \eta_2>0 \ \ \eta_3<0$,
 \\{ Case 2.2.3} \hspace{2mm} $\xi_1>0 \ \ \xi_2>0 \ \ \xi_3<0$ \ and \ $\eta_1<0 \ \
 \eta_2<0 \ \ \eta_3>0$.

 Likewise we also find one of $|\xi_1\eta_2-\xi_2\eta_1|$
 and $|\xi_1\eta_2+\xi_2\eta_1+2(\xi_1\eta_1+\xi_2\eta_2)|$ is controllable ( see Table 2 ).
Case 2.2.1 and Case 2.2.2 could be dealt with in the same way as Case 2.1.1 and Case 2.1.3
respectively, while Case 2.2.3 could be dealt with in the same way as Case 2.1.2. We omit
the detail.
\begin{table}[]
\centering
 \caption{Subcases of 2.2 and the corresponding controllable quantities}
  \begin{tabular}{|l|c|c|c|c|c|c|c|c|c|}
  \hline
  \hspace{4mm}Case                                          & \multicolumn{3}{c|}{Case
2.2.1}
 & \multicolumn{3}{c|}{Case 2.2.2}
 & \multicolumn{3}{c|}{Case 2.2.3}
 \\ \hline
  $\hspace{7mm}\iota$                                       & \multicolumn{3}{c|}{1
 \hspace{9mm}   2   \hspace{8mm}    3}
 & \multicolumn{3}{c|}{1  \hspace{8mm}   2   \hspace{8mm}    3}
 & \multicolumn{3}{c|}{1  \hspace{8mm}   2   \hspace{8mm}    3}
 \\ \hline
  $\hspace{7mm}\xi_{\iota}$                                 & \hspace{2mm}$+$
\hspace{2mm}
 &\hspace{2mm} $+$     \hspace{2mm}                         &\hspace{2mm} $-$   \hspace{2mm}
 & \hspace{2mm}$+$    \hspace{2mm}                      &\hspace{2mm} $+$     \hspace{2mm}
 &\hspace{2mm} $-$   \hspace{2mm}                                 & \hspace{2mm}$+$
 \hspace{2mm}                      &\hspace{2mm} $+$     \hspace{2mm}
 &\hspace{2mm} $-$   \hspace{2mm}                      \\ \hline
  $\hspace{7mm}\eta_{\iota}$                                & \hspace{2mm}$+$
\hspace{2mm}
 &\hspace{2mm} $-$     \hspace{2mm}                         &\hspace{2mm} $-$   \hspace{2mm}
 & \hspace{2mm}$-$    \hspace{2mm}                      &\hspace{2mm} $+$     \hspace{2mm}
 &\hspace{2mm} $-$   \hspace{2mm}                                 & \hspace{2mm}$-$
 \hspace{2mm}                      &\hspace{2mm} $-$     \hspace{2mm}
 &\hspace{2mm} $+$   \hspace{2mm}                         \\ \hline
  $H_2$ or  $H_3$                                 &
 \multicolumn{3}{c|}{$H_3(\zeta_2,\zeta_3)\sim N_1^2$} &
 \multicolumn{3}{l|}{\hspace{3mm}$H_3(\zeta_1,\zeta_2)\sim N_1^2$} &
 \multicolumn{3}{l|}{$\hspace{3mm}H_2(\zeta_1,\zeta_2)\sim N_1^2$} \\ \hline
 \end{tabular}
  \end{table}
Finally, if $|\xi_2|\lesssim |\xi_1-\xi_3|$ and $|\xi_1|\gg|\xi_3|$ or $|\xi_1|\gg|\xi_3|$.
Without loss of generality, we assume that $|\xi_1|\gg|\xi_3|$. Then $|\xi_1|\sim|\xi_2|\sim
N_1\sim |\eta_3|$, thereby $H_3(\zeta_1,\zeta_3)=|\xi_1\eta_3-\xi_3\eta_1|\sim N_1^2$.
\\{\bf Case 2.3} \hspace{2mm} $|\xi_3|\ll|\xi_1|\sim|\xi_2|\sim N_1\sim |\eta_3|$

If $|\eta_1-\eta_2|\ll N_1$, then $|\eta_1|\sim |\eta_2|\sim |\eta_3|\sim N_1$, we have
$$|\xi_1\xi_2\xi_3+\eta_1\eta_2\eta_3|\sim |\eta_1\eta_2\eta_3| \sim N_1^3$$
which contradicts to $|\xi_1\xi_2\xi_3+\eta_1\eta_2\eta_3| \ll  \gamma_0N_1^3$.
Thus $|\eta_1-\eta_2|\sim N_1$.

We also observe that
 \begin{align}
 &|\xi_1\eta_3+\xi_3\eta_1+2(\xi_1\eta_1+\xi_3\eta_3)| \notag \\
 =&|\xi_3(\eta_1+2\eta_3)+\xi_1(\eta_3+2\eta_1)|\notag \\
 =&|\xi_3(\eta_1+2\eta_3)+\xi_1(\eta_1-\eta_2)|\notag \\
 \sim &|\xi_1(\eta_1-\eta_2)|\sim N_1^2. \nonumber
\end{align}
 That's to say, in this case,
$$|\xi_1\eta_3-\xi_3\eta_1|\cdot|\xi_1\eta_3+\xi_3\eta_1+2(\xi_1\eta_1+\xi_3\eta_3)|\sim
N_1^3.$$
Thus, for low frequency we use the nonlinear Loomis--Whitney inequality directly and for
high frequency we use bilinear Stichartz estimate. The proof is obvious, hence we omit the
detail.

In conclusion, no matter for high modulation and low modulation (including Case 1 and Case
2), we get
\begin{align}
&\bigg|\iint_{\Omega_r^{*}}\frac{\sum^3_{{\iota}=1}m^2(\zeta_{\iota})(\xi_{\iota}+
\eta_{\iota})}
{m(\zeta_1)m(\zeta_2)m(\zeta_3)}\widehat{v}_{1}(\lambda_1)\widehat{v}_{2}(\lambda_2)
\widehat{v}_{3}(\lambda_3)\bigg|\notag \\
\lesssim
&(\gamma_0^{-\frac{1}{2}}N^{-\frac{1}{2}}+N^{-\frac{1}{4}})(L_1L_2L_3)^{\frac{1}{2}}
\|\widehat{v}_{1}\|_{L^{2}}
\|\widehat{v}_{2}\|_{L^{2}}
\|\widehat{v}_{3}\|_{L^{2}}.\nonumber
\end{align}
Thus we have proved the trilinear estimate.

\section*{Acknowledgments}
 The first author M.S is partially supported by China Postdoctoral Science Foundation grant
2019M650872. 

\end{document}